\documentclass[12pt]{article}
\usepackage{amsfonts,amsmath}
\usepackage{amssymb,latexsym}
\usepackage[dvips]{epsfig}
\usepackage{amscd,xypic}
\xyoption{all}

\def\R{\mathbb R}

\def\C{\mathbb C} 
\def\Z{\mathbb Z}

\def\lra{\longrightarrow}
\def\Hom{\mathrm{Hom}}
\def\RHom{\mathrm{RHom}}

\newcommand{\mf}{\mathfrak} 
\newcommand{\mc}{\mathcal}
\renewcommand{\define}{\stackrel{\mbox{\scriptsize{def}}}{=}}

\newcommand\F{\mathcal F}
\newcommand\I{\mathcal I}
\renewcommand\L{\mathcal L}
\newcommand\PP{\mathbb P}
\newcommand{\rt}[1]{\stackrel{#1\,}{\rightarrow}}
\newcommand{\Rt}[1]{\stackrel{#1\,}{\longrightarrow}}
\newcommand\To{\longrightarrow}
\newcommand\into{\hookrightarrow}
\newcommand\comp{{}_{{}^\circ\!}}
\renewcommand\_{^{}_}
\newcommand\tr{\operatorname{tr}}
\newcommand\id{\operatorname{id}}
\renewcommand\hom{\mathcal Hom}
\newcommand\Ext{\operatorname{Ext}}
\newcommand\Bl{\operatorname{Bl}}
\newcommand\ev{\operatorname{ev}}
\newcommand{\oplusop}[1]{{\mathop{\oplus}\limits_{#1}}}
\def\drawing#1{\begin{center}\epsfig{file=#1}\end{center}}

\makeatletter \@addtoreset{equation}{section} \makeatother

\newtheorem{theorem}[equation]{Theorem}
\newtheorem{lemma}[equation]{Lemma}
\newtheorem{prop}[equation]{Proposition}
\newtheorem{defn}[equation]{Definition}
\newtheorem{cor}[equation]{Corollary}
\newenvironment{proof}{\noindent\emph{Proof}.}{\hfill $\square$\\}
\newenvironment{example}[1]{\noindent\textbf{#1.}}{\\}
\newenvironment{rmk}{\emph{Remark}:}{\\}
\newenvironment{rmks}{\emph{Remarks}:}{\\}

\title{\vspace{-1cm} Braid cobordisms, triangulated categories, and flag varieties}
\author{Mikhail Khovanov and Richard Thomas}
\date{\empty}

\begin{document}
\maketitle

%%%%%%%%%%%%%%%%%%%%%%%%%%%%%%%%%%%%%%%%%%%%%%%%%%%%%%%%%%%%%%%%%%%%%%%%%%%

\begin{abstract} We argue that various braid group 
actions on triangulated categories should be extended 
to projective actions of the category of braid cobordisms 
and illustrate how this works in examples. We 
also construct actions of both the affine braid group and the braid cobordism
category on  the derived category of coherent sheaves on the cotangent bundle
to the full flag variety.  
\end{abstract}
\thispagestyle{empty}

\renewcommand\contentsname{\vspace{-1cm}}
\tableofcontents

\section{Introduction}  

There are many known and conjectural braid group actions on 
triangulated categories, including 
\begin{itemize} 
\item Actions generated by a chain 
of spherical objects \cite{KS, ST, RZ, HK}. 
\item Action on the derived categories of (constructible) sheaves of vector spaces on flag varieties \cite{Ro1, Ro3}, and related actions on $D^b(\mc{O}_0),$ where $\mc{O}_0$ is a regular block of 
the highest weight category for $\mf{sl}_n,$ and on its subcategories \cite{St2}.
\item Actions on categories of complexes of matrix factorizations \cite{KR}.  
\item Actions on Fukaya-Floer categories of various symplectic manifolds, 
such as quiver varieties \cite{SS}.  
\item Actions on derived categories of coherent 
sheaves on quiver varieties and other Calabi-Yau manifolds. 
\end{itemize} 

We begin by giving a brief survey of these and related actions in Sections~\ref{sec-survey} and \ref{sec-bc}. This is meant to convey the diverse areas in which braid
group actions arise, but it is not necessary to understand or even read these
sections to follow the rest of the paper.
By a \emph{weak action} of a group $G$ on a category $\mc{C}$ 
we mean an assignment of an invertible functor $F_g: \mc{C} \lra \mc{C}$ to 
each $g\in G$ such that $F_ g \circ F_h \cong F_{gh}$ for all $g,h\in G.$ It follows 
that $F_1\cong \mathrm{Id},$ where $1\in G$ is the unit element. 
A weak action is upgraded to a (genuine) action if the functor isomorphisms
$F_ g \circ F_h \cong F_{gh}$ are 
selected to satisfy the associativity constraint, which says that all 
diagrams below are commutative 
\begin{equation} 
\begin{CD} 
  F_f \circ F_g \circ F_h    @>{\cong}>>  F_{fg} \circ F_h   \\
   @VV{\cong}V   @VV{\cong}V   \\
   F_f \circ F_{gh}  @>{\cong}>>  F_{fgh} 
\end{CD}  
\end{equation} 
Deligne \cite{De} was the first to discuss braid group action on categories,
implicitly emphasizing the 
difference between weak and genuine actions, and  gave a neat criterion 
for when a weak action of a braid group lifts to a genuine 
action.  
Some of the braid group actions reviewed below are genuine and the
others are believed to be genuine. 

An action of a  group $G$ by exact functors on a triangulated category $\mc{C}$
induces an action of $G$ on the Grothendieck group $K(\mc{C}).$ The exact
functor $F_g$ descends to an endomorphism, denoted $[F_g],$ of 
the abelian group $K(\mc{C}).$  Isomorphisms 
$F_g\circ F_h \cong F_{gh}$ become equalities between maps 
$[F_g][F_h] = [F_{gh}]$ and we obtain an action of $G$ on  $K(\mc{C}).$
Often $K(\mc{C})$ is a free abelian group and, after tensoring it 
and all maps $[F_g]$ with $\C,$ we get a complex representation of $G,$ 
an action of $G$ on the complex vector space $K(\mc{C})\otimes \C.$ 
The operator in $K(\mc{C})\otimes \C$ associated with $g$ is still 
denoted $[F_g].$ In some cases (for instance, when $\mc{C}$ is a triangulated 
category of complexes of graded modules over a graded ring), $K(\mc{C})$ 
is naturally a $\Z[q,q^{-1}]$-module, and $[F_g]$ are module maps. 
Then, to get a complex representation, 
 we tensor $K(\mc{C})$ with $\C$ over $\Z[q,q^{-1}],$ selecting 
$q\in \C$ to be a generic complex number. Alternatively, one can work 
with the $\Z[q,q^{-1}]$-module $K(\mc{C})$ itself. 

In this paper we specialize to the case when $G$ is the $n$-stranded braid
group $\mathrm{Br}_n.$ Quite a few well-known complex representations of the braid group can be lifted to actions of $\mathrm{Br}_n$ on triangulated categories, including the Burau representation \cite{KS, ST}, representations that factor through various representations of the Hecke algebra of the symmetric group (by restricting the braid group action on the homotopy category of a regular block of $\mc{O}$ for $\mathfrak{sl}_n$ to the subcategories considered in \cite{KMS}), their $q=1$ specializations, etc. 

The existence of such actions of $\mathrm{Br}_n$ on triangulated 
categories is nontrivial, and their potential remains underexplored. 
The synthesis of the braid group, a manifestly topological object, 
and triangulated categories, which are strongly rooted in algebra, 
looks unusual and perplexing. What can one do with this symbiosis? 

\medskip

Suppose we are given a braid group action on $\mc{C}.$ To 
each braid $g$ there is associated a functor $F_g.$ Given two 
elements $g,h$ of the braid group, we can form the set of 
natural transformations $\mathrm{Hom}(F_g,F_h)$ 
from $F_g$ to $F_h$ in $\mc{C}.$ Together, they form a monoidal
category with elements of $\mathrm{Br}_n$ as objects 
and $\mathrm{Hom}(F_g,F_h)$  as morphisms from $g$ to $h.$ 
This information about the braid group action on $\mc{C}$ 
is lost upon passing to the Grothendieck group. 

At the same time, there exists a monoidal category of purely 
topological origin with elements of the braid group 
$\mathrm{Br}_n$ as objects. We call it the category  
of \emph{braid cobordisms} and denote it by $\mathcal{BC}_n.$ 

We can think of an $n$-stranded braid as a cobordism 
from $n$ fixed points on a plane to itself, such that the projection onto 
the vertical direction has no critical points, see below. Isotopy 
classes of $n$-stranded braids constitute the braid group $\mathrm{Br}_n.$ 
Sometimes, isotopy classes are themselves called braids. 

\drawing{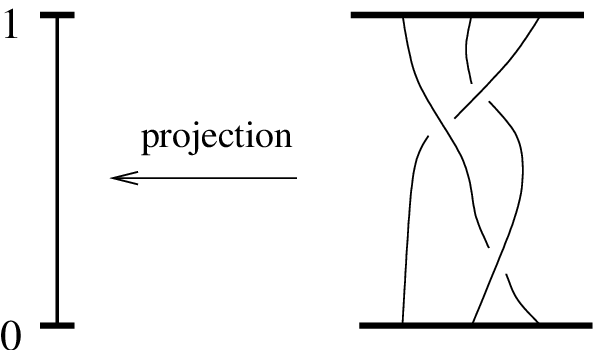} 

\medskip

Thus, a braid is a embedding of a 1-manifold which is a disjoint union 
of closed intervals into 
$\R^2\times [0,1]$ so that the projection onto $[0,1]$ has no 
critical points, and the boundary of the 1-manifold is mapped to 
$2n$ points 
$$(1,0,0), (2,0,0), \dots, (n,0,0), (1,0,1), (2,0,1), \dots, (n,0,1)\} \in 
 \R^2\times [0,1].$$ 
The first $n$ points on this list belong to the bottom $\R^2,$ i.e. 
to $\R^2\times \{0\},$ the last $n$ points to $\R^2\times \{1\}.$ 

\begin{figure} \drawing{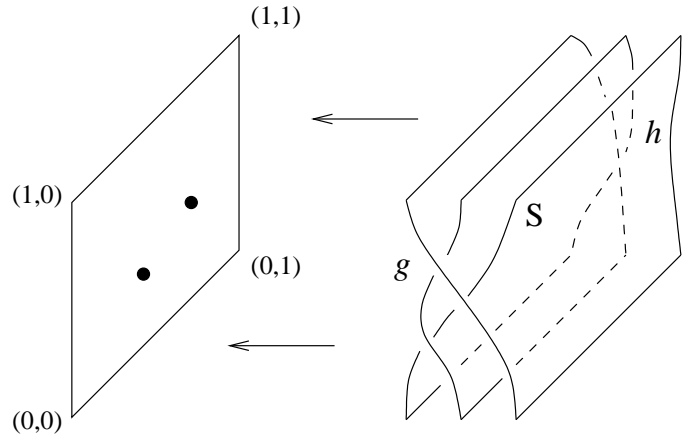} \caption{A depiction of a braid 
cobordism (mostly of its boundary), and its projection onto the unit square. 
Dots in the square are the simple branch points of the projection. The image of 
$g$ is the edge $[0,1]\times \{0\}$ of the square. The image of $h$ is 
the edge $[0,1]\times \{1\}$ of the square. Corners of $S$ are mapped 
to the vertices of the square.} \label{picasso1} 
\end{figure}

A braid cobordism between braids $g,h\in \mathrm{Br}_n$ is 
 a compact surface $S$ with boundary and corners, 
smoothly and properly embedded in $\R^2\times [0,1]^2,$ such 
that 
\begin{itemize} 
\item the boundary of $S$ is the union of four 1-manifolds  
\begin{eqnarray*}
 S \cap (\R^2\times [0,1]\times \{0\}) & = & g, \\
 S\cap (\R^2\times [0,1]\times \{1\}), & = & h, \\
S\cap (\R^2\times\{0\} \times [0,1]) & =& \{1,2, \dots, n\} \times \{0\} \times [0,1], \\
 S\cap (\R^2\times\{1\} \times [0,1]) & =& \{1,2, \dots, n\} \times \{1\} \times [0,1]. 
 \end{eqnarray*} 
 \item The projection of $S$ onto $[0,1]^2$ is a branched covering with 
simple branch points only.
\end{itemize} 
A braid cobordism is schematically depicted in figure~\ref{picasso1}. 

\begin{figure} \drawing{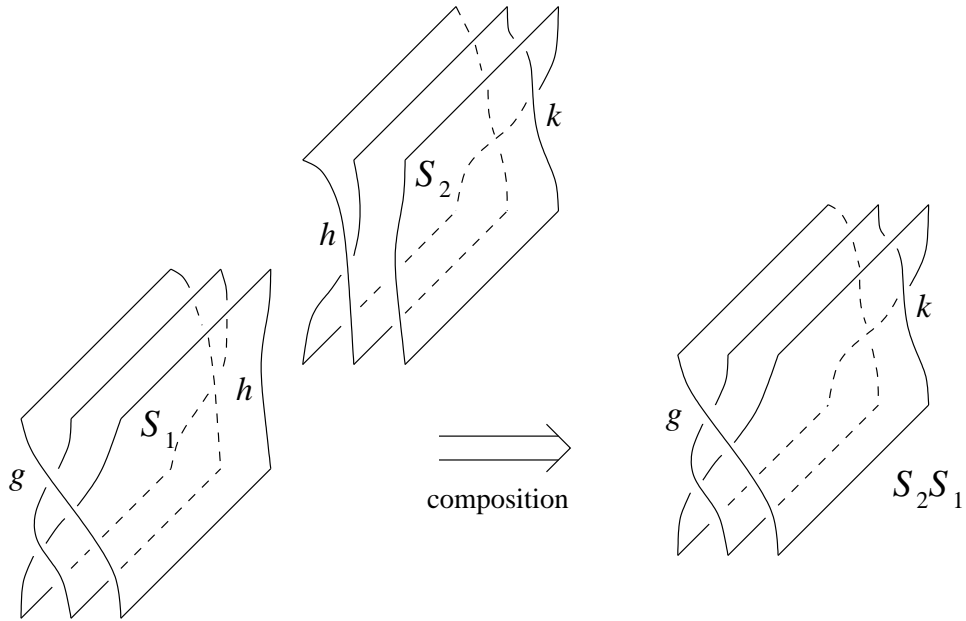} \label{picasso2} 
\caption{Composition of braid cobordisms} \end{figure}

Braid cobordisms which are isotopic rel boundary via an isotopy through braid
cobordisms are called equivalent. 
Let $\mc{BC}_n$ be the category with $n$-stranded braids 
as objects and isotopy classes (rel boundary) of braid cobordisms as morphisms.
Composition of a morphism $S_1$ from $g$ to $h$ with a morphism 
$S_2$ from $h$ to $k$ is the braid cobordism $S_2 S_1$ given by 
concatenating $S_2$ and $S_1$ along their common boundary $h,$ see 
figure~\ref{picasso2}. 

Note that in $\mc{BC}_n$ two isotopic braids are isomorphic, but not equal.
The category $\mc{BC}_n$ is monoidal, since given a morphism $S_1$ from 
$g_1$ to $h_1$ and a morphism $S_2$ from $g_2$ to $h_2,$ we can concatenate
them to a morphism 
$S_2 \circ S_1$ from $g_2g_1$ to $h_2 h_1.$ The concatenation is done 
by gluing the top portion of $\partial (S_1)$ to the bottom portion of 
$\partial(S_2).$  The monoidal structure is 
not strictly associative, since braids $g_3(g_2g_1)$ and $(g_3g_2)g_1$ 
are isotopic but not equal. 

In addition, the operation of placing braid cobordisms in parallel 
gives us a bifunctor 
\begin{equation} \label{eq-parallel}  
 \mc{BC}_n \times \mc{BC}_m \lra \mc{BC}_{n+m}.  
\end{equation} 

Braid cobordisms are called \emph{simple braided surfaces} in \cite[Section 16]{Ka}, and originally appeared in Rudolph \cite{Ru}.   
Braid cobordisms from the trivial braid to itself, 
with corners smoothed out, are called \emph{simple surface braids} in the 
literature \cite{Ka}, \cite[Section 3]{CS2},  and were introduced by O.~Viro. 

In Section~\ref{sec-bc1} we recall a combinatorial 
version of $\mc{BC}_n,$ following \cite{CS1}, \cite[Section 3]{CS2}. 
A braid cobordism can be described diagrammatically, via 
Carter and Saito's \emph{braid movies}. This leads to a 
category, denoted $\mc{BC}^c_n,$ with braid words as objects and braid movies,
modulo \emph{braid movie moves,} as morphisms. We call $\mc{BC}^c_n$ the
category of \emph{combinatorial braid cobordisms}. The categories $\mc{BC}_n$
and $\mc{BC}^c_n$ are equivalent. 

In the beginning of Section~\ref{sec-survey} we define rings $A_n,$ 
suitable triangulated categories $\mc{C}(A_n)$ of complexes of graded $A_n$-modules,  
and recall how to construct a weak braid group action on $\mc{C}(A_n).$ 
With every braid word $g$ we associate a functor $F_g$ such 
that $F_g\cong F_h$ if $g,h$ define the same braid in 
$\mathrm{Br}_{n+1}.$ 
Let $\mathrm{Fun}(A_n)$ be the category with 
objects--exact endofunctors in $\mc{C}(A_n)$ and morphisms--natural 
transformations of functors.  

\begin{theorem} \label{first-theorem} 
The braid group action on $\mc{C}(A_n)$ 
extends to a projective mon\-oidal functor from the category 
$\mc{BC}^c_{n+1}$ of combinatorial braid cobordisms to 
the category $\mathrm{Fun}(A_n).$ 
\end{theorem} 

To every braid cobordism $S$ between braid words $g$ and $h$ 
we assign a natural transformation $F_S: F_g \lra F_h,$ 
well-defined up to overall minus sign. The presence of the 
minus sign explains why we use the word \emph{projective} 
in the statement of the theorem. There are equalities 
of natural transformations 
$$F_{S_2}  F_{S_1} =\pm F_{S_2  S_1}$$ 
for any pair of braid cobordisms $S_1$ from $g$ to $h$ and 
$S_2$ from $h$ to $k.$ 

The functor $F$ intertwines the monoidal structures of the 
two categories: $F_{gh} \cong F_g \circ F_h$ and 
$F_{S_2} \circ F_{S_1} = \pm F_{S_2 \circ S_1}$ for 
any cobordisms $S_1, S_2$ between $g_1,h_1$ and $g_2,h_2,$ 
respectively. 

Composing our functor with the equivalence $\mc{BC}_{n+1}\cong \mc{BC}^c_{n+1},$ 
we obtain a projective monoidal functor from the topological category 
$\mc{BC}_{n+1}$ of braid cobordisms to the algebraic category 
$\mathrm{Fun}(A_n)$ of natural transformations between 
exact functors in $\mc{C}(A_n).$ 

\medskip

This example serves as the first illustration of the principle: 

\medskip

\emph{Interesting braid group actions on triangulated categories extend to projective 
representations of the category of braid cobordisms.}

\medskip

In Sections~\ref{sec-survey} and \ref{sec-bc} 
we list many additional examples when a braid group action on a 
triangulated category can be nontrivially extended to braid cobordisms.  
This multitude of examples motivates us to elevate the above 
principle to the definition of categorification of 
a braid group action.  From here on by a categorification of a braid 
group representation we mean a projective functor from $\mc{BC}_n$ to 
the category of exact endofunctors in a triangulated category, 
as outlined above and defined in Definition~\ref{d-one} of 
Section~\ref{sec-bc1} below. 

\medskip

Our requirement that the representation $F$ of braid cobordisms be only 
a projective functor is similar to considering weak braid group actions only, 
while requiring that $F$ be a functor should be analogous to considering  
genuine braid group actions. We expect that most, if not all, 
of the projective functors $\mc{BC}_n\stackrel{F}{\lra}\mathrm{Fun}(\mc{C})$
can be turned into genuine functors.  The diagram below summarizes 
how the authors currently think about braid group actions on 
triangulated categories, each arrow denoting a structural upgrade.  

\drawing{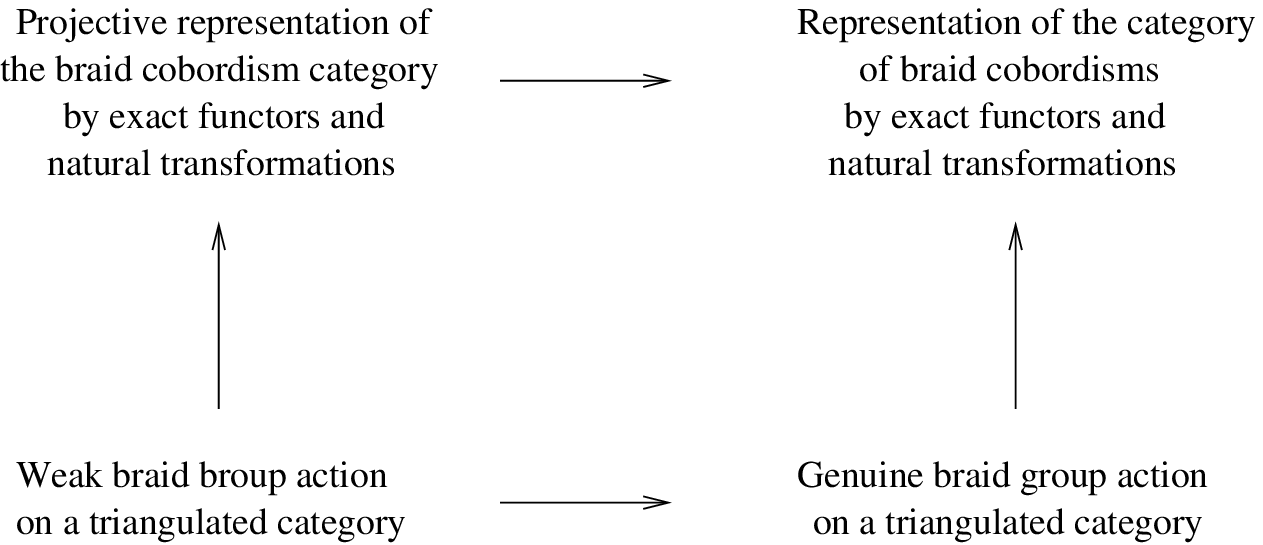}

At the end of Section~\ref{sec-bc-RR} we sketch the relation 
between Rouquier's definition of a categorification of a braid group 
action and ours. 
\medskip

In the last two sections we turn to coherent sheaves, where the proofs become
more technical but the ideas are essentially the same.
In Section \ref{flagsection} we exhibit an action of Br$_n$ on
$D^b(T^*Fl)$, the derived category of coherent sheaves on the cotangent bundle to the variety of full flags in $\C^n$.
The functors we use are family versions of the Dehn twists about $-2$-curves
of \cite{KS, ST}, as considered in \cite{Ho, Sz1}. Locally, away from some closed subsets of $T^*Fl$, the subschemes we twist about form a family
$A_{n-1}$-chain in the sense of \cite{Sz2},
but, since they all contain the zero section $Fl\subset T^*Fl$, the proof of the braid relations is much more complicated than in \cite{ST, Sz2} (due
to the extra term $U_{ij}$ in Proposition \ref{H}, for instance).

Given an $A_n$-chain of $-2$-curves in a surface, it is fairly easy to
show that the structure sheaf of the whole chain (twisted by an appropriate
line bundle) gives a spherical object whose Dehn twist extends the braid
group action of \cite{ST} to an \emph{affine} braid group action \cite[Example 3.9]{ST}. While such a simple description does not generalise to $T^*Fl$,
one can show that the extra twist is given by a certain product of the original
generators conjugated by tensoring by a line bundle. The resulting formula
does generalise (\ref{Tdef}), so in Section \ref{affine} we
extend our action to one of the affine braid group on $D^b(T^*Fl)$.

In Section \ref{surfaces} we extend the braid group action \cite{ST} on the derived category of coherent sheaves of certain surfaces to a projective action of the braid cobordism category.

\begin{theorem} Fix a smooth quasiprojective complex surface $X$ with an $A_n$-chain of $-2$-curves. Then the weak braid group action on $D^b(X)$
of \emph{\cite{ST}} extends to a projective action of $\mc{BC}^c_{n+1}$.
\end{theorem}

Finally in Section \ref{flagcobordism} we show how to do the same for
$D^b(T^*Fl)$. We expect these two generalisations of the braid group action on $D^b(T^*Fl)$ to be compatible, giving a genuine action of the affine braid cobordism category.

\begin{theorem} There is a projective action of $\mc{BC}^c_n$ on the
derived category of coherent sheaves on the cotangent bundle of the full
flag variety.
\end{theorem}

When the second author spoke about the braid group action on $D^b(T^*Fl)$
in Leeds in November 2005, Rapha\"el Rouquier informed us of another way to produce a braid
group action (which is hard to make explicit and so difficult to confirm
is the same as ours). Start with the braid group action on the derived category of constructible sheaves on the flag variety (see for example \cite{Ro1}).
This is induced by kernels of constructible sheaves on the product. Their
associated mixed Hodge modules (see for instance page 18 of \cite{Ro2})
interpolate between (derived categories of) constructible sheaves
and $\C^*$-equivariant coherent sheaves on the cotangent bundle. This gives
$\C^*$-equivariant kernels generating a braid group action on the derived category of coherent sheaves on the cotangent bundle of the flag variety.

Recently Bezrukavnikov, Mirkovi\'c and Rumynin \cite{BMR} have also produced an affine braid group action (for any group, not just of type $A_n$) by very different and sophisticated finite characteristic methods which we cannot claim to understand. Assuming that their action, which also seems to be hard
to make explicit, coincides with the explicit action of Section \ref{flagsection} then we show it extends to an action of the braid cobordism category.

\medskip
\noindent \textbf{Acknowledgements.} 
M.\,K. would like to thank Tom Mrowka 
and Peter Ozsv\'ath for the opportunity to lecture on some of 
the material in this paper at the IAS/Park City summer school 2006 
and Scott Carter for enlightening email discussions.
R.\,P.\,T. would like to thank Martin Bridson, Tom Bridgeland, Daniel Huybrechts, Ivan Smith, Simon Willerton and most especially Paul Seidel
for valuable discussions.
M.\,K. gratefully acknowledges NSF partial support via grant DMS-0407784.
R.\,P.\,T. is partially supported by the Royal Society and the Leverhulme Trust.

%%%%%%%%%%%%%%%%%%%%%%%%%%
%%
%%  BRIEF SURVEY OF BRAID GROUP ACTIONS 
%%
%%%%%%%%%%%%%%%%%%%%%%%%%%%%

\section{Examples of actions of braid groups and braid 
cobordisms on triangulated categories} \label{sec-survey}

\begin{example}{1. Spherical objects}
Let $A_n,$ for $n>2,$ be the quotient of the path ring of the quiver 

\drawing{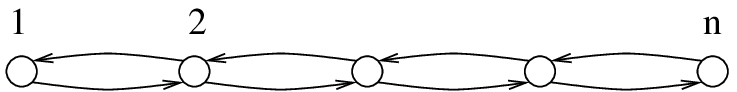} 

\noindent by the relations 
\begin{eqnarray*} 
 (i|i+1|i+2) & = & 0, \hspace{0.2in} 1\le i \le n-2,  \\
 (i|i-1|i-2) & = & 0, \hspace{0.2in} 3\le i \le n, \\
 (i|i-1|i) & = & (i|i+1|i), \hspace{0.2in} 2\le i \le n-1.
\end{eqnarray*} 
We define $A_2$ as the quotient of the above quiver ring, 
for $n=2,$ by the relations $(1|2|1|2)=0=(2|1|2|1).$ 
Let $A_1 = \Z[X_1]/(X_1^2).$ 

The ring $A_n$ is a free abelian group of rank $4n-2$ and 
has a basis given by $(i), 1\le i\le n,$ which are 
length zero paths and minimal idempotents,  
all possible paths of length one (which are $(i|i\pm 1)$), and 
$$X_i \define (i|i-1|i), \hspace{0.2in} i\ge 2, \hspace{0.2in} 
 X_1 = (1|2|1).$$  
Denote this basis by $B(A_n).$ 
The center of $A_n$ is a free abelian group with the basis 
$$ 1=\sum_{1\le i\le n} (i), \hspace{0.1in} X_1, \hspace{0.1in} 
X_2, \hspace{0.1in} \dots \hspace{0.1in} X_n.$$ 
The ring $A_n$ is naturally graded by path length. Let 
$$P_i \define A_n (i), \hspace{0.2in} _iP\define (i)A_n. $$ 
These are indecomposable left, respectively right, projective 
$A_n$-modules. Consider the complexes of $A_n$-bimodules
\begin{eqnarray*} 
 R_i : &\quad&  0 \lra P_i \otimes \hspace{0.03in} _iP 
 \stackrel{\beta_i}{\lra} A_n \lra 0, \\
R_i':  &&  0 \lra A_n \stackrel{\gamma_i}{\lra}
P_i \otimes \hspace{0.03in} _iP \lra 0,
\end{eqnarray*} 
with $A_n$ placed, in both complexes, in cohomological degree $0.$ 
Here $\beta_i$ is the composition of the inclusion and multiplication 
$$ P_i \otimes \hspace{0.03in} _iP  \subset A_n\otimes A_n \lra A_n,$$ 
while $\gamma_i$ is determined by 
$$ \gamma_i(1) = \sum \alpha_1 \otimes \alpha_2 ,$$ 
the sum being over all $\alpha_1, \alpha_2\in B(A_n)$ such that 
$\alpha_2 \alpha_1 = X_i.$ Thus, if $1< i < n,$
\begin{eqnarray*} 
\gamma_i(1)  & = & (i-1|i)\otimes(i|i-1) + (i+1|i)\otimes (i|i+1)  \\
      & & + (i)\otimes (i|i+1|i) + (i|i+1|i) \otimes (i). 
\end{eqnarray*} 
The projectives $P_i, \hspace{0.03in} _iP$ are examples of 
\emph{spherical} objects, so called since homologically they 
resemble spheres, 
$$\mathrm{End}(P_i) \cong \Z[X_i]/(X_i^2)\cong 
\mathrm{H}^{\ast}(\mathbb{S}^k,\Z)$$
(where $k=2$ if we use the path length grading on $A_n$), 
and satisfy a duality property, which in this example manifests itself in $A_n$ being a symmetric ring. 

Consider the category of complexes of $A_n$-bimodules up to 
chain homotopies (the \emph{homotopy category}, for short). 
Recall from \cite{KS} that the following isomorphisms hold in 
this category: 
\begin{eqnarray} 
   R_i \otimes R_i'  & \cong & A_n\ \cong\ R_i' \otimes R_i ,  \label{eq-1} \\
   R_i \otimes R_{i+1} \otimes R_i & \cong & 
   R_{i+1} \otimes R_i \otimes R_{i+1}, \label{eq-2}\\
   R_i \otimes R_j & \cong & R_j \otimes R_i, \hspace{0.2in} |i-j|>1. \label{eq-3}
\end{eqnarray} 
The tensor products above are taken over $A_n.$ Let $\mc{C}(A_n)$ 
be one of the triangulated categories formed out of complexes 
of left $A_n$-modules, such as the derived category $D^b(A_n\mathrm{-mod}),$
the homotopy category of complexes of $A_n$-modules, etc.  
    Tensoring $M \in \mc{C}(A_n)$ with the complexes 
$R_i, R_i'$ produces complexes of $A_n$-modules 
$$ R_i\otimes_{A_n} M, \hspace{0.1in} R_i'\otimes_{A_n} M\in \mc{C}(A_n).$$ 
We denote the resulting endofunctors in $\mc{C}(A_n)$ by $\mc{R}_i, \mc{R}_i':$
$$\mc{R}_i(M) \define R_i\otimes_{A_n} M, \hspace{0.2in}
 \mc{R}_i'(M) \define R_i'\otimes_{A_n} M.$$
The earlier bimodule isomorphisms translate into functor isomorphisms 
 \begin{eqnarray} 
   \mc{R}_i \circ \mc{R}_i'  & \cong & \mathrm{Id} \cong \mc{R}_i' \circ \mc{R}_i ,  
 \label{brrel-1} \\
   \mc{R}_i \circ \mc{R}_{i+1} \circ \mc{R}_i & \cong & 
   \mc{R}_{i+1} \circ \mc{R}_i \circ \mc{R}_{i+1}, \label{brrel-2} \\
   \mc{R}_i \circ \mc{R}_j & \cong & \mc{R}_j \circ \mc{R}_i, \hspace{0.2in} |i-j|>1,
   \label{brrel-3} 
\end{eqnarray} 
implying that the $(n+1)$-stranded braid group $\mathrm{Br}_{n+1}$ acts weakly
on $\mc{C}(A_n).$ It is proved in \cite{KS} that this action 
is faithful (faithfulness in the case $n=2$ was shown in \cite{RZ}). 
Theorem~\ref{first-theorem} states that this action extends to braid cobordisms. \medskip

Braid cobordisms are a subset of  tangle cobordisms \cite{BL}. 
Tangle cobordisms constitute a 2-category $\mathbb{T},$ with nonnegative 
integers as objects, tangles as morphisms, and tangle cobordisms as 2-morphisms. 
There are several known projective 2-functors from $\mathbb{T}$ to the 
2-category of triangulated categories. In the latter triangulated categories are objects, 
exact functors are 1-morphisms and natural transformations are 2-morphisms. Any 
such 2-functor restricts to a projective functor from the category of 
$n$-stranded braid cobordisms to the category of functors and natural transformations
in some triangulated category. Examples {\bf 2-6} below are of this nature. \end{example}

\begin{example}{2. The rings $H^n$} 
A certain family of rings $H^n$ was introduced  
in \cite{Kh1} and used in \cite{Kh1, Kh2} to define a projective 2-functor 
from a suitable version of the category $\mathbb{T}$ of tangle cobordisms 
to the 2-category of exact functors (and natural transformations) between
the categories $\mc{C}(H^n),$ over various $n\ge 0.$ Here $\mc{C}(H^n)$ 
is defined similarly to $\mc{C}(A_n).$ For instance, $\mc{C}(H^n)$ could be the derived category $D^b(H^n\mathrm{-mod}).$ 

In this construction, there is a sign indeterminacy in the natural transformation
associated with a tangle cobordism. 

Restricting from tangles to braids and from tangle cobordisms to 
braid cobordisms yields a projective functor from $\mc{BC}_{2n}$ to 
$\mathrm{Fun}(H^n),$ the category of exact functors in $\mc{C}(H^n)$ 
and natural transformations between them. To a $2n$-stranded braid 
we assign a complex of $H^n$-bimodules; tensoring with this 
complex is an endofunctor in $\mc{C}(H^n).$ To a braid cobordism 
there is assigned a homomorphism between complexes of $H^n$-bimodules, 
inducing a natural transformation between the corresponding functors. 

If we take $\mc{C}(H^n)$ 
to be the homotopy category of complexes of finitely-generated graded projective
$H^n$-modules (the rings $H^n$ are graded), then its Grothendieck 
group, after tensoring with $\C$ over $\Z[q,q^{-1}],$ become isomorphic 
to the space of $U_q(\mf{sl}_2)$-invariants in $V^{\otimes 2n},$ the 
$2n$-th tensor power of the defining representation of the quantum 
group $U_q(\mf{sl}_2).$ Thus, we obtain a categorification, 
in the sense of Definition~\ref{d-one}, of the braid group action  
on $\mathrm{Inv}(V^{\otimes 2n}).$
\end{example}

\begin{example}{3} Bar-Natan's invariant of tangle cobordisms \cite{BN} can be 
restricted to braid cobordisms and viewed in the framework of 
triangulated categories which are similar to but more general than $\mc{C}(H^n).$
His construction also restricts to a categorification of the braid group action on $\mathrm{Inv}(V^{\otimes 2n}).$ 
\end{example}

\begin{example}{4} An invariant of tangle cobordisms via highest 
weight categories was constructed by Stroppel \cite{St1, St2}, following 
conjectures in \cite{BFK}. Restricting her invariant to braids and braid cobordisms 
results in a projective functor from $\mc{BC}_n$ to 
$\mathrm{Fun}(D^b(\mc{O}^n)),$ the category of exact 
endofunctors in the derived category of $\mc{O}^n,$ the 
latter a direct sum of certain parabolic subcategories of a regular 
block of the highest weight category for $\mf{sl}_n.$ 
In this example the projectivity (indeterminacy) takes 
values in $\C^{\ast},$ due to $\C$ being the ground field. 
The Grothendieck group of (the graded version of) $\mc{O}^n$ 
is isomorphic to $V^{\otimes n},$ and the projective functor 
$$ \mc{BC}_n \lra \mathrm{Fun}(D^b(\mc{O}^n))$$
is a categorification of the braid group action on $V^{\otimes n}.$ 
\end{example}

\begin{example}{5} A categorification \cite{Kh4} of the quantum $\mf{sl}_3$
link  invariant was extended to tangles and tangle cobordisms by 
MacKaay and Vaz \cite{MV}. After restricting their structure to braid cobordisms
one should obtain a categorification of the action of the 
$3n$-stranded braid group on $\mathrm{Inv}(V^{\otimes 3n}),$ 
the $3n$-th tensor power of the fundamental representation $V$ 
of quantum $\mf{sl}_3.$ More generally, their construction
conjecturally categorifies spaces of invariants in arbitrary tensor products
of $V$ and its dual $V^{\ast}.$ The braid group $\mathrm{Br}_n$ 
acts on $V^{\otimes n},$ and hence on 
$\mathrm{Inv}(V^{\otimes n}\otimes W),$
 where $W$ is any tensor product of $V$'s and $V^{\ast}$'s. 
Following \cite{MV}, one should also obtain a categorification 
of the braid group action on $\mathrm{Inv}(V^{\otimes n}\otimes W)$, by 
adding dummy strands. In this example, the projectivity is only in the minus sign.
\end{example}

\begin{example}{6} Yet another invariant of tangle cobordisms was introduced
in \cite{KR}, based on matrix factorizations. A matrix factorization of $f\in
R,$ where 
$R$ is a polynomial algebra, consists of a pair of free $R$-modules 
$M^0, M^1$ and $R$-module maps 
$$ M^0 \stackrel{d}{\lra} M^1 \stackrel{d}{\lra} M^0 $$
such that $d^2(m) =fm$ for any $m\in M^0,M^1.$ 
The category of matrix factorizations up to chain homotopy is 
triangulated. Given a braid word $g,$ the construction of \cite{KR} 
associates to it a complex of matrix factorizations in 
$2n$ variables $x_1, \dots, x_n, y_1,\dots,y_n$ with 
$$f = \sum_{i=1}^n x_i^{m+1}-y_i^{m+1} ,$$ 
where $n$ is the number of strands. The isomorphism class 
of this complex in the category of complexes of matrix factorizations 
up to chain homotopy is an invariant of the braid. A braid cobordism 
between braids $g,h$ induces a homomorphism of corresponding complexes, 
well-defined up to rescaling by nonzero rational numbers.  One can
easily restate this construction in the language of functors, resulting 
in yet another categorification of a braid group action subject 
to the constraints of Definition~\ref{d-one}. We have no clue 
what representation of the braid group this categorification 
categorifies. It should contain a subrepresentation isomorphic 
to $\mathrm{Inv}(V^{\otimes n}\otimes (V^{\ast})^{\otimes n}),$ 
where $V$ is the defining representation of quantum $\mf{sl}_m$
and the braid group acts on the first $n$ factors.  
\end{example}

\begin{example}{7. Highest weight categories and flag varieties} 
We start with the triangular decomposition 
$\mf{sl}_n= \mf{n}_+ \oplus \mf{h} \oplus\mf{n}_-$ of the Lie 
algebra $\mf{sl}_n$, where $\mf{n}_+,$ respectively $\mf{n}_-,$ 
are the Lie algebras of strictly upper-triangular, respectively 
lower-triangular, matrices, and $\mf{h}$ the Lie algebra of traceless 
diagonal matrices. The category $\mc{O}$ is a full subcategory of 
the category of finitely-generated $\mf{sl}_n$-modules. 
A finitely-generated $\mf{sl}_n$-module $M$ 
belongs to the highest weight category $\mc{O}$ if and only if $\mf{h}$ 
acts semisimply and $\mf{n}_+$ locally-nilpotently on $M,$ see \cite{BG}.  
The trivial one-dimensional $\mf{sl}_n$ 
representation defines a central character $\chi: Z\lra \C,$ a
homomorphism from the center of $U(\mf{sl}_n)$ to $\C.$ Define 
$\mc{O}_0,$ the trivial central character block of $\mc{O},$ 
as the full subcategory of $\mc{O}$ consisting of modules annihilated 
by some power of $\mathrm{ker}(\chi)\subset Z.$ 
This category is equivalent to the category of finite-dimensional 
modules over some finite-dimensional Koszul $\C$-algebra $B_n,$ 
although a simple description of $B_n$ is not known even for $n=5.$  

The category $\mc{O}_0$ has $n!$ isomorphism classes of simple 
modules. There are exact self-adjoint endofunctors $\theta_i$ in 
$\mc{O}_0$ (called \emph{translation across the $i$-th wall}, $1\le i \le n-1$) 
which come with natural transformations 
$\theta_i \lra \mathrm{Id}$ and $\mathrm{Id}\lra \theta_i$. 
Let $\mc{C}(\mc{O}_0)$ be either the derived category of $\mc{O}_0$ 
or the category of complexes of objects in $\mc{O}_0$ up to chain 
homotopies. Taking an object $M\in \mc{C}(\mc{O}_0)$ to 
the total complex of 
$$ 0 \lra \theta_i M \lra M \lra 0$$ 
is an invertible functor $\mc{R}_i$ in $\mc{C}(\mc{O}_0),$ with 
the inverse functor $\mc{R}_i'$ taking $M$ to the total complex of 
 $$ 0 \lra M \lra \theta_i M \lra 0.$$ 
It follows from \cite{MS1}, \cite[Proposition 10.1]{MS2} that the
functors $\mc{R}_i$ satisfy the braid relations (\ref{brrel-1})-(\ref{brrel-3})
and define a weak braid group action on $\mc{C}(\mc{O}_0)$.

\medskip

A similar braid group action in the Koszul dual case of Zuckerman functors 
(rather than translation functors) acting on $D^b(\mc{O}_0)$
was given a detailed treatment  by Rouquier \cite{Ro1, Ro3}, 
who showed, in particular, that the action is genuine.  
The derived category $D^b(\mc{O}_0)$ embeds into the derived category $D^b(Y)$  
of sheaves of vector spaces on the variety $Y$ of full flags in $\mathbb{C}^n.$ For each 
$i=1, \dots, n-1$ there exists a correspondence $Y_i \subset Y \times Y$ 
which consists of pairs of flags that coincide except at the $i$-th term. The convolution 
with the constant sheaf on $Y_i$ extended by $0$ to the entire $Y\times Y$ is 
an invertible functor $\mc{R}_i$ in $D^b(Y).$ These functors generate a genuine 
action of the braid group on $D^b(Y),$ see \cite{Ro1}. 

This construction admits a variation where, for each decomposition $\mu=(\mu_1,
\dots, \mu_k)$ of $n,$ we consider the variety $Y_{\mu}$ of partial flags
in $\C^n,$ of 
dimensions $\mu_1, \mu_1+\mu_2, \dots, n - \mu_k$, 
and construct an invertible functor from $D^b(Y_{\mu})$ to $D^b(Y_{s_i \mu}),$
where $s_i\mu$ is the decomposition given by transposing $\mu_i$ and $\mu_{i+1}.$
The functor is given by convolution with $Y_i \subset Y_{\mu}\times Y_{s_i\mu},$
where the latter  consists of pairs of flags that coincide except at the
$i$-th term, 
while the intersection of the $i$-th terms is the smallest possible.  These
functors, which appeared in \cite[Proposition 7]{Kh3}, satisfy braid relations,
and the same argument as in Rouquier \cite{Ro1} implies that the braid group
action on the derived category of sheaves on the disjoint union of $Y_{s
\mu},$ over all 
elements $s$ of the symmetric group $S_{k+1},$ is genuine. For generic 
$\mu,$ this action can only be extended trivially to braid cobordisms 
with branch points, due to the absence of nonzero natural transformations
between $F_1$ and $F_{\sigma_i}.$ 
\end{example}

A closely related example of braid group categorification via Rouquier
complexes will be discussed below in Section~\ref{sec-bc-RR}. 

%%%%%%%%%%%%%%%%%%%%%%%%%
%%
%% BRAID COBORDISMS 
%%
%%%%%%%%%%%%%%%%%%%%%%%%%%%

\section{Invariants of braid cobordisms} 
\label{sec-bc}

\subsection{Braid cobordisms and braid group categorifications}
\label{sec-bc1}

\subsubsection*{Combinatorial braid cobordisms}
\begin{figure}  \drawing{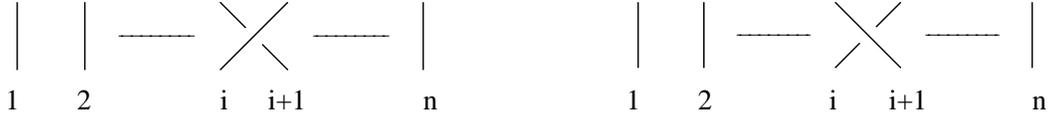}\caption{On the left: braid $\sigma_i.$ On 
 the right: braid $\sigma_i^{-1}.$}  \label{picasso3}\end{figure}

The category $\mc{BC}_n$ of 
$n$-stranded braid cobordisms admits a combinatorial version which 
we denote $\mc{BC}^c_n.$ Its objects are $n$-stranded braid words, 
i.e. arbitrary finite sequences of generators $\sigma_i,$ $1\le i\le n-1,$  
and their inverses $\sigma_i^{-1},$ see figure~\ref{picasso3}.

 A morphism between braid 
words is a finite sequence of transformations of two types: 

\medskip

I. Reidemeister moves of braid words, including 
 \begin{eqnarray*} 
    \tau_1 \sigma_i \sigma_i^{-1} \tau_2 & \leftrightarrow & 
   \tau_1 \tau_2, \\
   \tau_1 \sigma_i \sigma_j \tau_2 & \leftrightarrow & 
   \tau_1 \sigma_j \sigma_i \tau_2, \hspace{0.2in} |i-j|>1, \\ 
  \tau_1 \sigma_i \sigma_{i+1} \sigma_i \tau_2 & \leftrightarrow & 
   \tau_1 \sigma_{i+1}\sigma_i \sigma_{i+1} \tau_2. 
\end{eqnarray*} 
$\tau_1$ and $\tau_2$ are arbitrary braid words. The empty word is allowed. 
 Each move has several versions, and we won't write all of them, instead 
referring the reader to \cite[Section 3]{CS2} for details.  For 
instance, in the second move above we could have substituted $\sigma_i^{-1}$ 
for $\sigma_i,$ and, independently, $\sigma_j^{-1}$ for $\sigma_j.$ 

\medskip

II. Addition or removal of a single $\sigma_i$ or $\sigma_i^{-1}$ to or from a 
word: 
$$ \tau_1 \tau_2 \leftrightarrow \tau_1 \sigma_i^{\pm 1} \tau_2.$$ 

\begin{figure}  \drawing{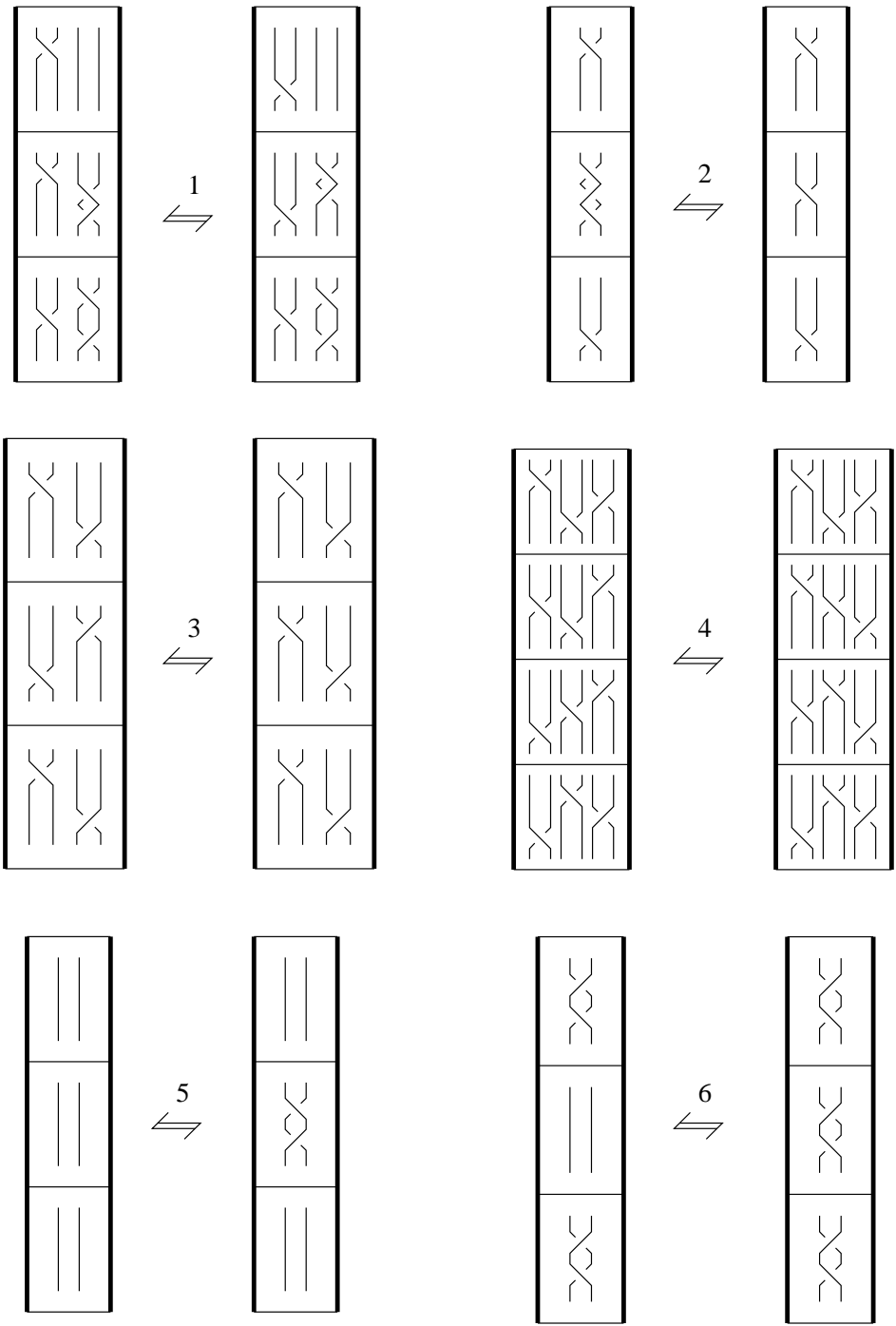}\caption{Carter and Saito's braid movie moves 1-6} 
 \label{picasso4}\end{figure}

\begin{figure}  \drawing{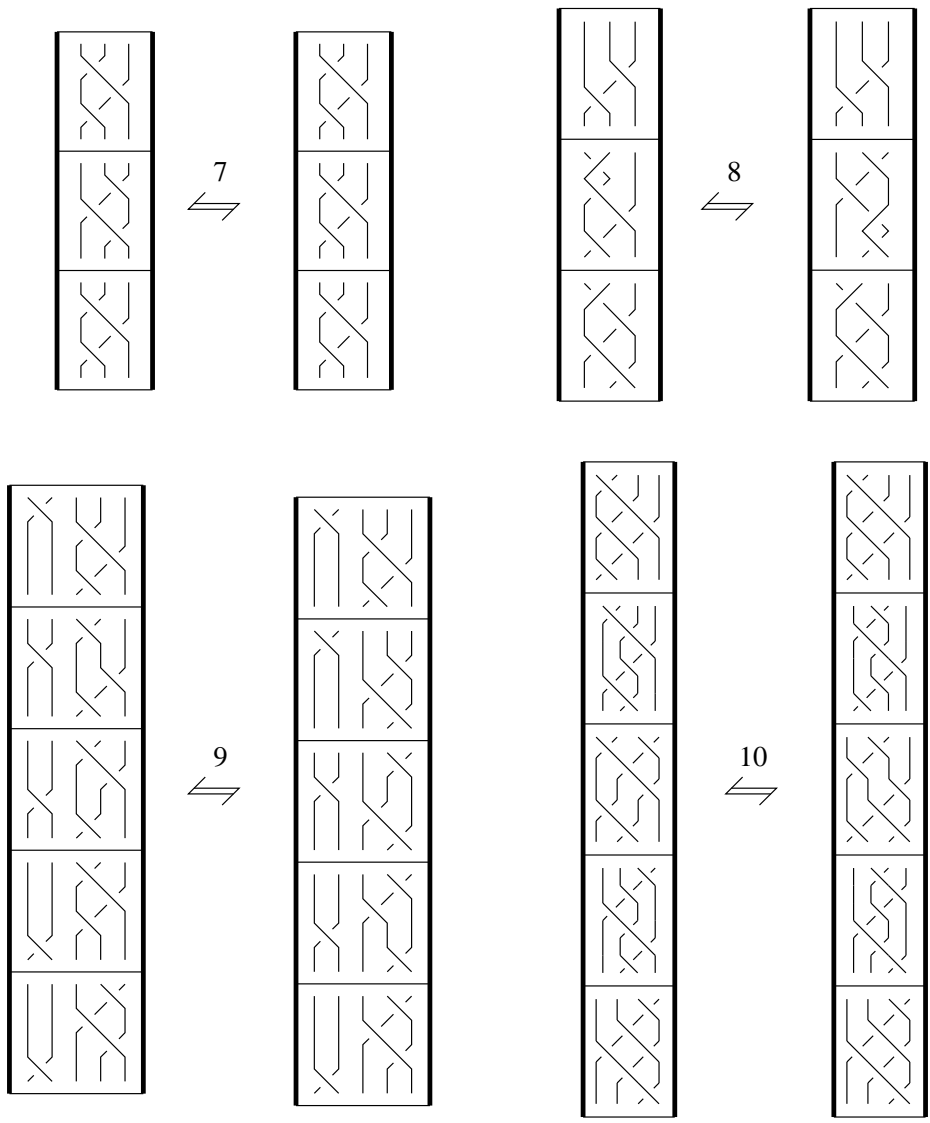}\caption{Braid movie moves 7-10}
 \label{picasso5} \end{figure}

\begin{figure}  \drawing{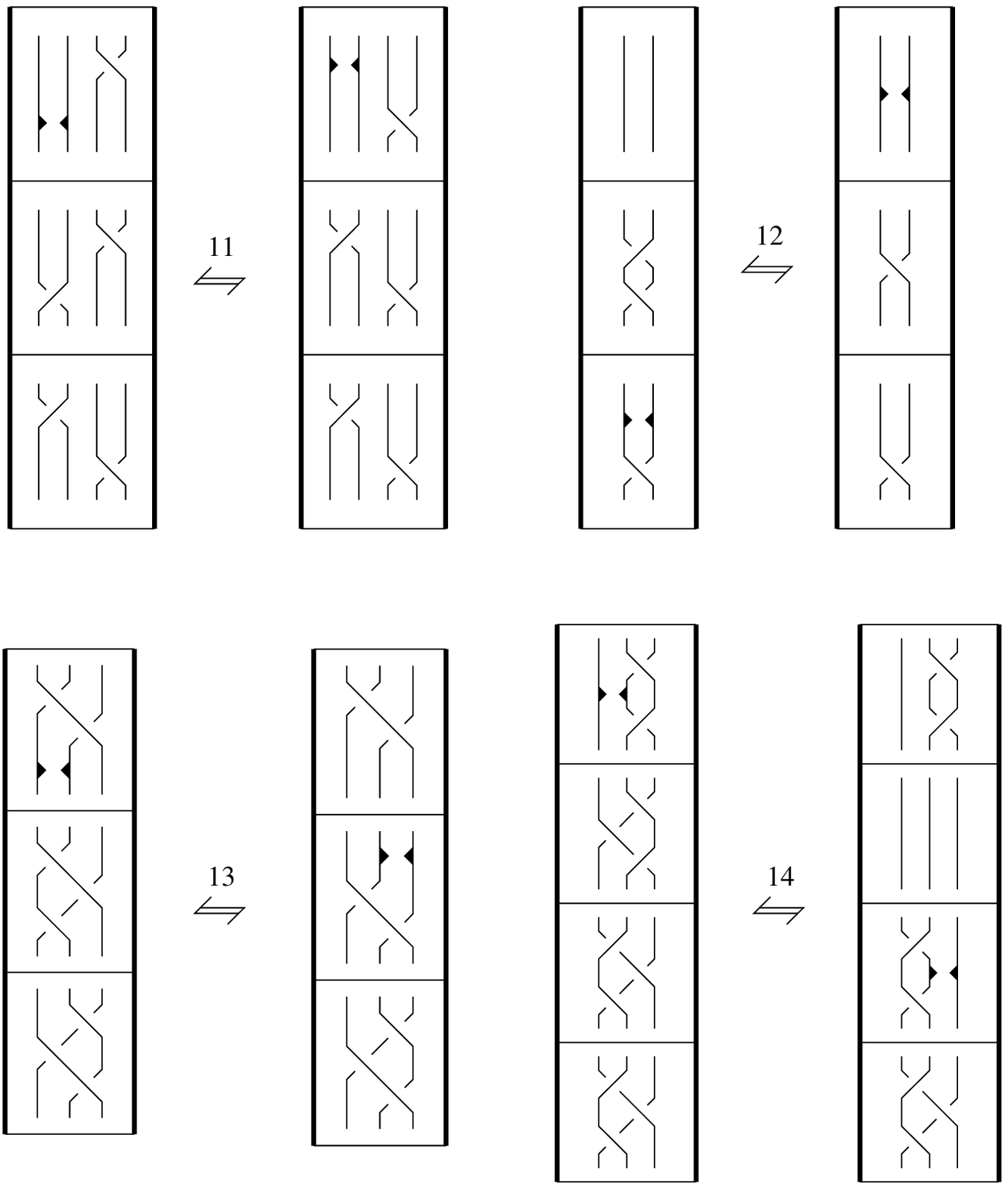}\caption{Braid movie moves 11-14} 
 \label{picasso6} \end{figure}

A braid \emph{movie} is a finite sequence of transformations of types I and II. 
Morphisms in the category $\mc{BC}_n^c$ are braid movies modulo
equivalence relations of two types: movie moves, shown in 
figures~\ref{picasso4}-\ref{picasso6}, and locality moves. The first ten 
movie moves are equations on sequences of Reidemeister moves of braids, 
while each sequence in moves 11 through 14 contains a unique type II transformation. 
A type II transformation is indicated by placing two little black triangles where 
a pair of adjacent vertical lines turns into $\sigma_i^{\pm 1}.$ In the figures,
movies run from top to bottom, and represent a morphism in $\mc{BC}_n^c$;
moves of these movies run left to right, representing equivalences between morphisms.


Each movie move has several variants, given by changing overcrossings to 
undercrossings and vice versa, whenever possible; reading each of 
the two sequences from bottom to top rather than from top to bottom 
(the latter is our default); reflecting all diagrams in a movie move about 
a vertical axis (if this results in a different version of the move), etc. 
Again, we refer to Carter and Saito \cite[Section3]{CS2} for complete 
explanations. 

\begin{figure}  \drawing{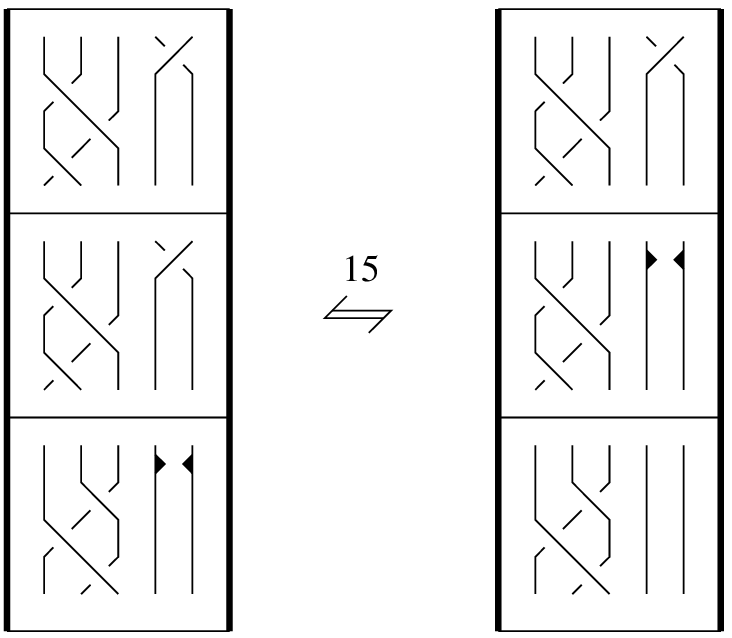}\caption{A locality movie move} 
 \label{picasso7} \end{figure}

Locality moves are 
$$ (\tau_1 \tau_2, \tau_1'\tau_2, \tau_1'\tau_2')  = 
    (\tau_1 \tau_2, \tau_1\tau_2', \tau_1'\tau_2'),$$
where each of the pairs $\tau_1, \tau_1'$ and $\tau_2,\tau_2'$ 
are related by a single braid move, of either type. An example is 
depicted in figure~\ref{picasso7}, where $\tau_1=1,$ the trivial 
braid word, $\tau_1'=\sigma_{i+3}^1,$ 
$\tau_2=\sigma_i\sigma_{i+1}\sigma_i,$ and 
$\tau_2'=\sigma_{i+1}\sigma_i\sigma_{i+1}.$  We will refer 
to locality moves as braid movie moves of type 15. 

Morphisms between braid words $w_1$ and $w_2$ are finite 
sequences of braid words related by type I and II transformations 
modulo the equivalence relation generated by braid movie moves and 
locality relations. 

\begin{figure}  \drawing{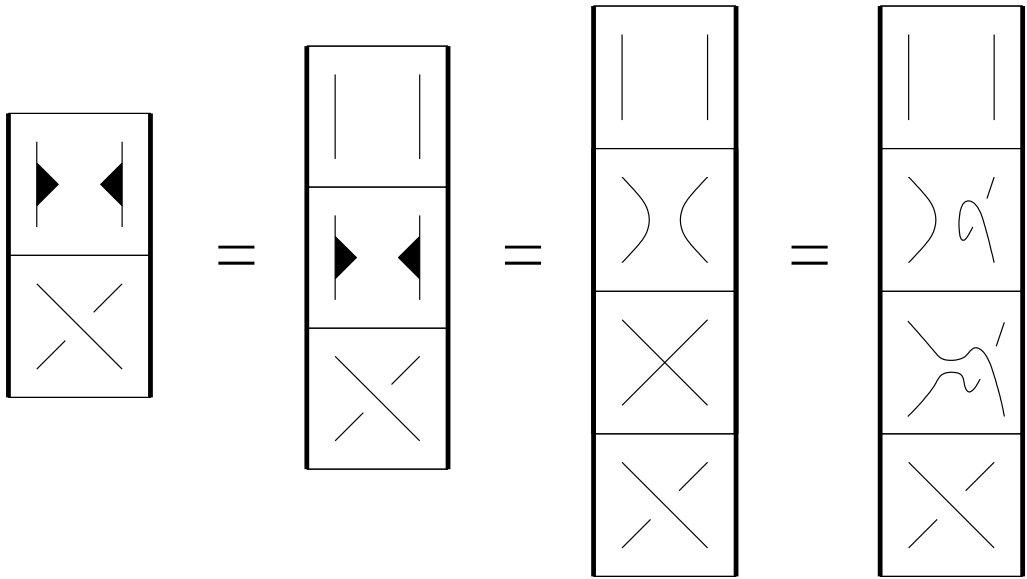}\caption{A type II braid transformation}
 \label{picasso8} \end{figure}

To explain the meaning of type II transformations, we describe 
a functor from $\mc{BC}_n^c$ to $\mc{BC}_n,$ the category 
of braid cobordisms. To each braid word we assign a braid in 
the usual way. To type I transformations we assign braid cobordisms 
given by isotopies of braids. Type II transformations
correspond to elementary braid cobordisms which contain a single 
simple branch point when projected onto $[0,1]^2.$ Our cryptic 
way of depicting type II transformations is deciphered in figure~\ref{picasso8}.
These cobordisms embed only in four dimensions rather than three
-- only after applying the twist on the right hand side of figure~\ref{picasso8}
in the fourth dimension can we embed the rest of the cobordism into three dimensions as the standard
simple cobordism between a 1-manifold and its surgery (the second and third
frames of the last movie in figure~\ref{picasso8}) provided by a Morse function
on a surface with an index one critical point (saddle point).

Since each braid movie move corresponds to a braid cobordism the above
construction is a functor 
from $\mc{BC}_n^c,$ the category of combinatorial braid cobordisms, to 
$\mc{BC}_n,$ the category of braid cobordisms. 

\begin{prop} This functor is an equivalence of categories: $\mc{BC}_n^c
\cong\mc{BC}_n.$
\end{prop} 

\begin{proof} See \cite{CS2} for a proof.
\end{proof}

\subsubsection*{Polarization of braid cobordisms}
Recall that a braid is \emph{positive} if 
it can be represented as a product of $\sigma_i,$ for $1\le i \le n-1.$ Likewise, a braid is called \emph{negative} if it can be written as a product of $\sigma_i^{-1},$
for $1\le i \le n-1.$ For braid cobordisms there also exists an analogue of positivity. We say that  a simple branch point of a braid cobordism is \emph{positive} 
if, when represented combinatorially via a type II move, it increases braid
positivity. Namely, positive type II moves are 
\begin{eqnarray*} 
  \tau_1 \tau_2 & \lra & \tau_1 \sigma_i \tau_2,  \\
  \tau_1 \sigma_i^{-1} \tau_2 & \lra & \tau_1 \tau_2. 
\end{eqnarray*} 
Negative type II moves are 
\begin{eqnarray*} 
  \tau_1 \tau_2 & \lra & \tau_1 \sigma_i^{-1} \tau_2,  \\
  \tau_1 \sigma_i \tau_2 & \lra & \tau_1 \tau_2. 
\end{eqnarray*} 
Each type II move is thus either positive or negative. Notice that type II
moves appear only in braid movie moves 11-15, 
and the property of being positive or negative is preserved by these movie
moves. 

Given a morphism $S$ in $\mc{BC}_n^c,$ the numbers $p_+(S),$ 
respectively $p_-(S),$ of positive, respectively negative, type II transformations
are invariants of $S.$ Likewise, given a morphism $S$ in $\mc{BC}_n,$ 
the numbers $p_+(S),$ respectively $p_-(S),$ of positive, respectively negative,
double branch points are invariants of $S.$ A braid cobordism is a braid isotopy 
if and only if  $p_+(S)=p_-(S)=0.$ We say that a morphism $S$ in 
$\mc{BC}_n^c$ or $\mc{BC}_n$ is positive if $p_-(S)=0$ and negative 
if $p_+(S)=0.$ Positive braid cobordisms constitute a subcategory in 
$\mc{BC}_n$ (and the same for negative cobordisms). Positive braid movies constitute a subcategory in $\mc{BC}_n^c$ (and the same for negative braid movies). Our 
equivalence of categories $\mc{BC}_n^c\cong \mc{BC}_n$ takes
positive, respectively negative, braid movies to positive, respectively negative,
braid cobordisms. 

\begin{figure}  \drawing{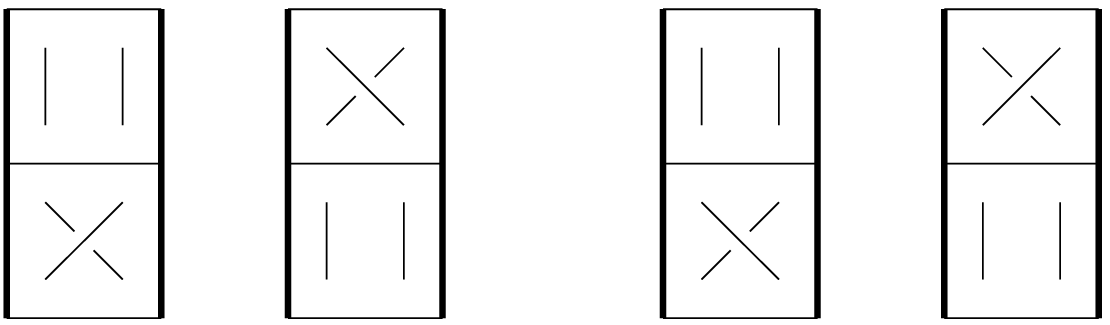}\caption{On the left: two movies of
a positive branch point. On the right: two movies of a negative branch point.}
 \label{picasso9} \end{figure} 

\begin{prop} Any braid cobordism can be written, in more than one way, 
as a composition of a positive and a negative braid cobordism (in either order). 
\end{prop} 

\begin{figure}  \drawing{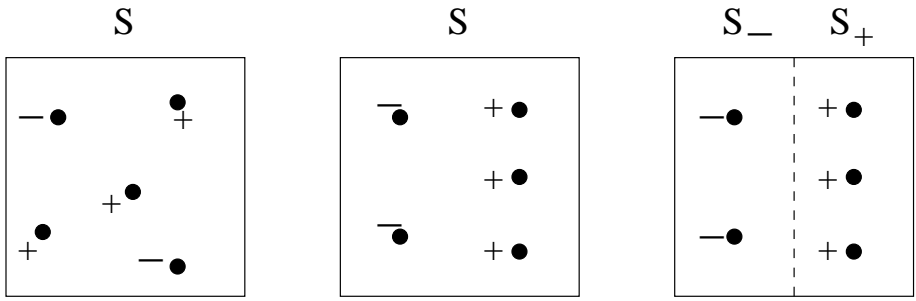}\caption{Moving branch points} 
 \label{picasso10} \end{figure}

\begin{proof} Project a braid cobordism $S$ onto the unit square, as in figure~\ref{picasso1}.
Next, isotop branch points so that the negative branch points  are to the left of 
the positive ones, see figure~\ref{picasso10}, and lift this isotopy to 
an isotopy of $S.$ We can now write $S$ as the composition of two cobordisms, 
$S=S_+ S_-,$ with $S_+$ positive and $S_-$ negative (this decomposition is 
not unique even when $S$ has 
no double points). Likewise, we can write $S$ as $S_-' S_+'$ where $S_-'$ is 
negative and $S_+'$ is positive.
\end{proof} 

A braid cobordism is positive if and only if it can be realized, \`a la Rudolph \cite{Ru}, as a piece of an 
algebraic curve in $\C^2.$ The local model of a simple positive cobordism
is the algebraic curve
\begin{equation} \label{z2z} 
z_2^2=z_1
\end{equation}
in $\C^2\cong\R^4$ with its projection
to $(x_1,y_1)=($\,Re\,$z_1$, Im\,$z_1)$. For any fixed $(x_1,y_1)$ the $z_2$
coordinate gives two points in $\C$ which wind round each other as $x_1$
varies, giving a braid for $y_1\ne0$. When $y_1=0$ the two strands cross.
This shows the essential four-dimensionality of positive type II transformations.
Similarly the simple negative type II transformation has local model $z_2^2=
\bar z_1$.

\subsubsection*{Categorification of the braid group} 

\begin{defn} \label{d-one} A categorification of a complex representation
$V$ of the braid group $\mathrm{Br}_n$ is a (projective) monoidal functor $F$ 
from $\mc{BC}_n$ to the category of exact functors $\mathrm{Fun}(\mc{C})$
in a triangulated category $\mc{C},$ and an isomorphism 
$V\cong K(\mc{C})\otimes \C $ of braid group representations. 
\end{defn} 

Thus, to each braid $g$ there is assigned an exact endofunctor $F_g$ in $\mc{C}$
and to each braid cobordism $S$ from $g$ to $h$ a natural transformation
$F_S: F_g \lra F_h,$ well-defined up to rescaling by invertible 
elements of the ground ring of $\mc{C}.$ These must satisfy consistency 
conditions. The linear operator $[F_g]$ on $K(\mc{C})\otimes \C$ must 
coincide with the action of $g$ on $V$ under the above isomorphism of 
these two complex vector spaces. 

A genuine action of a braid group on $\mc{C}$ extends in a trivial way 
to a categorification in our sense, by making $F_S=0$ whenever $S$ 
has a double point and defining $V$ to be $K(\mc{C})\otimes \C.$ 
Likewise, given a categorification as above, we 
can modify it, by setting $F_S=0$ if, say, $S$ has a negative branch point,
and otherwise doing nothing to $F_S.$ 
More generally, given a categorification and two elements $\lambda_+, \lambda_-$
of the ground ring of $\mc{C},$ we can rescale $F_S$ by $\lambda_+^{p_+(S)}
\lambda_-^{p_-(S)}$ to get another categorification.  
 
Naturally, we'd like our categorifications to carry nontrivial information
about braid cobordisms, and, whenever possible, we will avoid trivial or semi-trivial 
categorifications. By a semi-trivial categorification of a braid group representation
we mean one with $F_S=0$ for all $S$ with a negative branch point, or 
for all $S$ with a positive branch point. 

\subsection{Categorification from the rings $A_n$} \label{An}

We now construct a categorification of the reduced Burau representation via
the rings $A_n.$ The braid group action on the homotopy category $\mc{C}(A_n),$
described in the 
previous section, lifts to the action on the homotopy category of 
complexes of \emph{graded} $A_n$-modules. Indeed, $A_n$ has the natural 
grading by path length.  The bimodules $A_n$ and $P_i \otimes \hspace{0.03in}
_i P$ are graded, and the bimodule maps $\beta_i$ and $\gamma_i$ have 
degrees $0$ and $2$ respectively. Hence, we can upgrade $R_i$ and  $R_i'$
to complexes of graded bimodules, shifting the grading of $P_i \otimes \hspace{0.03in}
 _i P$ in $R_i'$ down by $2$ to make the differential preserve the grading:
$$R_i':  \hspace{0.2in}  0 \lra A_n \stackrel{\gamma_i}{\lra}
P_i \otimes \hspace{0.03in} _iP\{-2\} \lra 0. $$
Tensoring with $R_i$ and $R_i'$ induces a braid group action on the 
category of complexes of graded $A_n$-modules up to chain homotopy. 
We restrict to bounded complexes and denote the resulting category also by
$\mc{C}(A_n).$
It is easy to see that on the Grothendieck group the braid group action descends
to the reduced Burau representation of the braid group. 

To extend to braid cobordisms, we assign a natural transformation of 
functors to each morphism in $\mc{BC}_{n+1},$ i.e. to each braid movie. This transformation 
will be defined up to an overall minus sign, resulting in a projective functor
from $\mc{BC}_{n+1}$ to $\mathrm{Fun}(\mc{C}(A_n)).$ The functors associated
to braid words are given by tensor products with complexes of graded bimodules.
The latter complexes are tensor products of the complexes $R_i$ and $R_i'.$
Denote by $R(\sigma)$ the complex of graded bimodules associated with 
the braid word $\sigma,$ and by $F_{\sigma}$ the functor of tensoring with
$R(\sigma).$ 

A homomorphism of complexes of bimodules induces a natural transformation of the 
corresponding functors, and our natural transformations will be of this form. 

To a Reidemeister move of braid words we associate an isomorphism in the homotopy 
category of complexes of graded $A_n$-bimodules, namely, isomorphisms  
(\ref{eq-1}), (\ref{eq-2}), (\ref{eq-3}). Each of these isomorphisms is unique
up to an overall minus sign. Indeed, suppose we are given two invertible complexes of graded bimodules 
over a graded ring $A$ and an isomorphism $u$ between these complexes. Then
any isomorphism between these complexes has the form $m(z)u$ where 
$m$ is the multiplication by a degree $0$ invertible element of the center of $A$ 
(see \cite[Section 3]{Kh2}, for instance). The degree $0$ part of $A_n$ is just $\Z,$ 
its invertible elements are $1$ and $-1,$  and the above claim follows. 

To type II transformations we associate the following homomorphisms of 
complexes of bimodules. 

The homomorphism $w_1$ from $R(\emptyset)$ to $R(\sigma_i)= R_i$:   
\begin{equation} \label{positivemove1} 
\begin{CD} 
 0    @>>>  0 @>>>  A_n @>>> 0  \\
   @VVV   @VVV   @VV{\mathrm{id}}V  @VVV   \\
  0   @>>>  P_i \otimes \hspace{0.03in} _iP @>{\beta_i}>> A_n @>>> 0.\!
\end{CD}  
\end{equation} 
Here $R(\emptyset)\cong A_n$ is the identity bimodule complex, associated with the 
empty braid word $\emptyset,$ corresponding to the trivial braid. 

The homomorphism $w_2$ from $R(\sigma_i^{-1})$ to $R(\emptyset)$:  
\begin{equation} \label{positivemove2}
\begin{CD} 
 0   @>>> A_n @>{\gamma_i}>>  P_i \otimes \hspace{0.03in} _iP @>>> 0    \\ 
      @VVV          @VV{\mathrm{id}}V                   @VVV  @VVV   \\
 0    @>>>  A_n @>>>  0 @>>> 0.\!
\end{CD}  
\end{equation} 

Notice that the above two homomorphisms are assigned to elementary movies of \emph{positive} branch points. In both cases, the map is the identity on $A_n.$ 

\medskip

The homomorphism $w_3$ from $R(\sigma_i)$ to $R(\emptyset)$: 
\begin{equation} \label{negativemove1}
\begin{CD} 
   0   @>>>  P_i \otimes \hspace{0.03in} _iP @>{\beta_i}>> A_n @>>> 0  \\
      @VVV   @VVV   @VV{\delta_i}V  @VVV   \\
   0  @>>>  0 @>>>  A_n @>>> 0,\!
\end{CD}  
\end{equation} 
where $\delta_i(1) = X_{i-1}- X_{i+1}.$ This condition determines the bimodule
homomorphism $\delta_i$ uniquely. Notice that $\delta_i \beta_i=0.$ 

The homomorphism $w_4$ from $R(\emptyset)$ to $R(\sigma_i^{-1})$: 
\begin{equation} \label{negativemove2}
\begin{CD} 
 0    @>>>  A_n @>>>  0 @>>> 0  \\
   @VVV   @VV{\delta_i}V   @VVV  @VVV   \\
  0   @>>> A_n @>{\gamma_i}>>  P_i \otimes \hspace{0.03in} _iP @>>> 0.\!
\end{CD}  
\end{equation} 

Thus, to the four braid movies in figure~\ref{picasso9} we assign the 
homomorphisms $w_1,$ $w_2,$ $w_3,$ $w_4$ correspondingly (from left to right).

\medskip

Given a braid movie $S$ from braid word $g$ to braid word $h$ we 
assign to $S$ the homomorphism $R(S): R(g) \lra R(h)$ of complexes 
of graded bimodules, the latter a composition of homomorphisms, described 
above, corresponding to elementary braid movies.  This homomorphism 
preserves cohomological degree of complexes, and has internal degree 
$2p_-(S),$ since the bimodule homomorphism $\delta_i,$ corresponding to 
a movie of a negative branch point, has degree $2.$ \medskip

\noindent\emph{Proof of Theorem \ref{first-theorem}}.
We show that the above assignment 
extends to a projective functor from $\mc{BC}_{n+1}$ to 
$\mathrm{Fun}(\mc{C}(A_n)).$ This is equivalent to checking that the 
assignment is compatible with the braid movie moves 1--15. 
Each braid movie move contains two movies $S_1,S_2$ from a braid word $g$ to a braid word $h$ and we need to prove the identity $R(S_2) =\pm R(S_1),$
i.e. that $R(S_2)$ is homotopic to either $R(S_1)$ or $-R(S_1).$ 

Movie moves 1--10 are easy to take care of. In each of these moves the
movies $S_1$ and $S_2$ are invertible (have no branch points), being 
a sequence of braid Reidemeister moves, and 
the corresponding bimodule homomorphisms $R(S_1), R(S_2)$ are 
invertible. Two isomorphisms $R(S_1), R(S_2)$ between invertible 
graded $A_n$-bimodules $R(g), R(h)$ differ by multiplication by 
an invertible degree $0$ central element of $A_n.$ Hence $R(S_2) = \pm R(S_1)$ 
(this argument was used in \cite{Kh2, BN} and \cite{St1}
when proving that certain invariants of tangle cobordisms are 
well-defined up to rescaling).  

Movie move 15 (locality moves) follow at once from the following observation. 
Given four bimodules over a ring $A$ and bimodule homomorphisms 
$u: N_1 \lra N_2$ and $v: N_3 \lra N_4,$ the bimodule homomorphism 
$u\otimes v : N_1 \otimes N_3 \lra N_2 \otimes N_4$ can be written 
as $(u\otimes \mathrm{Id}) \circ (\mathrm{Id} \otimes v),$ as 
well as $(\mathrm{Id}\otimes v) \circ (u \otimes \mathrm{Id}).$ 

Movie move 11 follows easily, since the canonical isomorphisms of bimodule complexes 
$ R_i \otimes R_j \cong R_j \otimes R_i,$  
$R_i \otimes R_j' \cong R_j' \otimes R_i,$ etc. are compatible with the
homomorphisms $w_1,\dots, w_4.$ 

This leaves us with braid movie moves 12--14. 
Each of the two movies in each of these three moves contains a 
unique type II transformation, which could be either positive or 
negative. 

We first deal with the simpler case of a positive type II transformation. 
Then the two movies $S_1, S_2$ in any of the movie moves 12--14 
go from a braid word $g$ to a braid word $h$ obtained from $g$ 
by either inserting a positive braid $\sigma_i$ or deleting a 
negative braid $\sigma_i^{-1}.$ Each movie $S_1, S_2$ is 
a composition of Reidemeister braid moves and a single 
homomorphism (either $w_1$ or $w_2,$ see above) of bimodule 
complexes. It's easy to check that the space of grading-preserving 
homomorphisms from $R(\emptyset)$ to $R(\sigma_i)$ is 
$\Z,$ with generator $w_1.$ Likewise, the space of 
grading-preserving homomorphisms from $R(\sigma_i^{-1})$ 
to $R(\emptyset)$ is $\Z,$ with generator $w_2.$ This 
implies that the space of grading-preserving homomorphisms 
from $R(g)$ to $R(h)$ is $\Z,$ with $R(S_1), R(S_2)$ both 
being generators. Hence $R(S_2) = \pm R(S_1),$ which takes 
care of movie moves 12--14 in the positive case. 

\medskip

All possible versions of the movie move 12 with the negative 
type II transformation follow by direct computation. A sample 
computation for the version in figure~\ref{picasso6} is included below. 

The maps $F(S_1), F(S_2)$ go from the bimodule $A_n$ to 
$R_i'.$ The first map is given by the composition written below 
\begin{equation} 
\begin{CD} 
 0    @>>> 0 @>>>  A_n @>>> 0  @>>> 0     \\ 
      @VVV    @VVV      @VV{(\mathrm{id},\gamma_i,0)^t}V    @VVV  @VVV   \\
 0    @>>> P @>{d}>>  
  A_n \oplus P\{-2\} \oplus P @>{d}>> P\{-2\}  @>>> 0    \\ 
      @VVV     @VVV       @VV{((X_{i-1}-X_{i+1})\mathrm{id},0,0)}V     
                @VV{\mathrm{id}}V  @VVV   \\
 0   @>>>  0 @>>>  A_n @>{\gamma_i}>> P\{-2\} @>>> 0,\!  
\end{CD}  
\end{equation} 
where $P$ is a shorthand for $P_i \otimes \hspace{0.03in} _iP,$
while the map $F(S_2)$ is just $w_3.$ Clearly, $F(S_1)= F(S_2).$ 

\medskip

To treat movie moves 13--14 with a negative type II transformation
we need to do some preliminary work. 

Given an endofunctor $F$ on $\mc{C}(A_n)$ and a central element 
$a\in A_n,$ we denote by $l_a$ the natural transformation $F\lra F$ 
given by left multiplication by $a:$ 
$$ \spreaddiagramcolumns{.8pc}\xymatrix{
F \cong \mathrm{Id}\circ F \rto^{a\circ\mathrm{id}_F} & 
   \mathrm{Id}\circ F \cong F.}
$$ 
Likewise, denote by $r_a$ the natural transformation $F\lra F$ 
given by right multiplication by $a:$ 
$$ \spreaddiagramcolumns{.8pc}\xymatrix{
F \cong F\circ \mathrm{Id} \rto^{\mathrm{id}_F\circ a} & 
   F\circ \mathrm{Id} \cong F.}
$$ 
Suppose now that $a$ is a linear combination of $X_j$s, 
 $$ a= \sum_{j=1}^n a_j X_j,$$ 
and take $F=R_i$ or $R_i'.$  

\begin{prop} \label{lara}
For any $a$ as above and $1\le i\le n,$ 
 \begin{eqnarray*}  
 (l_a - r_{a-a_i(X_{i-1}+2 X_i + X_{i+1})})R_i  & = & 0 \\
\mathrm{and}\quad (l_a - r_{a-a_i(X_{i-1}+2 X_i + X_{i+1})})R_i'  & = & 0
 \end{eqnarray*} 
in $\mathrm{Fun}(A_n).$ 
\end{prop} 
In other words, the endomorphism 
$l_a - r_{a-a_i(X_{i-1}+2 X_i + X_{i+1})}$ of the complexes $R_i$
and $R_i'$ of bimodules is zero homotopic. 

\begin{proof} Consider the homotopy $h$ from $R_i$ to itself 
which takes $A_i$ to $P_i \otimes \hspace{0.02in} _iP$ via 
$\gamma_i$ and $P_i \otimes \hspace{0.02in} _iP$ to $0.$ 
Then the zero homotopic endomorphism $hd+dh$ of $R_i$ is given on 
generators $(i)\otimes (i) \in P_i \otimes \hspace{0.02in} _iP$ 
and $1\in A_n$ by 
 \begin{eqnarray*} 
 (hd+dh)((i)\otimes (i)) & = & X_i \otimes (i) + (i) \otimes X_i, \\
 (hd+dh)(1)  & = & X_{i-1}+2X_i + X_{i+1}. 
 \end{eqnarray*} 
Next, we compute 
 \begin{eqnarray*} 
(l_a - r_{a-a_i(X_{i-1}+2 X_i + X_{i+1})})((i)\otimes(i)) & = & 
   a_i X_i \otimes (i) +  a_i (i) \otimes X_i , \\
(l_a - r_{a-a_i(X_{i-1}+2 X_i + X_{i+1})})(1) & = & 
   a_i(X_{i-1}+2 X_i +X_{i+1}). 
\end{eqnarray*} 
Therefore, 
$$(l_a - r_{a-a_i(X_{i-1}+2 X_i + X_{i+1})})R_i= 
 a_i(hd+dh)R_i ,$$ 
and the left hand side is homotopic to zero. 
A similar argument works for $R_i'$, using $\beta_i$ in the homotopy.
\end{proof}

\begin{rmk} When $i=1$ or $i=n$ we should omit $X_{i-1},$ 
respectively $X_{i+1},$ from the above formulas and computations. 
\end{rmk}

\begin{rmk} The braid group action on $\mc{C}(A_n)$ descends 
to an action on the center of $A_n$ which happens to factor through 
the action of the symmetric group and is given by the above formulas.
The generator $s_i$ of the symmetric group takes $a$ to 
$a-a_i(X_{i-1} + 2 X_i + X_{i+1}).$  
This action also appears, in a similar context,  in \cite[Section 6]{FKS}. \end{rmk}

\begin{cor} \label{a-cor} The following endomorphisms of 
the bimodule complexes $R_i$ and $R_i'$ are null-homotopic: 
$$ l_{X_{i-1}+X_i}+r_{X_i+X_{i+1}}, \hspace{0.2in} 
    l_{X_{i}+X_{i+1}}+r_{X_{i-1}+X_{i}}, \hspace{0.2in} 
  l_{X_{i-1}-X_{i+1}} + r_{X_{i+1}-X_{i-1}} $$ 
\end{cor} 

\begin{figure}  \drawing{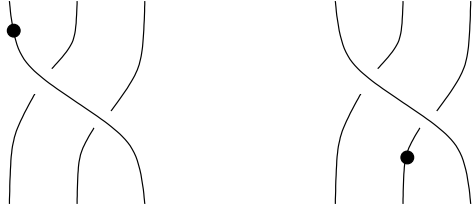}\caption{The endomorphisms $l_{X_1}$ and $r_{X_1+X_2}$ of $R_1' \otimes R_2'$ } 
 \label{picasso11} \end{figure}

\medskip

\begin{figure}  \drawing{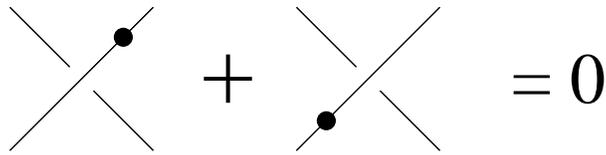}\caption{The endomorphisms 
$l_{X_i+X_{i+1}}$ and $-r_{X_{i-1}+X_{i}}$ of $R_i$ are 
chain homotopic.}  
 \label{picasso12} \end{figure}
We depict the endomorphisms $l_{X_{i-1} +X_{i}}$ and 
$r_{X_{i-1}+X_i}$ of a bimodule complex associated 
with a braid word 
by a dot placed near the end of the $i$-th top, respectively bottom, 
strand; see an example in figure~\ref{picasso11}. 
One of the cases of the corollary is written in figure~\ref{picasso12}. 

\medskip

\begin{figure}  \drawing{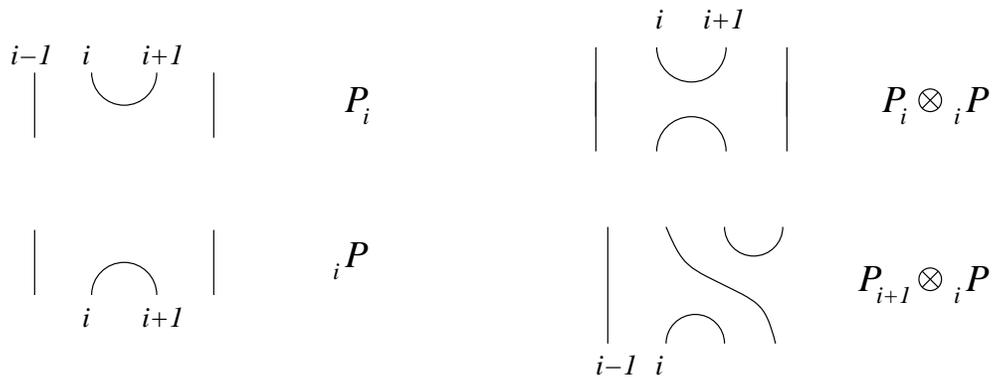}\caption{Graphical notation 
for projective modules $P_i, \hspace{0.02in} {}_iP$ and their tensor products.}
 \label{picasso13} \end{figure}

We now consider the instance of braid movie move 13 depicted in 
figure~\ref{picasso6} and, for convenience, introduce additional 
graphical notations in figure~\ref{picasso13}.  

\begin{figure}  \drawing{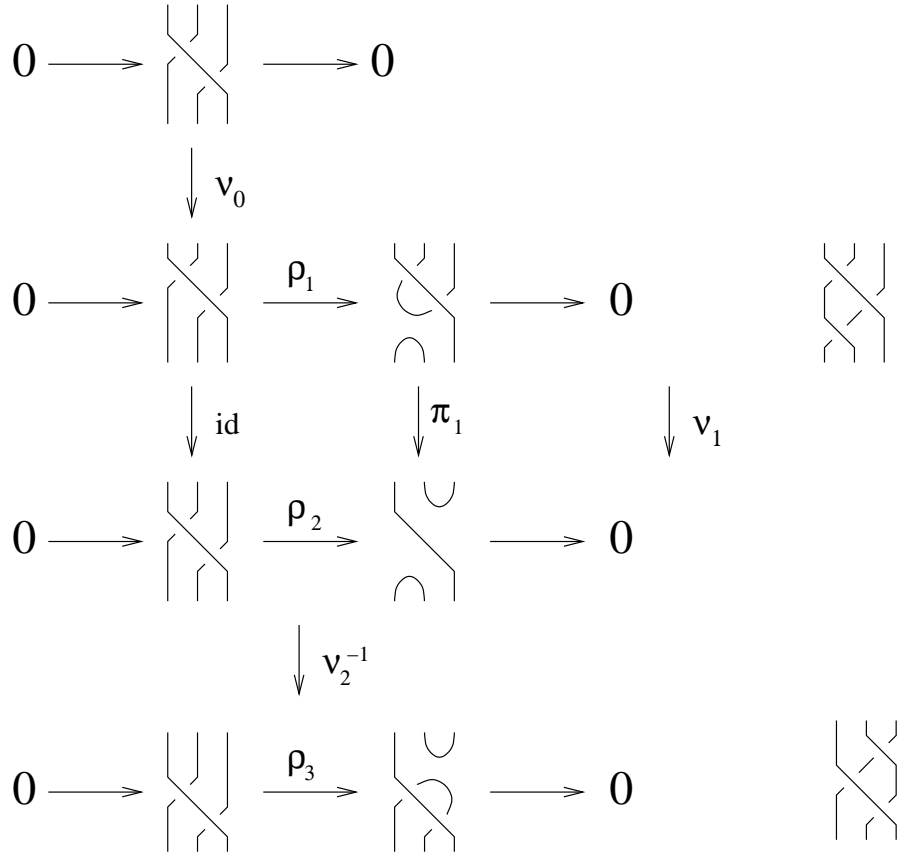}\caption{The map $F(S_1)$ in move 13 presented in graphical form.} 
 \label{picasso14} \end{figure}

\begin{figure}  \drawing{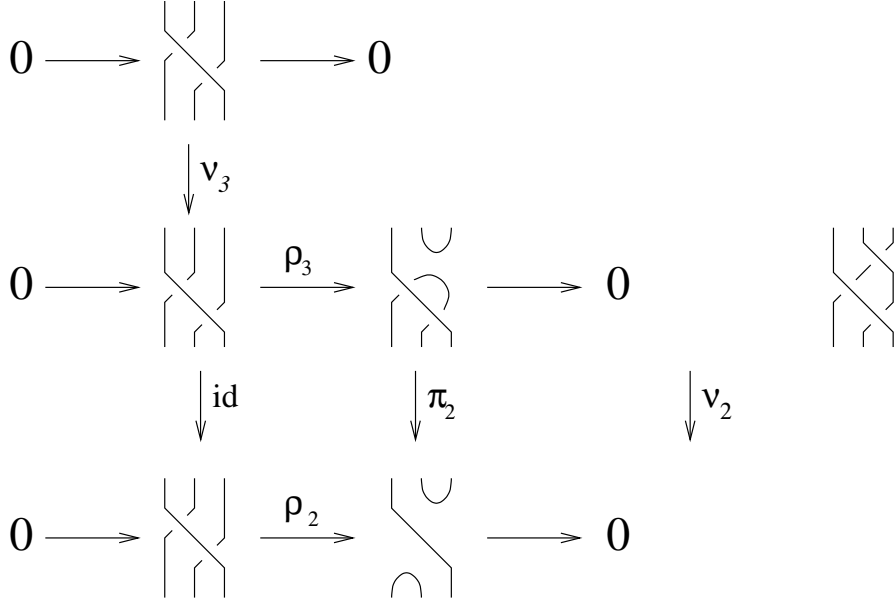}\caption{The map $\nu_2F(S_2)$ in move 13 written in graphical form.} 
 \label{picasso15} \end{figure}

Next, we rewrite $F(S_1)$ in graphical form, as in figure~\ref{picasso14}.
The map $\nu_0$ is right multiplication by $X_{i-1}-X_{i+1}.$  
The second row in the figure presents the 
complex $R(\sigma_i^{-1}\sigma_{i+1}^{-1} 
 \sigma_i^{-1})$ as the cone of the map of complexes 
 $$ 0 \lra R(\sigma_i^{-1}\sigma_{i+1}^{-1} ) \stackrel{\rho_1}{\lra} 
   R(\sigma_i^{-1}\sigma_{i+1}^{-1} ) \otimes_{A_n} P_i \otimes \hspace{0.03in}
    _iP \lra 0, $$ 
where $\rho_1= \mathrm{Id}\otimes \gamma_i.$ The map $\pi_1$ is 
a homotopy equivalence of complexes, given by taking the quotient of 
$R(\sigma_i^{-1}\sigma_{i+1}^{-1} ) \otimes_{A_n} P_i\otimes \hspace{0.03in} _iP$ by a contractible direct summand. The map $\nu_1$ 
is a homotopy equivalence. The map $\nu_2^{-1}$ is the inverse of the quotient
map $\nu_2,$ depicted in the lower half of figure~\ref{picasso15}. 
The map $\nu_2$ is the 
identity on $R(\sigma_i^{-1}\sigma_{i+1}^{-1} )= R_i' \otimes R_{i+1}'$ 
and the quotient map 
$\pi_2$ (analogous to $\pi_1$) on the second term of the cone. Since $\nu_2$
is a homotopy equivalence, its inverse is well-defined. 

\begin{figure}  \drawing{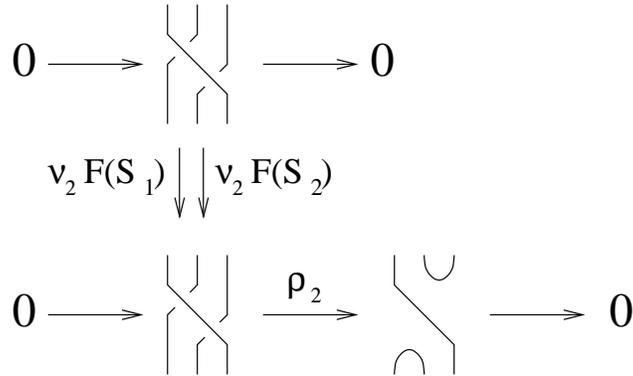}\caption{The maps $\nu_2F(S_1)$ and
$\nu_2 F(S_2)$.} 
 \label{picasso16} \end{figure}

In figure~\ref{picasso15} we depict the composition of $F(S_2)$ with $\nu_2.$
The map $\nu_3$ is left multiplication by $X_i - X_{i+2}.$ 
The equation  $F(S_1)=\pm F(S_2)$ is equivalent to 
$\nu_2 F(S_1) = \pm \nu_2 F(S_2).$ The map $\nu_2 F(S_1)$ 
is right multiplication by $X_{i-1}-X_{i+1}$ and $\nu_2 F(S_2)$ is left
multiplication by $X_i - X_{i+2},$ see figure~\ref{picasso16}. The endomorphisms
$r_{X_{i-1}- X_{i+1}}$ and $-l_{X_i - X_{i+2}}$ of $R_i' \otimes R_{i+1}'$
are chain homotopic (apply Corollary~\ref{a-cor} several times), 
via a degree 2 homotopy $h.$ We need to show that 
$r_{X_{i-1}- X_{i+1}}+l_{X_i - X_{i+2}}$ is null-homotopic when viewed 
as a homomorphism from  $R_i' \otimes R_{i+1}'$ to the complex depicted at the bottom of figure~\ref{picasso16}.  

A direct but nontrivial computation, which requires explicitly writing 
down the complex $R_i' \otimes R_{i+1}'$ and which we leave to the 
reader,  shows that the relation 
$$dh + hd = r_{X_{i-1}- X_{i+1}}+l_{X_i - X_{i+2}}$$
implies  $\rho_2 h=0,$ as a map of abelian groups. Consequently, 
$$\nu_2 (F(S_1) + F(S_2))= r_{X_{i-1}- X_{i+1}}+l_{X_i - X_{i+2}}$$ 
is null-homotopic as a map between the complexes in figure~\ref{picasso16}.
Thus, $F(S_1) + F(S_2)=0$ and this instance of braid movie move 13 follows.

Three other cases of movie move 13 with a negative branch point and 
with the top frame being one of the braids $\sigma_i \sigma_{i+1},$ 
$\sigma_{i+1}\sigma_i,$ $\sigma_{i+1}^{-1}\sigma_i^{-1}$ (instead of 
$\sigma_i^{-1}\sigma_{i+1}^{-1}$) follow via nearly identical arguments. 
The remaining instances of movie move 13 with a negative branch point 
are obtained by reading the frames in the above four cases from bottom 
to top. These instances follow (on the topological level) from the first four cases 
and from some of the moves 1--12. 

\bigskip

\begin{figure}  \drawing{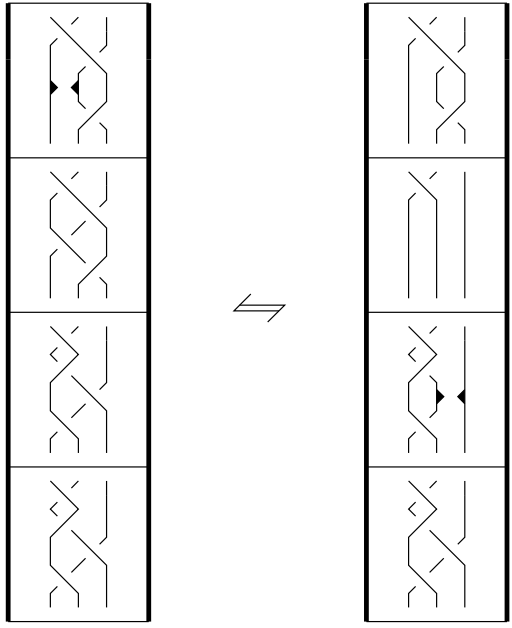}\caption{A modified braid movie move
14.} \label{picasso17} \end{figure}

Finally, we deal with braid movie move 14. The instance of this move shown
in figure~\ref{picasso6} is equivalent to the 
move depicted in figure~\ref{picasso17} (we formed the product of this move
with the identity cobordism on $\sigma^{-1}_i$). 

\begin{figure}  \drawing{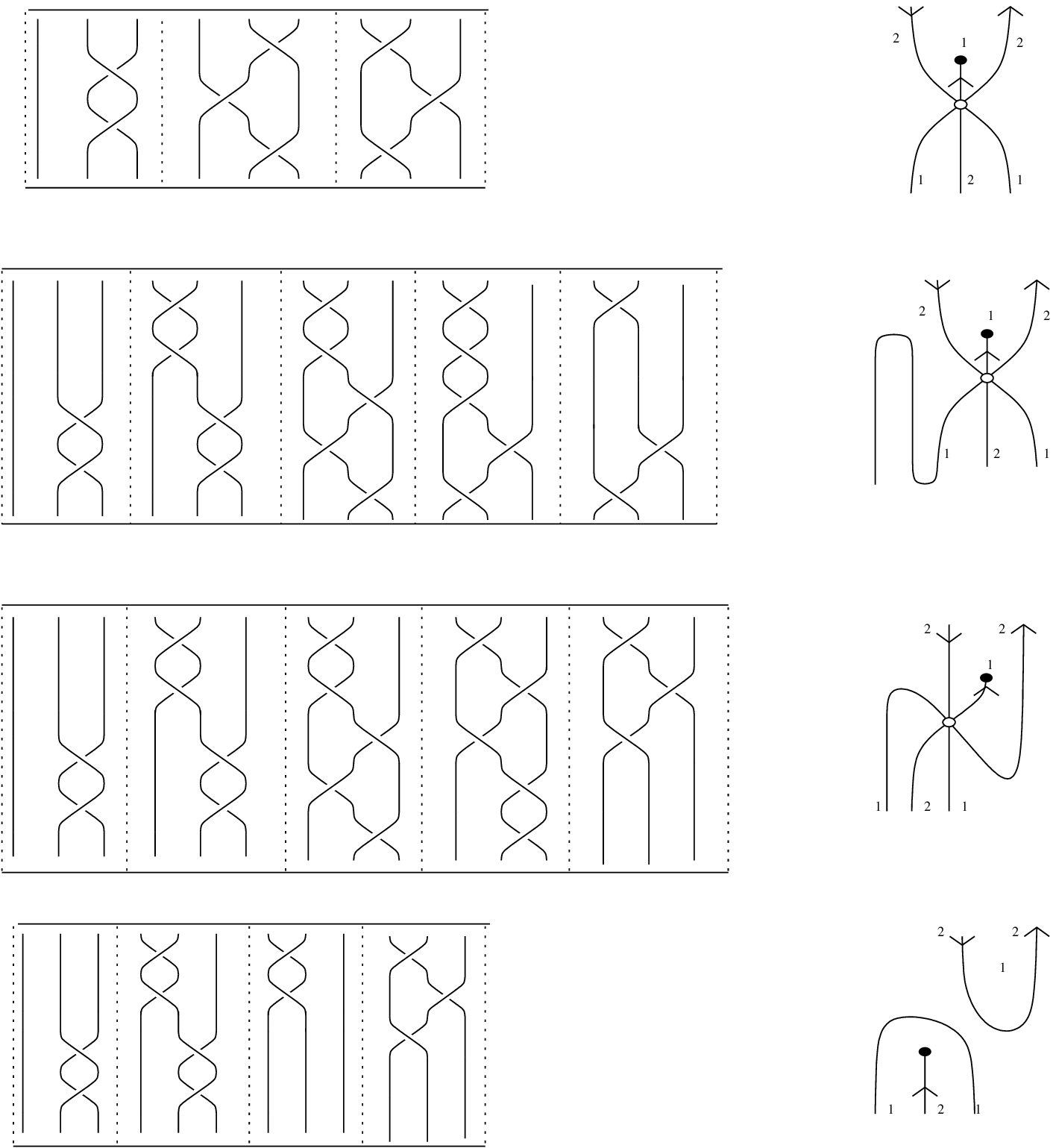}\caption{(courtesy of Scott Carter) Move 14 follows from the other moves. Diagrams on the right hand side are notations in \cite{CS2} for braid cobordisms.} \label{carter} \end{figure}

Look at the branch point indicated by the two triangles on the top strand of the left hand movie of figure~\ref{picasso17}. Move 13 allows us to drag this branch point up and across the third strand of the braid. Thus, the left hand movie is equivalent to the movie where the branch point is created to the right of this strand, just like in the other movie in figure~\ref{picasso17}.
These two movies are equivalent via a composition of several earlier moves. Other cases of move 14 follow in the same fashion, implying that move 14 is redundant. A complete derivation of move 14 from move 13 and other moves can be found in figure~\ref{carter}, kindly provided to us by Scott Carter.

Our proof of Theorem~\ref{first-theorem} ends here. \hfill$\square$\medskip

\subsection{Categorification via Rouquier complexes}
\label{sec-bc-RR}

Consider the polynomial 
algebra $k[x_1, \dots, x_n]$ over a field $k$ and its subalgebra 
$A= k[x_1-x_2, x_2-x_3, \dots, x_{n-1}-x_n]$ of polynomials in 
$x_i -x_{i+1},$ $1\le i\le n-1.$ The symmetric group $S_n$ acts 
on the set $\{x_1, \dots, x_n\}$ by permuting the indices, and 
we consider the induced action on $A.$ Let $A^i$ be the subalgebra 
of polynomials invariant under the action of the permutation $s_i=(i,i+1).$ 
We make $A$ and $A^i$ graded by giving each $x_j$ degree $1,$ 
and consider the following complexes $R_i$ and $R_i'$ of graded $A$-bimodules:
\begin{eqnarray*} 
R_i: && 0 \lra A\otimes_{A^i} A\stackrel{m}{\lra} A \lra 0 \lra 0 , \\
R_i': &\quad& 0 \lra 0 \lra A \stackrel{\eta_i}{\lra} 
A\otimes_{A^i} A\{ -1\} \lra 0 , 
\end{eqnarray*} 
with $A$ in cohomological degree $0$ in both complexes, $m$ the 
multiplication map, and 
$$\eta_i(1) = (x_i -x_{i+1})\otimes 1 + 1 \otimes (x_i - x_{i+1}).$$ 
Rouquier \cite{Ro1} introduced these complexes (for an arbitrary Coxeter group) 
and showed that they induce 
a genuine braid group action on the homotopy category of complexes of 
graded $A$-modules. 

We now extend Rouquier's construction to an invariant of braid cobordisms.
To the two movies of a positive double point (see figure~\ref{picasso9})
we assign the homomorphisms 
\begin{equation} 
\begin{CD} 
 0    @>>>  0 @>>>  A_n @>>> 0  \\
   @VVV   @VVV   @VV{\mathrm{id}}V  @VVV   \\
  0   @>>>  A \otimes_{A^i} A @>{m}>> A @>>> 0  
\end{CD}  
\end{equation} 
and 
\begin{equation} 
\begin{CD} 
 0   @>>> A @>{\eta_i}>>  A\otimes_{A^i} A @>>> 0    \\ 
      @VVV          @VV{\mathrm{id}}V                   @VVV  @VVV   \\
 0    @>>>  A_n @>>>  0 @>>> 0  
\end{CD}  
\end{equation} 
of the corresponding bimodule complexes. To a negative double point 
we should assign a homomorphism from $A$ to (possibly shifted) $R_i'.$ 
Since $\eta_i$ is injective, the only bimodule homomorphism from $A$ 
to $R_i'$ is zero, while any homomorphism from $A$ to $R_i'[-1]$ is 
null-homotopic. Thus, we are forced to assign the zero homomorphism of 
complexes to any braid cobordism which contains a negative branch point. 

\begin{prop} The above assignment is a semi-trivial categorification 
of a braid group representation. 
\end{prop} 

\begin{proof}
Indeed, we know from \cite{Ro1} that the complexes $R_i, R_i'$ produce a
braid group action. Arguments analogous to the ones in the first half of
the proof 
of Theorem~\ref{first-theorem} imply that the braid movie move relations
hold and our data extends to a projective monoidal functor from $\mc{BC}_n^c$
to $\mathrm{Fun}(\mc{C}(A)),$ where $\mc{C}(A)$ is the homotopy 
category of complexes of graded $A$-modules. If we restrict to bounded complexes
of finitely-generated graded $A$-modules, the Grothendieck group will 
be a rank one free $\Z[q,q^{-1}]$-module with generator $[A],$ and 
$\sigma_i$ will act by $\sigma_i [A] = - q [A].$
\end{proof}

\begin{rmk} 
The categorification is semi-trivial since any negative braid cobordism has
the trivial invariant. Such terminology doesn't do justice to this  
categorification, since the bimodules in $R_i$ and their tensor products carry an enormous amount of information, can be used to 
reconstruct a regular block of the category $\mc{O}$ for $\mf{sl}_n,$ see \cite{So},    
and serve as building blocks for a triply-graded link homology theory \cite{Kh5}. \end{rmk}

\begin{rmk} Since we can work over 
$\Z$ rather than over $k$, the projectivity is at most in the minus sign. \end{rmk}

The above example of a braid group action led Rouquier to define the 
notion of a braid group categorification (working in the context of 
Coxeter braid groups; we specialize his construction to the case of the usual 
braid group $\mathrm{Br}_n$). For each braid $\sigma$ choose 
a representative $R(\sigma)$ in the isomorphism class of complexes 
of graded $A$-bimodules assigned to $\sigma.$ Rouquier considers a 
monoidal category $\mc{B}_{S_n}$ with objects $\sigma$ and graded 
bimodule homomorphisms $\mathrm{Hom}(R(\sigma),R(\tau))$ as 
morphisms from $\sigma$ to $\tau.$ He then defines a braid group 
categorification as a monoidal functor from $\mc{B}_{S_n}$ to 
$\mathrm{Fun}(\mc{C}),$ for a triangulated category $\mc{C}.$ 

Recall that our definition of a categorification of a braid group representation 
requires a (projective) monoidal functor from $\mc{BC}_n$ to 
$\mathrm{Fun}(\mc{C}).$ In both categories, $\mc{BC}_n$ and 
$\mc{B}_{S_n},$ isomorphism classes of objects are in bijection 
with braids. For $\mc{B}_{S_n}$ this is a nontrivial result, which says 
that $R(\sigma) \cong R(\tau)$ implies $\sigma=\tau,$ see \cite[Remark 3.9]{Ro1}. 
The morphisms in the two categories are quite different, though. In 
$\mc{BC}_n$ they are braid cobordisms, a purely topological construct, 
while in $\mc{B}_{S_n}$ they are homomorphisms of certain bimodules. 

A few paragraphs above we defined a (projective) monoidal functor 
from $\mc{BC}_n$ to  $\mc{B}_{S_n},$ trivial on any cobordism which 
contains a negative branch point. Using this functor we can convert 
any braid group categorification in the sense of Rouquier into 
a semi-trivial categorification of a braid group representation in our sense.
This conversion is the closest relation that we were able to find between
Rouquier's notion of categorification of braid group representations and ours.

\subsection{Categorification in Fukaya-Floer categories} \label{FukF}

\subsubsection*{A semi-trivial action of the braid cobordism category}

Braid group Br$_n$ actions arise naturally in symplectic geometry as groups
of symplectomorphisms generated by Dehn twists about $A_{n-1}$-chains of Lagrangian spheres \cite[Appendix]{Se2}. Following \cite{KS} we describe the local model as monodromy on Milnor fibers of $k$-dimensional $A_{n-1}$ singularities (for any
$k$). Passing to the derived Fukaya-Floer category one gets
braid group actions on triangulated categories, categorifying the usual monodromy
representations on (co)homology. We then describe an extension
to a semi-trivial representation of the braid cobordism category.
Throughout we are assuming that all the technicalities involved in the existence of the (derived) Fukaya-Floer category are overcome; since our local examples are \emph{exact} symplectic manifolds with $c_1=0$, Paul Seidel's
forthcoming book \cite{Se6} justifies this carefully (and allows us to use
Floer cohomology consisting of finite dimensional graded $\C$-vector spaces, rather than relatively or periodically graded modules over Novikov rings), but there are probably parts of the discussion below that available technical results do not yet quite fully justify.

Thinking of $\mathrm{Br}_n$ as the fundamental group of the configuration space $C_n$ of $n$ \emph{unordered} points $\mathbf{p}=\{p_j\}_{j=1}^n$ in $\C$, we can associate a $k$-dimensional symplectic manifold to every point
$\mathbf{p}\in C_n$:
\begin{equation} \label{Xp}
X_\mathbf{p}=\left\{(x_1,\ldots,x_k,t)\in\C^k\times\C\colon \sum_{i=1}^kx_i^2=
\prod_{j=1}^n(t-p_j)\right\}.
\end{equation}
We give $X_\mathbf{p}$ the symplectic structure $\omega$ that it inherits from $\C^{k+1}$.
Via the $t$ variable $X_\mathbf{p}$ fibers over $\C$ with fibers affine
quadrics $\sum_{i=1}^kx_i^2=a$ (which are symplectically just $T^*S^{k-1}$s) except at the critical points $t=p_j$ where the $S^{k-1}$s
collapse to points in the singular fibers $\sum_{i=1}^kx_i^2=0$.

Therefore for every loop in $C_n$ we get an $S^1$-family of such symplectic manifolds with a fiberwise symplectic closed 2-form on the total space.
There is a natural connection on the bundle whose horizontal subspaces
are the orthogonals, with respect to the 2-form, to the fibers; see \cite{Se4} or \cite{Se5}. The induced symplectic parallel transport along the fibration
gives us monodromy in the symplectomorphism group Aut$(X_\mathbf{p},\omega)$.
(There is a technical issue here due to the noncompactness of the fibers;
this is resolved in a related example in \cite[Section 4A]{SS}.)
Isotopic loops give hamiltonian isotopic symplectomorphisms (since the
curvature of the connection is hamiltonian) producing a homomorphism
$$
\mathrm{Br}_n=\pi_1(C_n)\to\pi_0(\mathrm{Aut}(X_\mathbf{p},\omega)).
$$
In particular the positive generators of $\mathrm{Br}_n$ can be described
as follows. Order the points $p_j\in\C$ and pick paths in $\C$ between consecutive
points, which do not intersect in their interiors. Over these paths lie canonical Lagrangian
spheres $S^k\subset X_\mathbf{p} $ (fibered by the $S^{k-1}$s in the $T^*S^{k-1}$ fibers of
$X_\mathbf{p}$). Then rotating two adjacent points anticlockwise through 180 degrees about the path between them gives rise to the monodromy transformation that is the Dehn twist about the associated Lagrangian sphere. This can be
thought of as monodromy around the singular variety one gets by collapsing
this Lagrangian sphere by moving the two points together. (In other words
we can fill in the $S^1$-family to a family over the disc $D^2$ whose central
fiber has an ordinary double point given by setting two $p_j$s to be equal
in (\ref{Xp}).)

Roughly speaking, symplectomorphisms of $X_\mathbf{p}$ act on the Fukaya-Floer category of $X_\mathbf{p}$ since the latter is a symplectic invariant
constructed out of the Lagrangian submanifolds of $X_\mathbf{p}$ and their Floer cohomologies. More precisely, since $c_1(X_\mathbf{p})=0$, there is an $A^\infty$-category whose objects are those \emph{graded} Lagrangian
submanifolds \cite{Se3} (carrying flat unitary bundles) whose Floer cohomology is well-defined \cite{FO3}.
Seidel also shows that the above Dehn twists furnish \emph{graded}
symplectomorphisms \cite{Se3}, and so act (as equivalences) on the derived Fukaya-Floer
category $D^b(Fuk(X_\mathbf{p}))$ of $X_\mathbf{p}$. Hamiltonian isotopies between symplectomorphisms induce natural transformations between their corresponding
functors on $D^b(Fuk(X_\mathbf{p}))$ (see below), so we end up with an induced
representation
$$
\mathrm{Br}_n=\pi_1(C_n)\to\pi_0(\mathrm{Aut}(X_\mathbf{p},\omega))
\to\mathrm{Auteq}(D^b(Fuk(X_\mathbf{p}))).
$$
In particular, the action of the Dehn twist $T_L$ about the Lagrangian sphere $L$ (the twist being independent of the grading on $L$) on the object $[L']\in D^b(Fuk(X_\mathbf{p}))$
can be described as follows \cite{Se4}. Think of the identity in
$HF^*(L,L')^*\otimes HF^*(L,L')$ as a morphism from $HF^*(L,L')\otimes[L]$ to $[L']$ in $D^b(Fuk(X_\mathbf{p}))$, and take its cone:
\begin{equation} \label{TL}
T_L[L']=\mathrm{Cone}\big(HF^*(L,L')\otimes L\to L'\big).
\end{equation}
(Cf. (\ref{tri1}), which is what this should correspond to under mirror symmetry.)

We explain how this can be extended to a semi-trivial representation of the braid cobordism
category. For any braid cobordism with a negative braid move we
set to zero the natural transformation between the two corresponding functors
on $D^b(Fuk(X_\mathbf{p}))$. Braid cobordisms with no positive or negative
moves (i.e. isotopies between two braids) induce natural transformations between the two functors associated to those braids (as claimed above)
as follows. The isotopy gives a family of symplectic manifolds $X_\mathbf{p}$
over the annulus, with the two boundary $S^1$-families being those associated
to the two braids. Call their monodromies $\psi$ and $\phi$. The family carries
a closed 2-form which is fiberwise symplectic; we add on a large multiple
of the pullback of a symplectic form on the annulus to produce a symplectic
total space. Then Seidel \cite{Se5} defines what he calls the relative Gromov invariant of this family, counting pseudoholomorphic sections (in an appropriate compatible almost complex structure) interpolating between horizontal sections
over each boundary (i.e. generators of the Floer cochain complexes of the
boundary symplectomorphisms $\psi$ and $\phi$) to produce a canonical element of the Floer cohomology $HF^*(\psi^{-1}\circ\phi)$
of the symplectomorphism $\psi^{-1}\circ\phi$. This element provides the
required natural transformation by pairing it with $1\in HF^0(\psi(L),\psi(L))$
using the cup product in Floer cohomology:
$$
HF^*(\psi(L),\psi(L))\otimes HF^*(\phi\circ\psi^{-1})\to
HF^*(\psi(L),\phi(L)).
$$
(This is best described, as Ivan Smith pointed out to us, on the product
$(X_\mathbf{p},\omega)\times(X_\mathbf{p},-\omega)$, where the above groups become $HF^*(\psi(L_1)\times\psi(L_2),\Delta_{X_\mathbf{p}})$,
$HF^*(\Delta_{X_\mathbf{p}},\Gamma_{\phi\circ\psi^{-1}})$, and
$HF^*(\psi(L_1)\times\psi(L_2),\Gamma_{\phi\circ\psi^{-1}})$ respectively, and the cup product becomes the ordinary cup product on Lagrangian Floer cohomology defined by holomorphic triangles. Here $\Delta_{X_\mathbf{p}}
\subset X_\mathbf{p}\times X_\mathbf{p}$ is the diagonal Lagrangian, and
$\Gamma_{\phi\circ\psi^{-1}}$ is the graph of $\phi\circ\psi^{-1}$.)

Given a braid cobordism which is a composition
of braid isotopies and positive braid moves, there is a corresponding fibration
over the annulus, the fibers of which are either manifolds
or have complex ordinary double points. This is because the family (\ref{Xp})
over $C_n$ extends by the same formulae over the partial compactification
$\overline{C}_n$ where pairs of points are allowed to collide, over which $X_p$ acquires an ordinary double point. The local model (\ref{z2z}) for
positive braid moves corresponds \cite{Ru} locally to a holomorphic map from the base (the $z_1$ axis in (\ref{z2z})) to $\overline C_n$ that intersects the double point locus transversally and positively. So for an appropriate
choice of oriented complex structure on our annulus the pullback family (\ref{Xp}), with the closed 2-form pulled back from $\C^{k+1}$, is a Lefschetz
fibration over the annulus in the sense of \cite[Definition 1.5]{Se5}.

Adding a large oriented symplectic form pulled back from the annulus, we
get a Lefschetz fibered symplectic manifold. Using pseudoholomorphic sections of this Lefschetz fibration with an appropriate almost complex structure, Seidel extends the Gromov invariant to such fibrations
\cite[Section 3]{Se5} to give a chain map between the Floer cochain complexes
of the two boundary symplectomorphisms. This is our natural transformation.

For instance, in the case of the braid cobordism with one positive braid move from the identity to the Dehn twist about a Lagrangian sphere $L$, we get the natural transformation from the identity to $T_L$ that comes from mapping
$[L']$ into the second term of (\ref{TL}) (inducing, for instance, the
horizontal map of \cite[Theorem 1]{Se5}). This is the analogue for Fukaya
categories of the positive natural transformations (\ref{moves}).

Since all of the positive braid movie moves of Section \ref{sec-bc1} arise from isotopies of maps of the annulus into $\overline{C}_n$ \emph{through maps which are orientedly transverse to the double point locus}, they give rise to
isomorphic natural transformations. This is because in Seidel's theory, two
homotopies between symplectomorphisms that are themselves homotopic induce
homotopic maps on Floer cochain complexes.

We do not currently know a natural way to extend this nontrivially to negative branch points, as the geometry involved would be locally anti-holomorphic.
Seidel and Smith have some ideas about how one might proceed, however.

\subsubsection*{The Fukaya-Floer category of $T^*Fl$}

In discussions with Paul Seidel in 1998 about the braid group action on the derived category of coherent sheaves on the cotangent bundle of the flag variety of Section \ref{flagsection}, he explained what should be the translation
of the action under mirror symmetry. The $N_i$ of Section \ref{flagsection} are coisotropic submanifolds (in the standard symplectic structure on $T^*Fl$) whose characteristic foliation is a
fibration by isotropic $S^2$s; therefore one can do generalised
Dehn twists about them (\cite[Section 1.4]{Se1}, \cite[Section 2.3]{Pe}), yielding symplectomorphisms $T_i$. Seidel conjectures that these satisfy the braid relations, and so give a representation of the braid group on the Fukaya category of $K$.

One should be able to extend this to a semi-trivial action of the braid cobordism
category in the same way as described above for the manifolds $X_{\mathbf
p}$. In place of the ordinary double point degenerations of the $X_{\mathbf
p}$ one would use normal crossings degenerations of $T^*Fl$ in which the $N_i$ collapse to the bases $K_i$ of their characteristic $S^2$ foliations. This would require a family version of Seidel's theory \cite{Se5}.

\subsubsection*{The examples of Seidel-Smith and Manolescu}

There are other naturally occurring families of symplectic manifolds fibering over $C_n$ whose Fukaya-Floer categories are related to knot theory and link homology \cite{Kh1, Kh3}. Seidel and Smith \cite{SS} consider a certain
nilpotent slice of the set of $\mathfrak{sl}_{2n}$ matrices, with the fibration
to $C_{2n}^{\,0}$ (configurations with center of mass zero) given by the set of eigenvalues of the matrix. The symplectic monodromy of the fibers
gives a representation of the braid group in the group of symplectomorphisms
up to isotopy, and the Seidel-Smith braid invariant is the Floer cohomology of a certain Lagrangian with its image under this representation; they then
show that this is in fact a link invariant.\

This example is closely related to the example of the Fukaya-Floer category
of $T^*Fl$ above,
by the deformation of adjoint orbits into the cotangent bundle of
the flag variety (by taking the usual simultaneous resolution of the degeneration
to the nilpotent cone). Manolescu \cite{Ma1} shows that it is also closely related to the example (\ref{Xp}), as it sits inside the Hilbert scheme of points
on the two dimensional version of (\ref{Xp}). (In this picture, Seidel-Smith's
Lagrangian is the product of $n$ of the Lagrangian $S^2$s described above,
fibered over $n$ non-intersecting arcs joining disjoint pairs of points in
$C_{2n}$.) In \cite{Ma2} he also produces other more complicated manifolds fibering over configuration space giving rise to more
complicated link invariants by the same procedure.

In all of these examples the fibrations extend to fibrations with singular
fibers over the partial compactification
$\overline{C}_n$ in which pairs of points are allowed to coincide. Instead
of a Lefschetz fibration with vanishing cycles of Lagrangian spheres collapsing over the discriminant locus to ordinary double points, the example of \cite{SS} exhibits a family version of this behaviour. That is, the singular fibers
contain a locus $S$ of singularities, and are locally holomorphically modelled on $S$ times by a surface ordinary double point. The smooth fibers contain coisotropic families of isotropic vanishing $S^2$s fibering over the copy of $S$ to which they collapse under
parallel transport to the singular fiber. Something similar happens in \cite{Ma1}
(with some extra singularities thrown in for good measure) and in \cite{Ma2},
with $S^2$s replaced by $\PP^k$s and the surface ordinary double points (locally
the contraction of $\PP^1\subset T^*\PP^1$) replaced by singularities modelled on the contraction of $\PP^k\subset T^*\PP^k$.

Thus an extension of Seidel's theory \cite{Se5} from Lefschetz fibrations
to fibrations with these fibered singularities would give a family
Gromov invariant defining a natural transformation between the symplectic
monodromy of a braid (and its associated functor on the Fukaya-Floer category)
and that of the braid composed with a simple positive braid cobordism. To
reiterate: such a positive braid move corresponds to a map of the annulus
into $\overline{C}_n$ with one transverse positively oriented intersection
with the double point locus of $\overline{C}_n$. Thus the map can be
made holomorphic at this point by an isotopy through maps with the same transversality
to the double point locus. Pulling back the family over $\overline{C}_n$
and pulling back a large positively oriented symplectic form from the base
we get a symplectic manifold fibering over the annulus with one singular
fiber with singularity of a fixed standard holomorphic fibered model (compatible with the symplectic form, which is K\"ahler in a neighbourhood of the singular
locus). Denote the symplectic monodromies around the two boundary circles
of the annulus by $\psi$ and $\phi$. Then an extension of \cite{Se5} to this setting would provide us with an element of $HF^*(\phi\circ\psi^{-1})$
and a natural transformation between the actions of $\phi$ and $\psi$ on the derived Fukaya-Floer categories of the fiber.

\subsection{More on polarization}

In our examples of categorifications there tend to 
exist big discrepancies between the invariants $F(S)$ of positive and negative
braid cobordisms $S.$ Often, the invariant is nontrivial on any 
positive braid cobordism but trivial on a negative braid cobordism with 
just a few branch points. The invariant described in Theorem~\ref{first-theorem}
serves as a model example. Recall that to any braid cobordism $S$ between
braid words $g$ and $h$ we assigned a homomorphism $\pm F(S)$ 
of $A_n$-bimodule complexes $R(g)$ and $R(h).$ 

\begin{prop} \label{pr-notzero} $F(S)$ is nonzero in the homotopy category for any 
positive braid cobordism $S.$ 
\end{prop} 

\begin{proof} Composing $S$ with the identity cobordism from $g^{-1}$ to
itself, we can assume that $g=1$ and $F(S)$ is a homomorphism 
from $A_n$ to the complex $R(h).$ The latter complex consists of 
a single bimodule $A_n$ and the sum of terms isomorphic to 
$P_i \otimes {\hspace{0.02in}} _jP,$ for various $i$ and $j.$ 
Thus, $R(h)$ decomposes, $R(h)\cong A_n \oplus Y(h),$ as a 
bimodule, not as a complex, where $Y$ is a direct sum of terms
 $P_i \otimes {\hspace{0.02in}} _jP.$

Any composition of bimodule maps 
$$ A_n \lra P_i \otimes {\hspace{0.02in}} _jP \lra A_n$$ 
takes $1\in A_n$ to the ideal $I$ of $A_n$ generated by $X_1, \dots, X_n.$
Hence, given any chain homotopy $\alpha$ from $A_n$ to $R(h),$ 
$d\alpha + \alpha d$ will take $1\in A_n$ to an element of the 
ideal $I$ plus an element of $Y(h).$ 

The map $F(S)$ is the composition of maps $F(S_k): R(g_k) \lra R(g_{k+1})$
 corresponding to the Reidemeister moves and positive branch points. Here
$1=g_0, g_1, \dots, g_m=h$ is the sequence of braid words in some movie 
of $S.$ Each of these 
maps $F(S_k)$ takes $1\in A_n\subset R(g_k)$ to $\pm 1 + y_k \in A_n\oplus Y(g_{k+1})$ 
(this follows by a simple case-by-case verification). Therefore,  
$F(S)$ takes $1$ to $\pm 1 + y \in A_n \oplus Y(h)$ and cannot 
be null-homotopic, since $\pm 1$ does not lie in the ideal $I.$
\end{proof} 

This argument fails if $S$ has a negative branch point, since 
then the corresponding map takes $1$ to $\pm (X_{i-1}-X_{i+1}),$ 
which \emph{does} belong to $I,$ plus an element of $Y(g_{k+1}).$ 
In fact, $F(S)=0$ if $S$ is the cobordism from the trivial braid 
to the braid word $\sigma_i^{-2}$ which creates two adjacent  
negative branch points, for the simple reason that $(X_{i-1}-X_{i+1})^2=0$
in $A_n.$ We don't know a necessary and sufficient topological 
condition for $F(S)$ to be zero. 

An extreme version of the above positive-negative braid 
cobordism discrepancy is manifested in semi-trivial categorifications, 
sometimes the only ones available to extend a braid group action. 
For instance, in the categorification via Rouquier complexes, 
$F(S)=0$ for any cobordism with a negative branch point, but 
the argument in the proof of Proposition~\ref{pr-notzero} seems 
to have its counterpart for Rouquier complexes, showing that 
$F(S)\not= 0$ for any positive braid cobordism. 

A similar positive-negative imbalance exists for the invariants of braid cobordisms coming from the rings $H^n$ and from complexes of matrix factorizations (examples 
{\bf 2} and {\bf 6} in Section~\ref{sec-survey}).  

Positive braid cobordisms can be realized via complex curves \cite{Ru} 
and negative cobordisms via antiholomorphic curves in $\C^2.$ We believe that the positive-negative asymmetry in our examples will 
ultimately prove to be of the same nature as 
in Donaldson theory, where his invariants exhibit a markedly 
different behaviour on algebraic surfaces (where they are always nontrivial) and on their antiholomorphic counterparts. The same notions
of positivity extend from complex to symplectic geometry, and conjecturally
Donaldson's invariants never vanish on symplectic manifolds (Taubes has proved
this for the Seiberg-Witten invariants). One way that positivity arises in symplectic geometry is in the monodromy
of Donaldson's symplectic Lefschetz fibrations involving only positive Dehn
twists. (This is directly related to the nonvanishing of Donaldson
invariants, or at least Seiberg-Witten invariants, by the work of Donaldson
and Smith on pseudoholomorphic multisections of such fibrations.) This manifests itself in that it was relatively easy to define positive braid moves on
Fukaya-Floer categories in Section \ref{FukF}, but harder
to define non-trivial representations of the negative moves.

\subsection{When Definition~\ref{d-one} fails}

There are several known and conjectural 
examples where a braid group (or a Coxeter braid group) 
acts naturally on a triangulated category $\mc{C}$ which decomposes nontrivially
into a direct product of its indecomposable pieces, 
$\mc{C}= \oplusop{\mu\in J} \mc{C}_\mu,$ for $J$ being some index set;
for instance, the set of weights of an irreducible representations of a simple
Lie algebra. The braid group action induces a nontrivial permutation 
homomorphism $\mathrm{Br}_n \lra S_J$ to the group of permutations of 
the set $J.$ Definition~\ref{d-one} is ill-suited for such actions, since
if 
$F_{\sigma_i} (\mc{C}_{\mu}) \subset \mc{C}_{\mu'}$ for $\mu'\not= \mu,$
the only natural transformation from the identity functor on $\mc{C}_{\mu}$
to $F_{\sigma_i}$ is zero. We don't know any good substitute for
Definition~\ref{d-one} 
which would deal with braid group actions of these type. Instead, we list
several such actions below. 

\begin{itemize} 
\item The braid group action on the derived category of sheaves on disjoint
unions of partial flag varieties in Section~\ref{sec-survey} example {\bf 7}. 
\item Examples {\bf 5} and {\bf 6} in Section~\ref{sec-survey} can be extended by looking at oriented tangles with ``out" orientation at 
$k_1$ endpoints and ``in" orientation at $k_2$ endpoints. 
The braid group on $k_1+k_2$ strands acts by braiding tangles along 
their endpoints. Triangulated categories then are
direct sums of categories, one for each lineup of in and out orientations,
and the braid group action induces a permutation action of the symmetric
group $S_{k_1+k_2}$ on the set of these lineups. 
The braid group action extends to an invariant of 
oriented tangles and tangle cobordisms, and a substitute for
Definition~\ref{d-one} 
in this case can presumably be found by restricting to a suitable subcategory of braid cobordisms.  
\item In a remarkable paper \cite{CR}, Chuang and Rouquier gave the 
definition and deep analysis of $\mf{sl}_2$-categorification and, 
as a result of their work, introduced certain derived equivalences $\Theta$
between blocks of cyclotomic Hecke algebras, symmetric and general linear
groups over finite fields, category $\mc{O},$ etc. It is natural to conjecture
that these equivalences organize into braid group actions on direct sums
of suitable blocks. In the cyclotomic Hecke algebra case these actions should
categorify braid group actions on irreducible $U_q(\mf{sl}_n)$-representations
and affine braid group actions on integrable $U_q(\widehat{\mf{sl}}_n)$-representations
(if the categories are not graded, set $q=1$), enriching Ariki's categorification
\cite{A} of these representations in the $q=1$ case. Derived equivalences similar to the 
ones in \cite{CR} should exists in derived categories of coherent sheaves on Nakajima varieties \cite{Nk}, giving rise to Coxeter braid group actions on these categories. 
\end{itemize}

\section{An affine braid group action on $D^b(T^*Fl)$}
\label{flagsection}

\subsection{Notation and statement of results}

Fix an $n$-dimensional complex vector space $V$, and denote by $Fl$ its full flag variety
$$
Fl=Fl(V)=\{0=V_0<V_1<\cdots<V_{n-1}<V_n=V\}.
$$
Let $\pi_i\colon Fl\to Fl_i$ denote the quotient given by forgetting the $i$th flag, so $Fl_i$ is the variety of nested subspaces $(V_j)_{j\ne i}$ with $\dim V_j=j$. $Fl$ therefore
carries the tautological bundles $V_j$, all of them except $V_i$ being the
pullbacks of tautological bundles (also denoted $V_j$) on $Fl_i$.
We let $L_i=\Lambda^iV_i^*$ denote the determinants of their duals.

We work on $K\define T^*Fl$ with the divisors $N_i\define\pi_i^*T^*Fl_i\subset K$ and their induced projections $p_i$ to $T^*Fl_i$. Letting $K_i$ denote $T^*Fl_i$, we get diagrams
\begin{equation} \label{piota}
\begin{array}{l}
\ N_i\stackrel{\iota\_i}{\into}K\\
\ \ \downarrow^{p_i} \\
\ K_i.
\end{array}
\end{equation}
Suppressing some obvious pullback maps for clarity, $N_i$ is the zero locus of the canonical section of the fiberwise cotangent
sheaf $T^*_{\pi_i}=\omega_{\pi_i}$ given by projecting the tautological section
of $T^*Fl$ on $K=T^*Fl$ to $T^*_{\pi_i}$. We always use $\omega$ to denote
(relative) canonical bundles. Since $Fl$ is the projective bundle
$\PP(V_{i+1}/V_{i-1})$ over $Fl_i$, we see that $T_{\pi_i}\cong
\hom(V_i/V_{i-1},V_{i+1}/V_i)$. Thus
\begin{equation} \label{omega}
\mc{O}(N_i)\cong\omega_{\pi_i}\cong\omega_{p_i}\cong L_i^{-2}\otimes L_{i-1}\otimes L_{i+1}.
\end{equation}

Define the exact functors
$$
a_i\define\iota\_{i*}p_i^*\colon D^b(K_i)\to D^b(K),
$$
and their right adjoints
$$
b_i\define p_{i*}\iota_i^!\colon D^b(K)\to D^b(K_i).
$$
Here $\iota_i^!=\mc{O}(N_i)\otimes\iota_i^*\,[-1]$
is the right adjoint of $\iota\_{i*}$
by Serre duality, and is such that $\iota\_{i*}\iota_i^!\cong\hom_{\mc{O}_K}(\mc{O}_{N_i},\ \cdot\ )$.

The right adjoint of $p_{i*}$ is $\omega_{p_i}\otimes p_i^*[1]$ by Serre duality,
and the right adjoint of $\iota_i^*$ is $\iota_{i*}$, so that of $\iota_i^!$
is $\iota_{i*}(-N_i)[1]$. Therefore the right adjoint of $b_i$ is
$\iota_{i*}(-N_i)[1](\omega_{p_i}\otimes p_i^*)[1]=\iota_{i*}p_i^*[2]=a_i[2]$ by
(\ref{omega}).

That is we have adjunctions $a_i\dashv b_i\dashv a_i[2]$, and so, in particular, the counit $\ev_i\colon a_ib_i\to\id$ and unit $\ev_i'\colon\id[-2]\to a_ib_i$.
We define functors
\begin{equation} \label{ui}
U_i\define a_ib_i\colon D^b(K)\to D^b(K), \qquad i=1,\ldots,n-1,
\end{equation}
and
\begin{equation} \label{ti}
T_i\define\text{Cone}\big(U_i\Rt{\ev_i}\id\big), \qquad
T_i'\define\text{Cone}\big(\id\Rt{\ev_i'}U_i[2]\big)[-1].
\end{equation}

In Section \ref{YB} we will prove the following.

\begin{theorem}
$T_iT_i'\cong\id\cong T_i'T_i$, and the $T_i$ satisfy the braid relations:
\begin{itemize}
\item $T_iT_j\cong T_jT_i$ \qquad\,for $|i-j|>1$, and
\item $T_iT_jT_i\cong T_jT_iT_j$ \ for $|i-j|=1$.
\end{itemize}
\end{theorem}

Then in Section \ref{affine} we will define an extra invertible functor
by
\begin{equation} \label{Tdef}
T\define(\L_1T_1T_2\ldots T_{n-2})T_{n-1}(T_{n-2}'\ldots T_2'T_1'\L_1^{-1}),
\end{equation}
where $\L_1$ denotes the functor $L_1\otimes(\ \cdot\ )$. Calling this $T_0=T=T_n$ in a cyclic ordering, it extends
the above braid group action to one of the affine braid group:

\begin{theorem} The functor $T$ commutes with $T_i$ for $2\le i\le n-2$, and braids with $T_1$ and $T_{n-1}$:
$$
T_0T_1T_0\cong T_1T_0T_1 \qquad\text{and}\qquad T_nT_{n-1}T_n\cong T_{n-1}T_nT_{n-1}.
$$
\end{theorem}

In fact we use Fourier-Mukai transforms (Section \ref{FMTs}) rather than functors, giving a representation of the braid group in the group of invertible
Fourier-Mukai transforms. Applying the Fourier-Mukai functor gives the above
(slightly weaker) results.

\subsection{Fourier-Mukai transforms} \label{FMTs}

In places, using Fourier-Mukai kernels makes this part of the paper rather a triumph of notation over clarity. For instance, the rest of this section derives the Fourier-Mukai kernels for
the above functors, and is long and formal. But it can be safely skipped: the later sections show that composition of these kernels satisfies the relations
of the braid group. This will satisfy most readers, without seeing a proof that the corresponding braid group of Fourier-Mukai functors is the one described above.

We will often suppress
pullback maps applied to line bundles, and pushforwards applied to structure sheaves when the map is an embedding: given $\iota\colon D\into K$ we often denote $\iota_*\mc{O}_D$ by $\mc{O}_{i(D)}$ or $\mc{O}_D$ to save on notation.

If $D_1$
and $D_2$ are divisors on spaces $X_1$ and $X_2$ respectively, we use the
notation $\mc{O}(D_1,D_2)$ for the line bundle $\mc{O}(D_1)\boxtimes\mc{O}(D_2)$ on any
fiber product $X_1\times_BX_2$ of $X_1$ and $X_2$. When $X_1=X=X_2$ we let
$\Delta_X$ (or sometimes just $\Delta$) denote the diagonal copy of $X$ in $X\times_BX$.

In several places we will have $A,B\subset X$ whose intersection $A\cap B$ is a Cartier divisor in $A$ with $\mc{O}_A(A\cap B)\cong L|_A$ for some line bundle $L$ over $A\cup B$. Then we will repeatedly use the connecting
homomorphism $\mc{O}_B(L)\to\mc{O}_A[1]$ of the standard exact sequence
$$
0\to\mc{O}_A\to\mc{O}_{A\cup B}(L)\to\mc{O}_B(L)\to0,
$$
i.e. the morphism representing the $\Ext^1$ class of the above extension.

$D^b(X)$ denotes the bounded derived category of coherent sheaves
on a smooth quasi-projective variety $X$, and all functors
between such categories will be the derived functors (though we omit the $\mathbf L$s and $\mathbf R$s). We denote the derived dual $\hom(\F,\mc{O})$
of an object $\F$ by $\F^\vee$.

Notoriously, cones are in general not functorial in derived categories because given an arrow one must pick a representative of its quasi-isomorphism class to take the cone on. For instance, working with the homotopy category of complexes of quasi-coherent injective sheaves with coherent cohomology, we have to pick a map of complexes in a homotopy class before taking its cone, and this choice is not functorial.

For us, however, just as in \cite{ST}, there is a functorial choice of map
since all of our maps are evaluation (and co-evaluation) maps. Thus
all of our cones will be chosen functorially.

We recommend \cite{Hu} for the theory of Fourier-Mukai transforms. A Fourier-Mukai transform is an exact functor $D^b(X_1)\to D^b(X_2)$ induced by an object (Fourier-Mukai kernel) $P\in D^b(X_1\times X_2)$ on the product. (So, importantly,
we \emph{do not insist on our Fourier-Mukai transforms being invertible}.)
Using the obvious projections
\begin{equation} \label{diagram}
\diagram
X_1\times X_2 \rto^(0.6){\pi_2} \dto^{\pi_1} & X_2 \\ X_1,\!\!
\enddiagram
\end{equation}
the functor is
$$
\Phi_P\define\pi_{2*}\big(P\otimes\pi_1^*(\ \cdot\ )\big).
$$
The structure sheaf $\mc{O}_\Delta$ of the diagonal $\Delta\subset X\times X$ represents the identity on $D^b(X)$, and, more generally, given a morphism $f\colon X\to Y$, the functors
$f_*$ and $f^*$ are represented by $\mc{O}_{(\id\times f)X}\in D^b(X\times Y)$ and $\mc{O}_{(f\times\id)X}\in D^b(Y\times X)$ respectively.
It is a foundational result of Orlov that in fact all fully faithful exact functors arise in this way for projective $X_i$, though we shall not need this (our manifolds are only quasi-projective, but all the functors we use have obvious Fourier-Mukai kernels). However, this natural
transformation $\Phi$ from $D^b(X_1\times X_2)$ to Fun($D^b(X_1),D^b(X_2))$ is neither full nor faithful (see for example \cite{Ca}),
though it can be made to be so by working with dg-categories instead \cite{To}.
Again, this will not affect us; all of the natural transformations that we use are induced by obvious maps of kernels.

The composition of two Fourier-Mukai transforms $D^b(X_1)\!\Rt{\Phi_P}\!
D^b(X_2)\!\Rt{\Phi_Q}\!D^b(X_3)$ is \cite[Proposition 5.10]{Hu}
$\Phi_{Q\star P}\cong\Phi_Q\circ\Phi_P\colon D^b(X_1)\to D^b(X_3)$, where $Q\star P\in D^b(X_1\times X_3)$ is the convolution
\begin{equation} \label{conv}
Q\star P=\pi\_{13*}\big(\pi_{12}^*P\otimes \pi_{23}^*Q\big).
\end{equation}
Here $\pi\_{ij}$ is the projection from $X_1\times X_2\times X_3$ onto $X_i\times
X_j$.

The right adjoint of the functor $\Phi_P$ is
\begin{equation} \label{dual}
\Phi_{P^\vee\otimes\omega_{X_1}[-x_1]},
\end{equation}
where we are identifying $X_1\times X_2$ and
$X_2\times X_1$ in the obvious way, suppressing the pullback $\pi_1^*$ on
the canonical bundle $\omega_{X_1}$, and setting $x_i=\dim X_i$.
By (\ref{conv}) the composition $\Phi_P\circ\Phi_{P^\vee\otimes\omega_{X_1}[-x_1]}$
is represented by the kernel
\begin{equation} \label{PP}
\left(P^\vee\otimes\omega_{X_1}[-x_1]\right)\star P=\pi_{13*}\left(
\pi_{12}^*P^\vee\otimes\omega_{X_1}[-x_1]\otimes\pi_{23}^*P\right),
\end{equation}
where $\pi_{13}\colon X_2\times X_1\times X_2\to X_2\times X_2$ is the projection
to the first and third factors.

Let $\delta\colon X_1\times X_2\into X_1\times X_2\times X_1$ be the product
of $X_1$'s diagonal and the identity on $X_2$.
Then the counit $\Phi_P\circ\Phi_{P^\vee\otimes\pi_1^*\omega_{X_1}[-x_1]}
\to\id$ of the adjunction is represented by the morphism of Fourier-Mukai kernels given by the following composition. We restrict (\ref{PP}) to the image of $\delta$ before applying $\pi_{13*}$, then take
the trace map and the isomorphism $\pi_{13*}\omega_{X_1}\cong\mc{O}[-x_1]$:
\begin{multline} \label{counit}
\pi_{13*}\left(\pi_{12}^*P^\vee\otimes\omega_{X_1}[-x_1]\otimes\pi_{23}^*P\right)
\to\pi_{13*}\left(\delta_*(P^\vee\otimes P)\otimes\omega_{X_1}[-x_1]\right)\\
\Rt{\tr}\pi_{13*}\left(\mc{O}_{\delta(X_1\times X_2)}\otimes\omega_{X_1}[-x_1]\right)
\cong\mc{O}_{\Delta}.
\end{multline}
\smallskip

Using these facts we can deduce the Fourier-Mukai kernels for our functors
$U_i$ and $T_i$. The reader unconcerned with the relationship to the functors (\ref{ui}, \ref{ti}) could omit this and proceed to the next section
where the functors are forgotten and only their Fourier-Mukai kernels (\ref{UFM},
\ref{TFM}) are used to give a braid group representation in the group of invertible Fourier-Mukai kernels.

\begin{lemma} \label{UFML}
The functor $U_i$ (\ref{ui}) has Fourier-Mukai kernel
\begin{equation} \label{UFM}
U_i\define\mc{O}_{(\iota_i\times\iota_i)(N_i\times_{K_i}N_i)}(N_i,0)[-1]\in D^b(K\times K).
\end{equation}
\end{lemma}

\begin{proof}
This is proved in a number of steps, at each stage composing Fourier-Mukai functors using (\ref{conv}). \smallskip
\begin{itemize}
\item $p_{i*}\colon D^b(N_i)\to D^b(K_i)$ is represented by the
FM kernel $\mc{O}_{(\id\times p_i)N_i}\in D^b(N_i\times K_i)$.
\item $p_i^*$ is represented by $\mc{O}_{(p_i\times\id)N_i}\in
D^b(K_i\times N_i)$.
\item Therefore, the composition $p_i^*p_{i*}$ has kernel
\begin{equation} \label{pp}
\qquad \pi_{13*}\big[\mc{O}_{((\id\times p_i)N_i)\times N_i}\otimes
\mc{O}_{N_i\times(p_i\times\id)N_i}\big],
\end{equation}
where $\pi_{13}\colon N_i\times K_i\times N_i\to N_i\times N_i$ is the projection.
\item The intersection of $((\id\times p_i)N_i)\times N_i$ (of codimension
$\dim K_i$) and $N_i\times
(p_i\times\id)N_i$ (codimension $\dim K_i$) in $N_i\times K_i\times N_i$ is $N_i\times_{K_i}N_i$, which is smooth of the expected codimension $2\dim
K_i$. Therefore the (derived) tensor product of the structure sheaves equals the structure sheaf of the (transverse) intersection. 
\item Since $\pi_{13}$ restricted to $N_i\times_{K_i}N_i\subset
N_i\times K_i\times N_i$ is an embedding,
mapping it into $N_i\times N_i$ as $\id\times_{K_i\!}\id$, (\ref{pp})
shows that the kernel for $p^*_ip_{i*}$ is
$\mc{O}_{(\id\times_{K_i\!}\id)(N_i\times_{K_i}N_i)}\in D^b(N_i\times
N_i)$.
\end{itemize}

Next we compose with $\iota_i^!$:
\begin{itemize}
\item $\iota_i^!=\iota_i^*(N_i)[-1]\colon D^b(K)\to D^b(N_i)$ is represented by the Fourier-Mukai kernel $\mc{O}_{(\iota_i\times\id)N_i}
(N_i,0)[-1]\in D^b(K\times N_i)$.
% where the divisor $N_i$ that we twist by is equal to the restriction of
% both $\mc{O}(N_i,0)$ and $\mc{O}(0,N_i)$ from $K\times N_i$ to $(\iota_i\times\id)(N_i)$.
\item Therefore the composition $p^*_ip_{i*}\iota_i^!$ is represented by
$$
\qquad \pi_{13*}\big[\mc{O}_{((\iota_i\times\id)N_i)\times N_i}(N_i,0)[-1]\otimes
\mc{O}_{K\times(\id\times_{K_i\!}\id)(N_i\times_{K_i}N_i)}\big],
$$
where $\pi_{13}\colon K\times N_i\times N_i\to K\times N_i$ is the projection.
\item The intersection of $((\iota_i\times\id)N_i)\times N_i$ and $K\times
(\id\times_{K_i\!}\id)(N_i\times_{K_i}N_i)$ in $K\times N_i\times
N_i$ is
$N_i\times_{K_i}N_i$, embedded via $(a,b)\mapsto(\iota_i(a),a,b)$.
This is smooth of the expected dimension, so the (derived) tensor product of the structure sheaves equals the
structure sheaf of the (transverse) intersection. 
\item Composing with $\pi_{13}$ embeds $N_i\times_{K_i}N_i$ in $K\times
N_i$ via $\iota_i\times\id$, so $p^*_ip_{i*}\iota_i^!$ has kernel
$\mc{O}_{(\iota_i\times\id)(N_i\times_{K_i}N_i)}(N_i,0)[-1]\in D^b(K\times
N_i)$.
\end{itemize}
Finally we can compose with $\iota_{i*}$, whose kernel is $\mc{O}_{(\id\times
\iota_i)N_i}\in D^b(N_i\times K)$, in the same way. The upshot is, unsurprisingly,
that the Fourier-Mukai kernel for $\iota_{i*}p^*_ip_{i*}\iota_i^!$ is the pushforward of the kernel for $p^*_ip_{i*}\iota_i^!$ from $D^b(K\times N_i)$
to $D^b(K\times K)$, i.e. $\mc{O}_{(\iota_i\times\iota_i)(N_i\times_{K_i}N_i)}
(N_i,0)[-1]$, as claimed.
\end{proof}

There is an obvious exact sequence
$$
0\to\mc{O}_\Delta\to\mc{O}_{\Delta\,\cup\,\left((\iota_i\times\iota_i)(N_i\times_{K_i}N_i)\right)}
(N_i,0)\to\mc{O}_{(\iota_i\times\iota_i)(N_i\times_{K_i}N_i)}(N_i,0)\to0,
$$
where $\Delta\stackrel{\iota\_\Delta}\into K\times K$ is the diagonal and
the last arrow is the restriction map. This gives the connecting homomorphism
\begin{equation} \label{par}
\mc{O}_{(\iota_i\times\iota_i)(N_i\times_{K_i}N_i)}(N_i,0)[-1]\to\mc{O}_\Delta.
\end{equation}

\begin{prop}
The natural transformation $U_i\Rt{\ev}\id$ is represented by the morphism
(\ref{par}) in $D^b(K\times K)$. Thus its cone, the functor $T_i$ (\ref{ti}), has Fourier-Mukai kernel
\begin{equation} \label{TFM}
T_i\define\mc{O}_{\Delta\,\cup\,(\iota_i\times\iota_i)(N_i\times_{K_i}N_i)}
(N_i,0)\in D^b(K\times K).
\end{equation}
\end{prop}

\begin{proof}
We first claim that the counit $p^*_ip_{i*}\to\id$ of the adjunction
$p^*_i\dashv p_{i*}$ is the evaluation map $\ev_{p_i}$ induced by the map
\begin{equation} \label{punit}
\mc{O}_{(\id\times_{K_i\!}\id)(N_i\times_{K_i}N_i)}\to\mc{O}_\Delta
\end{equation}
given by restriction to the diagonal $\Delta\subset N_i\times_{K_i}N_i\subset
N_i\times N_i$. This is just a notational chase and the reader might want to take it on trust, or observe that it does the right thing when applied to any object of $D^b(N_i)$.

We prove it using adjoints and the counits of (\ref{counit}). Since $N_i$
is noncompact the formulae for composition (\ref{conv}) and counit (\ref{counit}) will only work if we work relative to $K_i$ (notice $p_i\colon N_i\to
K_i$ does have compact fibers) and use relative, instead of absolute, Serre duality. That is we replace the diagram (\ref{diagram}) by
$$
\diagram
K_i\times_{K_i}N_i \rto^(0.65){\pi_2} \dto^{\pi_1} & N_i \\ K_i,\!\!
\enddiagram
$$
with respect to which $p_i^*$ is represented by the kernel $\mc{O}_{K_i\times_{K_i}N_i}$.
Since $K_i\times_{K_i}N_i\cong N_i$ and $\pi_1$ and $\pi_2$ are
just $p_i$ and $\id$ respectively, this statement is the tautology that
$$
\pi_{2*}\left(\mc{O}_{K_i\times_{K_i}N_i}\otimes\pi_1^*(\ \cdot\ )\right)=
\id\left(p_i^*(\ \cdot\ )\right)=p_i^*.
$$
Then, by the relative version of (\ref{dual}), its right adjoint $p_{i*}$
is represented by the functor
$$
\mc{O}_{N_i\times_{K_i}K_i}^\vee\otimes\omega\_{K_i/K_i}[0]=\mc{O}_{N_i\times_{K_i}K_i},
$$
since the relative dimension of $K_i$ over $K_i$ is zero.
Using, as usual, $\pi_{ij}$ to denote the projection of
$N_i\times_{K_i}K_i\times_{K_i}N_i\to N_i\times_{K_i}N_i$
onto its $i$th and $j$th factors, by (\ref{conv}) the composition $p^*_ip_{i*}$ is represented by
$$
\pi_{13*}\left(\pi_{12}^*\mc{O}_{N_i\times_{K_i}K_i}\otimes
\pi_{23}^*\mc{O}_{K_i\times_{K_i}N_i}\right)=\pi_{13*}
\mc{O}_{N_i\times_{K_i}K_i\times_{K_i}N_i}=\mc{O}_{N_i\times_{K_i}N_i},
$$
and the counit $p^*_ip_{i*}\to\id$ is represented by the composition (cf.
(\ref{counit}))
$$
\pi_{13*}\mc{O}_{N_i\times_{K_i}K_i\times_{K_i}N_i}
\to\pi_{13*}(\delta_*\mc{O})\cong\mc{O}_{\Delta},
$$
where $\delta\colon N_i\times_{K_i}K_i\into N_i\times_{K_i}K_i\times_{K_i}N_i$
is the product of the diagonal on $N_i$ and the identity on $K_i$. The trace map is just the identity since we only are pushing down the dimension zero identity map $K_i\to K_i$. Thus the natural transformation $p^*_ip_{i*}\to\id$ is represented by the restriction map $\mc{O}_{N_i\times_{K_i}N_i}\to
\mc{O}_\Delta$. Pushing this forward from $D^b(N_i\times_{K_i}N_i)$ to $D^b(N_i\times
N_i)$ gives (\ref{punit}). \medskip

Next we have to do the same for $\iota_{i*}$ and its right adjoint $\iota_i^!$.
To handle the noncompactness we this time work relative to $K$, considering
$N_i$ as a variety over $K$ (of relative dimension $-1$). So (\ref{diagram})
becomes
$$
\diagram
N_i\times_KK \rto^(0.65){\pi_2} \dto^{\pi_1} & N_i \\ K,\!\!
\enddiagram
$$
with $\iota_{i*}$ represented by the kernel $\mc{O}_{N_i\times_KK}$. That is
$N_i\times_KK\cong N_i$ and $\pi_1,\,\pi_2$ are $\id,\,\iota_i$, respectively,
so
$$
\pi_{2*}\left(\mc{O}_{N_i\times_KK}\otimes\pi_1^*(\ \cdot\ )\right)=
\iota_{i*}(\id)=\iota_{i*}.
$$
The Fourier-Mukai kernel for its right adjoint $\iota_i^!$ is then
$\mc{O}_{N_i\times_KK}^\vee\otimes\omega_{N_i/K}[-1]=\mc{O}_{N_i}(N_i)[-1]\in
D^b(K\times_KN_i)$. Working in $D^b(K\times_KN_i\times_KK)$ in the usual
notation (\ref{conv}) gives the kernel for the composition $\iota_{i*}\iota_i^!$
as
\begin{multline*}
\pi_{13*}\big(\pi_{12}^*\mc{O}_{N_i}(N_i)[-1]\otimes\pi_{23}^*\mc{O}_{N_i}\big)
=\pi_{13*}\mc{O}_{N_i}(N_i)[-1] \\ =\mc{O}_{N_i}(N_i)[-1]\in D^b(K\times_KK)=D^b(K).
\end{multline*}
(Pushing forward $D^b(K\times_KK)\to D^b(K\times K)$ gives its more usual
absolute, rather than relative, Fourier-Mukai kernel
$\mc{O}_{(\iota_i\times\iota_i)N_i}(N_i,0)[-1]$.)

Finally by (\ref{counit}) the counit $\iota_{i*}\iota_i^!\to\id$ is represented
by the following morphism of Fourier-Mukai kernels. Firstly restriction to
the product (over $K$) of $N_i$ and the diagonal of $K$ in $K\times_KN_i\times_KK$ is the identity (since both are just $N_i$). Using the fact that $\pi_{13}$
is $\iota_i$, we are left with the trace map
$$
\mc{O}_{N_i}(N_i)[-1]\Rt{\tr}\mc{O}_K.
$$
That is, we want to know what the relative Serre duality trace morphism
$\omega_{N_i/K}[-1]\Rt{\tr^{-1}}\mc{O}_K$ is in dimension $-1=\dim N_i-\dim
K$. (Equivalently we want to describe the natural evaluation map
$\mc{O}_{N_i}^\vee\to\mc{O}_K$.) It is of course the connecting homomorphism of the triangle
\begin{equation} \label{con}
\mc{O}_{N_i}(N_i)[-1]\to\mc{O}_K\to\mc{O}_K(N_i)\to\mc{O}_{N_i}(N_i).
\end{equation}
Pushing forward from $K=K\times_KK$ to $K\times K$ we find that $\iota_{i*}\iota_i^!\to\id$
is represented by
$$
\mc{O}_{(\iota_i\times\iota_i)N_i}(N_i,0)[-1]\to\mc{O}_\Delta,
$$
with the morphism the pushforward of the connecting homomorphism (\ref{con})
by the diagonal map $\iota_\Delta\colon\Delta\to K\times K$. \medskip

We are finally ready to describe the counit $\iota_{i*}p_i^*p_{i*}\iota_i^!
\to\id$. By (\ref{punit}), $\ev_{p_i}\colon\iota_{i*}p_i^*p_{i*}\iota_i^!\to
\iota_{i*}\iota_i^!$ is represented by the restriction map
\begin{equation} \label{pev}
\mc{O}_{(\iota_i\times\iota_i)(N_i\times_{K_i}N_i)}(N_i,0)[-1]
\to\mc{O}_{(\iota_i\times\iota_i)N_i}(N_i,0)[-1].
\end{equation}
Composing with the pushforward by $\iota_\Delta$ of the connecting homomorphism
of (\ref{con}) (representing $\iota_{i*}\iota_i^!\to\id$) we obtain the counit $\iota_{i*}p_i^*p_{i*}\iota_i^!\to\id$ (by the obvious composition property of units). So this is represented by
the connecting morphism of the natural exact triangle
\begin{multline} \label{tis}
\mc{O}_{(\iota_i\times\iota_i)(N_i\times_{K_i}N_i)}(N_i,0)[-1]\to \\
\mc{O}_\Delta\to\mc{O}_{\Delta\,\cup\,(\iota_i\times\iota_i)(N_i\times_{K_i}N_i)}
(N_i,0)\to\mc{O}_{(\iota_i\times\iota_i)(N_i\times_{K_i}N_i)}(N_i,0).
\end{multline}
Thus its cone, representing $T_i$, is $\mc{O}_{\Delta\,\cup\,\left((\iota_i
\times\iota_i)(N_i\times_{K_i}N_i)\right)}(N_i,0)$.
\end{proof}

Finally, similar and equally tedious working shows that the functor
$\id\Rt{\ev_i'}U_i[2]$ is represented by the morphism
\begin{equation} \label{ev'}
\mc{O}_\Delta\to\mc{O}_{(\iota_i\times\iota_i)(N_i\times_{K_i}N_i)}(N_i,0)[1]
\end{equation}
in $D^b(K\times K)$ which is the connecting homomorphism of the exact sequence
\begin{equation} \label{ted}
0\to\mc{O}_{(\iota_i\times\iota_i)(N_i\times_{K_i}N_i)}(N_i,0)\to
\mc{O}_{\Delta\,\cup\,(\iota_i\times\iota_i)(N_i\times_{K_i}N_i)}(L_i^{-1},L_i)
\to\mc{O}_\Delta\to0,
\end{equation}
whose second arrow is restriction to $\Delta$.
Here we using the claim that $\Delta\,\cap\,(\iota_i\times\iota_i)
(N_i\times_{K_i}N_i)=\Delta_{N_i}\define(\iota_i\times\iota_i)N_i$ is a divisor in $N_i\times_{K_i}N_i$ in the class of (the restriction to
$N_i\times_{K_i}N_i$ of) $(L_i\otimes L_{i+1}^{-1}\otimes L_{i-1}^{-1},L_i)$. In turn this is $(L_i^{-1}(-N_i),L_i)$ by (\ref{omega}).

The claim is proved by noting that the pushdown of $(L_i\otimes L_{i+1}^{-1}
\otimes L_{i-1}^{-1},L_i)$ from $N_i\times_{K_i}N_i$ to $K_i$ is
\begin{equation} \label{iden}
\big(V_{i+1}\big/V_{i-1}\big)^*\otimes \big(V_{i+1}\big/V_{i-1}\big)^*
\otimes\Lambda^2\big(V_{i+1}\big/V_{i-1}\big)\,\cong\ \mc End
\big(V_{i+1}\big/V_{i-1}\big),
\end{equation}
since the pushdown of $(V_i/V_{i-1})^*=L_i\otimes L_{i-1}^{-1}$ from $N_i$
to $K_i$ is $(V_{i+1}/V_{i-1})^*$ (as everything is a pullback from the
pushdown $Fl\to Fl_i$, which is the projective bundle $\PP(V_{i+1}/V_{i-1})$
over $Fl_i$). Under this isomorphism, $\id\in\mc End(V_{i+1}/V_{i-1})$
gives a global section of (\ref{iden}) corresponding to a global section of $(L_i\otimes L_{i+1}^{-1}
\otimes L_{i-1}^{-1},L_i)$ which has zero locus
the diagonal $\Delta_{N_i}$ (this can be checked fiberwise, where it is just
the diagonal $\PP^1\subset\PP^1\times\PP^1$).

Taking cones then, $T_i'$ (\ref{ti}) has kernel
\begin{equation} \label{T'FM}
\mc{O}_{\Delta\,\cup\,(\iota_i\times\iota_i)(N_i\times_{K_i}N_i)}(L_i^{-1},L_i).
\end{equation}

\subsection{First results}

We now forget all about functors and work exclusively with the kernels
(\ref{UFM}, \ref{TFM}, \ref{T'FM}). To streamline notation we suppress
inclusion maps when they are the obvious ones, for instance writing
$\mc{O}_{N_i\times_{K_i}N_i}(N_i,0)[-1]$ for $U_i$ (\ref{UFM}) and
$\mc{O}_{\Delta\,\cup\,(N_i\times_{K_i}N_i)}(N_i,0)$ for $T_i$ (\ref{TFM}).

\begin{prop} \label{spl}
There is an exact triangle of kernels
\begin{equation} \label{UUU}
U_i\Rt{a}U_i\star U_i\Rt{b}U_i[-2]
\end{equation}
in $D^b(K\times K)$.

It is split by either of $\,\id\star\ev'_i$ or $\ev'_i
\star\id\colon U_i[-2]\to U_i\star U_i$, and similarly by either of $\,\id
\star\ev_i$ or $\ev_i\star\id\colon U_i\star U_i\to U_i$, giving isomorphisms $U_i\star U_i\cong U_i\oplus U_i[-2]$.
\end{prop}

\begin{rmk}
This is one place where using functors would have been easier and more
informative. In the notation of the introduction, the natural transformation
induced by $a$ is the unit $\id\to b_ia_i$ of $a_i\dashv b_i$, premultiplied
by $a_i$ and postmultiplied by $b_i$. A similar procedure produces $b$ from the counit $b_ia_i\to\id[-2]$ of $b_i\dashv a_i[2]$, and the statements about splitting then follow from standard facts about composing units and counits of adjoint functors.
\end{rmk}

\begin{proof}
In the usual notation (\ref{conv}),
$$
U_i\star U_i=\pi_{13*}\left(\mc{O}_{N_i\times_{K_i}N_i\times K}(N_i,0,0)[-1]\otimes
\mc{O}_{K\times N_i\times_{K_i}N_i}(0,N_i,0)[-1]\right).
$$
The intersection of the two supports is $N_i\times_{K_i}N_i\times_{K_i}N_i$,
which is smooth but has excess intersection dimension one. Therefore the
derived tensor product of the two structure sheaves has both Tor$_0$ (the structure sheaf of the intersection)
and Tor$_1$. By a standard local Koszul resolution calculation the latter is given by the structure sheaf of the intersection tensored with the excess conormal bundle, which by adjunction is
$$
\mc{O}(-N_i+\omega_{K_i},-N_i,0)\otimes\mc{O}(0,-N_i,-N_i+\omega_{K_i})\otimes
\mc{O}(N_i-\omega_{K_i},N_i,N_i-\omega_{K_i})\vspace{-3mm}
$$ \begin{equation} \label{tor1}
\hspace{85mm} \cong\mc{O}(0,-N_i,0).
\end{equation}
So we have an exact triangle
\begin{multline} \label{UU}
\pi_{13*}\mc{O}_{N_i\times_{K_i}N_i\times_{K_i}N_i}(N_i,0,0)[-1]\Rt{a}
U_i\star U_i \\
\Rt{b}\pi_{13*}\mc{O}_{N_i\times_{K_i}N_i\times_{K_i}N_i}(N_i,N_i,0)[-2],
\end{multline}
where $N_i\times_{K_i}N_i\times_{K_i}N_i$ is embedded in $K\times K\times
K$ in the obvious way via $\iota_i\times\iota_i\times\iota_i$. On restriction
to this, $\pi_{13}$ is a $\PP^1$-bundle over $N_i\times_{K_i}N_i$: the
pullback of $p_i\colon N_i\to K_i$ to the base $N_i\times_{K_i}N_i$. Thus
it has relative canonical bundle the pullback of $\omega_{p_i}\cong\mc{O}(N_i)$,
which is $\mc{O}(0,N_i,0)$. Since $p_{i*}\mc{O}\cong\mc{O}$ and $p_{i*}\omega_{p_i}\cong
\mc{O}[-1]$ by Serre duality, (\ref{UU}) becomes
$$
\mc{O}_{(\iota_i\times\iota_i)(N_i\times_{K_i}N_i)}(N_i,0)[-1]\to
U_i\star U_i \to\mc{O}_{(\iota_i\times\iota_i)(N_i\times_{K_i}N_i)}(N_i,0)[-3],
$$
which by (\ref{UFM}) gives (\ref{UUU}). \bigskip

The morphism $1\star\ev_i'\colon U_i[-2]\to U_i^2$ in $D^b(K\times K)$ is
\begin{multline}
\pi_{13*}\left(\mc{O}_{N_i\times_{K_i}N_i\times K}(N_i,0,0)[-1]\otimes
\mc{O}_{K\times\Delta}[-2]\right. \\ \left. \To \mc{O}_{N_i\times_{K_i}N_i\times
K}(N_i,0,0)[-1]\otimes\mc{O}_{K\times N_i\times_{K_i}N_i}(0,N_i,0)[-1]\right),
\label{b} \end{multline}
where the map is the identity on the left hand sides of the tensor products
and $\ev_i'$ of (\ref{ev'}) on the right hand sides.

We wish to show that the composition of this map with $U_i\star U_i\rt{b}
U_i[-2]$ is the identity. But $b$ was defined by projecting the second derived tensor product in (\ref{b}) to its degree zero (Tor$_0$) part. That is, the
composition is given by
\begin{equation} \label{bb}
\pi_{13*}\left(\mc{O}_{N_i\times_{K_i}\Delta_{N_i}}(N_i,0,0)[-3]\to
\mc{O}_{N_i\times_{K_i}N_i\times_{K_i}N_i}(N_i,N_i,0)[-2]\right),
\end{equation}
where the map is the connecting homomorphism of the natural divisor exact
sequence
$$
0\to\mc{O}(N_i,N_i,0)\to\mc{O}(N_i,N_i,0)(N_i\times_{K_i}\Delta_{N_i})\to
\mc{O}_{N_i\times_{K_i}\Delta_{N_i}}(N_i,0,0)\to0,
$$
(cf. (\ref{ted})) pushed forward from $N_i\times_{K_i}N_i\times_{K_i}N_i$ to $K\times K\times K$ and shifted by $[-3]$.

So the cone on (\ref{bb}) is
$$
\pi_{13*}\left(\mc{O}_{N_i\times_{K_i}N_i\times_{K_i}N_i}(N_i,N_i,0)
(N_i\times_{K_i}\Delta_{N_i})[-2]\right)\ \cong\ 0,
$$
because $\mc{O}(N_i,N_i,0)(N_i\times_{K_i}\Delta_{N_i})$ is a line bundle of degree $-1$ on the $\PP^1$-fibers of $\pi_{13}|_{N_i\times_{K_i}N_i\times_{K_i}N_i}$.

Thus the induced map in (\ref{bb}),
$$
\mc{O}_{N_i\times_{K_i}N_i}(N_i,0)[-3]\to\mc{O}_{N_i\times_{K_i}N_i}(N_i,0)[-3],
$$
must be a nonzero multiple of the identity, and so isomorphic to the identity,
as required.

The proof that $\ev_i'\star1$ also splits (\ref{UUU}) is very similar.
\medskip

To analyse $\ev_i\star1\colon U_i\star U_i\to U_i$ we tensor the following
exact sequence (the product of (\ref{tis}) with $K$) with
$\mc{O}_{K\times N_i\times_{K_i}N_i}(0,N_i,0)$:
$$
0\to\mc{O}_{\Delta\times K}\to\mc{O}_{\big(\Delta\,\cup\,(N_i\times_{K_i}N_i)\big)
\times K}(N_i,0,0)\to\mc{O}_{N_i\times_{K_i}N_i\times K}(N_i,0,0)\to0.
$$
Taking $\pi_{13*}$ of the connecting homomorphism and shifting by $[-2]$
gives $\ev_i\star1$. Now the intersection of $K\times N_i\times_{K_i}N_i$
with both of $\big(\Delta\,\cup\,(N_i\times_{K_i}N_i)\big)\times K$ and 
$N_i\times_{K_i}N_i\times K$ is the same ($N_i\times_{K_i}N_i\times_{K_i}N_i$),
so in the long exact sequence of Tors the final two terms are the same
with the restriction map between them an isomorphism. Thus the next two terms
to the left give the exact sequence
$$
\mathrm{Tor}_1\big(\mc{O}_{N_i\times_{K_i}N_i\times K},
\mc{O}_{K\times N_i\times_{K_i}N_i}\big)(N_i,N_i,0)\to
\mc{O}_{\Delta_{N_i}\times_{K_i}N_i}(N_i,0,0)\to0.
$$
Now recall (\ref{UU}) that the inclusion of this Tor$_1$ term into the full
derived tensor product is precisely what gives the map $a\colon U_i\to U_i\star U_i$; i.e. applying
$\pi_{13*}[-1]$ to the above two terms gives $(\ev_i\star1)\circ a$. We also
computed this Tor$_1$ in (\ref{tor1}); the upshot is that $(\ev_i\star1)\circ a$ is
$$
\pi_{13*}\left[\mc{O}_{N_i\times_{K_i}N_i\times_{K_i}N_i}\to
\mc{O}_{\Delta_{N_i}\times_{K_i}N_i}\right](N_i,0)[-1],
$$
with the map inside the brackets the \emph{surjection of structure sheaves}.
Thus this map is isomorphic to the restriction map; when we apply $\pi_{13*}$
we get
$$
\left[\mc{O}_{N_i\times_{K_i}N_i}\to\mc{O}_{N_i\times_{K_i}N_i}\right](N_i,0)[-1]
$$
on $K\times K$, with the map the identity. So $(\ev_i\star1)\circ a$ is isomorphic
to the identity map $U_i\to U_i$, and $\ev_i\star1$ splits (\ref{UUU}).
The proof for $1\star\ev_i$ is very similar.
\end{proof}

We are now ready to prove invertibility of our functors. This can by now be proved in many ways. It follows as a special
case of the results of Horja \cite{Ho} or similar ideas of Szendr\"oi \cite{Sz1},
for instance, and the quickest proof is probably by adapting Ploog's method \cite{Pl, Hu} for the invertibility of the spherical twists \cite{ST} to this family situation. We give a proof inspired by purely categorical methods
and the earlier invertibility proofs in \cite{KS, ST} and the more general result proved by Rouquier \cite{Ro1}, at the level of Fourier-Mukai transforms.

\begin{prop} \label{invert}
In $D^b(K\times K)$ we have quasi-isomorphisms $T_i\star T_i'\cong
\mc{O}_\Delta\cong T_i'\star T_i$.
\end{prop}

\begin{proof}
$T_i'\star T_i$ is the total complex of
\begin{equation} \label{tt}
\diagram
U_i[-1] \rto^{\ev_i} \dto_{\ev'_i\star\id} & \mc{O}_\Delta[-1] \dto^{\ev'_i}
\\ U_i\star U_i[1] \rto^(.56){\id\star\ev_i} & U_i[1].\!\!
\enddiagram
\end{equation}
We map the trivial cone $U_i[1]\rt{\id}U_i[1]$ into the bottom row of (\ref{tt})
by the map $a[1]$ on the first factor and the identity on the second. By Proposition \ref{spl} this is indeed a morphism in $D^b(K\times K)$ and taking the cone makes $T_i'\star T_i$ isomorphic to the total complex of
$$
\diagram
U_i \rto^{\ev_i} \dto_{b\comp(\ev'_i\star\id)} & \mc{O}_\Delta
\\ U_i,\!\!
\enddiagram
$$
where now the first $U_i$ is in degree 1 and the other two terms are in degree
0. By Proposition \ref{spl} the vertical arrow is isomorphic to the identity, so $T_i'\star T_i\cong\mc{O}_\Delta$. The proof that $T_i\star T_i'\cong\mc{O}_\Delta$ is very similar.
\end{proof}

Defining $Fl_{ij},\,i\ne j$, to be the quotient of $Fl$ given by forgetting both the $i$th and $j$th subspaces, we have natural maps
\begin{equation} \label{FF}
\spreaddiagramrows{-0.5pc}
\spreaddiagramcolumns{-0.5pc}
\diagram
 & Fl \ddto_{\pi_{ij}\!} \dlto_{\pi_i} \drto^{\pi_j} \\
 Fl_i \drto_{\rho_{ij}} && Fl_j\dlto^{\rho_{ji}} \\
 & Fl_{ij}.\!\!
\enddiagram
\end{equation}
$\pi_i$ and $\pi_j$ are $\PP^1$-bundles, and if $|i-j|>1$ then $\rho_{ij}$
and $\rho_{ji}$ are also; in this case the diagram is \emph{Cartesian}, i.e. $\pi_{ij}$ is the fiber product of $\rho_{ij}$ and $\rho_{ji}$. (For $|i-j|=1$,
$\rho_{ij}$ and $\rho_{ji}$ are $\PP^2$-bundles.)

We also define $N_{ij}$ to be the (transverse) intersection $N_i\cap
N_j$ inside $K$. For $|i-j|>1$ this is the pullback $\pi_{ij}^*T^*Fl_{ij}$. Finally we set $K_{ij}\define\rho_{ij}^*T^*Fl_{ij}$, a divisor in $K_i=T^*Fl_i$, the zero set of the canonical section of $T^*\rho_{ij}$. It is the image of $N_{ij}$ under $p_i\colon N_i\to K_i=T^*Fl_i$.

With this notation we can prove the following.

\begin{prop} \label{commute}
For $|i-j|>1$ we have $T_i\star T_j\cong T_j\star T_i\,\in D^b(K\times K)$.
\end{prop}

\begin{proof}
In the usual notation we calculate $T_i\star T_j$ as
\begin{equation} \label{calc}
\pi_{13*}\left(\mc{O}_{K\times\left(\Delta\cup(N_i\times_{K_i}N_i)\right)}
(0,N_i,0)\otimes\mc{O}_{\left(\Delta\cup(N_j\times_{K_j}N_j)\right)\times K}(N_j,0,0)\right).
\end{equation}
Both supports have codimension $\dim K$ in $K\times K\times K$, and their
intersection is
$$
(\id\times\id\times\id)K\ \cup\ (\Delta_{N_i}\times_{K_i}N_i)\ \cup\ 
(N_j\times_{K_j}\Delta_{N_j})\ \cup\ (N_j\times_{K_j}N_{ij}\times_{K_i}N_i).
$$
The last component $N_j\times_{K_j}N_{ij}
\times_{K_i}N_i$ is $N_{ij}\times_{K_{ji}}N_{ij}\times_{K_{ij}}N_{ij}$, which in turn equals $N_{ij}\times_{T^*Fl_{ij}}N_{ij}$ by pulling back $T^*Fl_{ij}$ to both sides of the isomorphism $Fl\times_{Fl_j}Fl\times_{Fl_i}Fl\cong
Fl\times_{Fl_{ij}}Fl$ (which amounts to the statement that (\ref{FF}) is Cartesian).

Therefore all four irreducible components are of codimension $2\dim K$,
the intersection is a local complete intersection and the derived tensor product of the structure sheaves is just the structure sheaf of the intersection.
Moreover, $\pi_{13}$ restricted to each component is an embedding into $K\times K$, yielding
\begin{equation} \label{Tij}
T_i\star T_j\cong\mc{O}_{\Delta\,\cup\,(N_i\times_{K_i}N_i)\,\cup\,
(N_j\times_{K_j}N_j)\,\cup\,(N_{ij}\times_{T^*Fl_{ij}}N_{ij})}(D_{ij}).
\end{equation}
We claim that $\mc{O}(D_{ij})$ is $\mc{O}(N_i+N_j,0)$, i.e. that
$\mc{O}(0,N_i,0)$ in (\ref{calc}) is isomorphic to $\mc{O}(N_i,0,0)$ on
restriction to the intersection. This is obvious on the first two irreducible components of the intersection, and on the
latter two it follows from the fact that $\mc{O}_{N_j}(N_i)$ is pulled back from $K_j$ (for $|i-j|>1$: it is the pullback of $\mc{O}_{K_j}(K_{ji})$).

So $D_{ij}$ and therefore also (\ref{Tij}) are symmetric in $i$ and $j$. \end{proof}

\subsection{The Yang-Baxter equation} \label{YB}

In this section we fix $j=i+1$ and prove the final braid relation $T_i\star T_j\star T_i\cong T_j\star T_i\star T_j$.

The geometry and notation become especially nasty in this section, and
the interested reader is advised to set $n=3$ and $i,j=1,2$, so that
$Fl=Fl(\C^3)$ fibers by $\pi_2$ over $Fl_2=\PP^2$ and by $\pi_1$ over
$Fl_1=\PP^{2*}$. In this case $M_{ij}$ below is just the zero section
$Fl\subset T^*Fl$, while $K_{ij}=\PP^2\subset T^*\PP^2=K_2$ and
$K_{ji}=\PP^{2*}\subset T^*\PP^{2*}=K_1$. Lemma \ref{blowup} is then just the Mukai flop. The geometry
that results in this special case is where all the action takes place in
the general case (which is a bundle over $T^*Fl_{ij}$ with fiber the $n=3$
geometry).

We refer to the diagram (\ref{FF}).
As before we define $N_{ij}$ to be the (transverse) intersection $N_i\cap
N_j$ inside $K$, and we set $M_{ij}\define\pi_{ij}^*T^*Fl_{ij}$. Since $|i-j|=1$
these are not equal, but $M_{ij}\subset N_{ij}$ is a
divisor, in fact an exceptional divisor:

\begin{lemma} \label{blowup}
The map $q_i\define p_i|_{N_{ij}}$ is the blow-up of $K_i=T^*Fl_i$ in
$K_{ij}\define\rho_{ij}^*T^*Fl_{ij}$, and $M_{ij}$ is the exceptional divisor.

Symmetrically, $q_j\define p_j|_{N_{ij}}$ is the blow-up of $K_j=T^*Fl_j$ in $K_{ji}=
\rho_{ji}^*T^*Fl_{ij}$ with the same exceptional divisor $M_{ij}\subset N_{ij}$.
\end{lemma}

\begin{proof}
Thinking of $N_{ij}$ as the zero set of the canonical section of $T^*\pi_j$
restricted to $N_i=\pi_i^*T^*Fl_i$, the functions that are linear on the
fibers of its projection to $Fl$ are
$$
\frac{\pi_i^*TFl_i}{T\pi_j}.
$$
These generate the ideal of the zero section $Fl$; the subset of those that
vanish on the divisor $M_{ij}=\pi_{ij}^*T^*Fl_{ij}\subset N_{ij}$ is
\begin{eqnarray*}
\ker\left(\frac{\pi_i^*TFl_i}{T\pi_j}\To\pi_{ij}^*TFl_{ij}\right) \!\!&=&\!\!
\frac{\pi_i^*T\rho_{ij}}{T\pi_j}\ =\ 
\frac{\hom(V_j/V_{i-1},V_{j+1}/V_j)}{\hom(V_j/V_i,V_{j+1}/V_j)} \\
\!\!&=&\!\! \hom\!\left(\!\frac{V_i}{V_{i-1}},\frac{V_{j+1}}{V_j}\!\right)=
\left(\!\frac{V_i}{V_{i-1}}\!\right)^{\!\!*}\!\otimes\!\left(\!
\frac{V_{j+1}}{V_j}\!\right)\!.
\end{eqnarray*}

Thinking of $\pi_i\colon Fl\to Fl_i$ as the projective bundle
$\PP(V_j/V_{i-1})\to Fl_i$ with tautological line bundle $\mc{O}_\PP(-1)=
V_i/V_{i-1}$, we see that the pushdown of
$(V_i/V_{i-1})^*=\mc{O}_\PP(1)$ is $(V_j/V_{i-1})^*$. Therefore, the pushdown of the above ideal is
$$
\left(\!\frac{V_j}{V_{i-1}}\!\right)^{\!\!*}\!\otimes\!\left(\!
\frac{V_{j+1}}{V_j}\!\right)=\ \hom\!\left(\!\frac{V_j}{V_{i-1}},\frac{V_{j+1}}{V_j}\!\right)
=T\rho_{ij}.
$$
But this is the subset of the linear functions $TFl_i$ on $K_i=T^*Fl_i$ that
vanish on $\rho_{ij}^*T^*Fl_{ij}=K_{ij}$. Thus we have shown that
$$
q_{i*}\big(\mc{O}_{N_{ij}}(-M_{ij})\big)=\I_{K_{ij}\subset K_i}.
$$
Since everything in sight is smooth, this is the defining relation of the blow-up of $K_i$ in $K_{ij}$.
\end{proof}

We will also need the following.

\begin{lemma} \label{blow2}
The natural projection
$M_{ij}\times_{K_{ji}}M_{ij}\to K_{ij}\times_{T^*Fl_{ij}}K_{ij}$
is the blow-up of $K_{ij}\times_{T^*Fl_{ij}}K_{ij}$ in the diagonal $K_{ij}$.
\end{lemma}

\begin{proof}
This is proved by pulling back $T^*Fl_{ij}$ to
\begin{equation} \label{blow3}
Fl\times_{Fl_j}Fl\ \cong\ \Bl_{\Delta_{Fl_i}}(Fl_i\times_{Fl_{ij}}Fl_i).
\end{equation}
To prove (\ref{blow3}) we exhibit a rational map $\xymatrix{
Fl_i\times_{Fl_{ij}}Fl_i \ar@{-->}[r] & Fl_j}$ that blows up the diagonal
$\Delta_{Fl_i}$. This map takes the $V_j$ subspaces of $V$ in the two copies of $Fl_i$ and intersects them to define the $V_i$ subspace of $V$ in $Fl_j$,
and so is ill-defined over the diagonal where this intersection is too big.

We define the map through sections of line bundles, and work relative to
$Fl_{ij}$. Think of $Fl_j\to Fl_{ij}$ as the projective bundle $\PP(V_{j+1}/
V_{i-1})$ with its tautological line bundle $\mc{O}_\PP(-1)=V_i/V_{i-1}$, and work with the line bundle $\Lambda^3(V_{j+1}/V_{i-1})\otimes\mc{O}_\PP(1)$ on $Fl_j$. Its pushdown
to $Fl_{ij}$ is
\begin{equation} \label{lambda2}
\Lambda^3(V_{j+1}/V_{i-1})\otimes(V_{j+1}/V_{i-1})^*\ \cong\ 
\Lambda^2(V_{j+1}/V_{i-1}).
\end{equation}

Similarly we think of $Fl_i\to Fl_{ij}$ as the projective bundle $\PP((V_{j+1}/
V_{i-1})^*)$ with tautological line bundle $\mc{O}_\PP(-1)=(V_{j+1}/V_j)^*$.
The pushdown on $\mc{O}_\PP(1)$ is then $V_{j+1}/V_{i-1}$.

Thus the sections of $\mc{O}(1,1)$ on $Fl_i\times_{Fl_{ij}}Fl_i$ over $Fl_{ij}$
are $(V_{j+1}/V_{i-1})\otimes(V_{j+1}/V_{i-1})$, containing the sections
(\ref{lambda2}) as the antisymmetric tensors. These latter sections vanish precisely (scheme theoretically) on the diagonal $\Delta_{Fl_i}$ and generate
its ideal, so we get the blow-up claimed.

The resulting map $\Bl_{\Delta_{Fl_i}}(Fl_i\times_{Fl_{ij}}Fl_i)\to Fl_j$
is the projection $Fl\times_{Fl_j}Fl\to Fl_j$ since a flag is determined by its projections to $Fl_i$ and $Fl_j$ (i.e. $\pi_i\times\pi_j\colon
Fl\to Fl_i\times Fl_j$ is an embedding).
\end{proof}

We compute $U_i\star U_j$ in the usual way (\ref{conv}):
$$
U_i\star U_j=\pi_{13*}\left(\mc{O}_{N_j\times_{K_j}N_j\times K}(N_j,0,0)[-1]
\otimes\mc{O}_{K\times N_i\times_{K_i}N_i}(0,N_i,0)[-1]\right).
$$
Both supports have codimension in $K\times K\times K$ equal to $\dim K$,
and their intersection
\begin{equation} \label{Z}
N_j\times_{K_j}N_{ij}\times_{K_i}N_i\ \subset\ K\times K\times K
\end{equation}
is smooth (being a $\PP^1\times\PP^1$-bundle over $N_{ij}$) of codimension
$2\dim K$. Thus the intersection is transverse, and
$$
U_i\star U_j=\pi_{13*}\Big(\mc{O}_{N_j\times_{K_j}N_{ij}\times_{K_i}N_i}
(N_j,N_i,0)\Big)[-2].
$$
On restriction to (\ref{Z}), $\pi_{13}$ maps $(a,b,c)$ to $(\iota_i(a),
\iota_j(c))\in K\times K$. A point of the image determines $a$ and $c$ (since
$\iota_i$ and $\iota_j$ are injections) and so $q_j(b)$ and $q_i(b)$ as $p_j(a)$
and $p_i(c)$ respectively. But this determines $b$ completely and so the
map is injective. In fact it is an embedding because $\iota_i,\ \iota_j$ and $q_i\times q_j\colon N_{ij}\into K_i\times K_j$ are (the latter because
$\pi_i\times\pi_j\colon Fl\into Fl_i\times Fl_j$ is). Therefore
\begin{equation} \label{ij}
U_i\star U_j=\mc{O}_{N_j\times_{K_j}N_{ij}\times_{K_i}N_i}(D)[-2]\ \in\
D^b(K\times K),
\end{equation}
where $N_j\times_{K_j}N_{ij}\times_{K_i}N_i$ is embedded into $K\times K$
as described above, and $D$ is the divisor thereon given by restricting
$\mc{O}(N_j,N_i,0)$ from $K\times K\times K$ via the inclusion (\ref{Z}).

Recall the description of the morphism $\ev_i$ as a connecting homomorphism
(\ref{par}). Chasing this through the above derivation of $U_i\star
U_j$ shows that $\id\star\ev_j\colon U_i\star U_j\to U_i$ is given by
(the shift by $[-2]$ of) the connecting homomorphism of the exact sequence
\begin{multline} \label{iorj}
0\to\mc O_{N_i\times_{K_i}N_i}(N_i,0)\to\mc O_{(N_i\times_{K_i}N_i)
\cup\_{N_{ij}\times_{K_i}N_i}(N_j\times_{K_j}N_{ij}\times_{K_i}N_i)}(D) \\
\to\mc O_{N_j\times_{K_j}N_{ij}\times_{K_i}N_i}(D)\to0
\end{multline}
that arises from the fact that $N_i\times_{K_i}N_i$ and
$N_j\times_{K_j}N_{ij}\times_{K_i}N_i$ intersect in $N_{ij}\times_{K_i}N_i$, which is a divisor in both. Here $\mc O(D)$ is the line bundle on the union
which is $\mc O_{N_i\times_{K_i}N_i}(N_i+N_j,0)$ glued to
$\mc O_{N_j\times_{K_j}N_{ij}\times_{K_i}N_i}(D)$ across $N_{ij}\times_{K_i}N_i$ by their common restriction to $\mc O_{N_{ij}\times_{K_i}N_i}
(L_i^{-2}+L_{i-1}+N_j,L_j)$
(by (\ref{omega}) and the fact that $L_j$ is pulled back from a line bundle
$L_j$ on $K_i$ and so can be dragged across the first fiber product).
\medskip

We introduce one more piece of notation. Just as the kernel $U_i$ represents a functor $\iota\_{i*}p_i^*p\_{i*}\iota_i^!$ arising naturally from the diagram (\ref{piota}), the diagram
$$
\begin{array}{l}
\ M_{ij}\stackrel{\iota\_{ij}}{\into}K\\
\ \ \downarrow^{p_{ij}} \\
\!\!T^*Fl_{ij}
\end{array}
$$
induces the functor $\iota\_{ij*}p_{ij}^*p\_{ij*}\iota_{ij}^!$ on $D^b(K)$. By much the same working as in Proposition \ref{UFML} this can be shown to be represented by the kernel
$$
U_{ij}\define\mc{O}_{M_{ij}\times_{T^*Fl_{ij}}M_{ij}}(\omega\_{Fl/Fl_{ij}},0)[-3]
\ \in\,D^b(K\times K).
$$

\begin{prop} \label{H}
There is an isomorphism $U_i\star U_j\star U_i\cong U_i[-2]
\oplus U_{ij}$, such that the composition
of $\,U_i\star U_j\star U_i\Rt{\ev_j}U_i\star U_i$ and the map $U_i\star
U_i\rt{b}U_i[-2]$ of Proposition \ref{spl} acts as the identity on $U_i[-2]$.
\end{prop}

\begin{proof} By (\ref{ij}) and (\ref{conv}), in the notation above $U_i\star U_j\star U_i$ is
$$
\pi_{13*}\left(\mc{O}\_{N_i\times_{K_i}N_i\times K}
(N_i,0,0)[-1]\otimes\mc{O}_{K\times(N_j\times_{K_j}N_{ij}\times_{K_i}N_i)}
(0,D)[-2]\right).
$$
Both supports have codimension $\dim K$ in $K\times K\times K$. Their intersection
\begin{eqnarray} \label{inter}
N_i\times_{K_i}N_{ij}\times_{K_j}N_{ij}\times_{K_i}N_i &\subset& K\times K\times K, \\
(a,b,d,c) \quad &\mapsto& \quad (a,b,c), \nonumber
\end{eqnarray}
has $b=d$ except over the exceptional locus
$M_{ij}\times_{K_j}M_{ij}\subset N_{ij}\times_{K_j}N_{ij}$ where the projections
$q_i\colon N_{ij}\to K_i$ and $q_j\colon N_{ij}\to K_j$ are not isomorphisms.
Therefore (\ref{inter}) is the union $A\cup B$ of two irreducible components; the first where $b=d$,
$$
A=N_i\times_{K_i}N_{ij}\times_{K_i}N_i,
$$
the second the closure of the locus where $b\ne d$,
\begin{equation} \label{B}
B=M_{ij}\times_{K_{ij}}M_{ij}\times_{K_{ji}}M_{ij}\times_{K_{ij}}M_{ij}.
\end{equation}
These both have codimension $2\dim K$, so that the derived tensor product of the structure sheaves is the structure sheaf of $A\cup B$, and
\begin{equation} \label{iji}
U_i\star U_j\star U_i\ \cong\ \pi_{13*}\big(\mc{O}\_{A\cup B}
(N_i,0,0)(0,D)\big)[-3].
\end{equation}
Set $E\define A\cap B=M_{ij}\times_{K_{ij}}M_{ij}\times_{K_{ji}}M_{ij}$, a divisor in both $A$ and $B$.
We will show that the inclusion of the natural subsheaf $\mc{O}_A(-E)
\oplus\mc{O}_B(-E)\into\mc{O}_{A\cup B}$ induces an isomorphism
$$
\pi_{13*}\big(\mc{O}\_A(N_i,0,0)(0,D)(-E)\big)\ \oplus\
\pi_{13*}\big(\mc{O}\_B(N_i,0,0)(0,D)(-E)\big) \vspace{-3mm}
$$
\begin{equation} \label{ideal}
\hspace{5cm} \Rt{\simeq} \pi_{13*}\big(\mc{O}\_{A\cup B}(N_i,0,0)(0,D)\big).
\end{equation}
It is therefore sufficient to show that
\begin{equation} \label{AAA}
\pi_{13*}\big(\mc{O}\_A(N_i,0,0)(0,D)(-E)\big)\to\pi_{13*}\big(\mc{O}\_A(N_i,0,0)(0,D)\big)
\end{equation}
and
\begin{equation} \label{BBB}
\pi_{13*}\big(\mc{O}\_B(N_i,0,0)(0,D)(-E)\big)\to\pi_{13*}\big(\mc{O}\_B(N_i,0,0)(0,D)\big)
\end{equation}
are isomorphisms.
\medskip

On restriction to $A$, $\pi_{13}$ projects out the $b$-variable to
$N_i\times_{K_i}N_i$. As such it is the pullback via $N_i\times_{K_i}N_i\to
K_i$ of the map $q_i\colon N_{ij}\to K_i$ (i.e. the basechange of $q_i$ from
$K_i$ to $N_i\times_{K_i}N_i$).

By Lemma \ref{blowup} $q_i$ is the blow-up of $K_i$ in $K_{ij}$ of codimension two, so $\omega_{q_i}$ is the exceptional divisor $\mc{O}_{N_{ij}}(M_{ij})$. Pulling back to $A$ we get $\mc{O}_A(E)$, where $E$ is the exceptional
divisor of $\pi_{13}|_A$; i.e. $E=M_{ij}\times_{K_{ij}}M_{ij}\times_{K_{ij}}
M_{ij}=A\cap B$.

Since $\omega_{K_i}\cong\mc{O}_{K_i}$, this relative canonical bundle is also $\omega_{N_{ij}}$, which by adjunction is $\mc{O}(N_i+N_j)|_{N_{ij}}$, i.e.
$\mc{O}(0,N_1+N_j,0)$ on pullback to $A$. But since $b=d$ on $A$ we also have
$\mc{O}(0,D)|_A=\mc{O}(0,N_i+N_j,0)$ from its definition (\ref{ij}). The upshot is
that
\begin{equation} \label{nij}
\mc{O}_{N_{ij}}(M_{ij})\cong\mc{O}(N_i+N_j)|_{N_{ij}}
\quad\text{and}\quad \mc{O}_A(0,D)=\mc{O}_A(E).
\end{equation}

It follows that (\ref{AAA}) is $\pi_{13*}\big(\mc{O}_A(N_i,0,0)\big)\to
\pi_{13*}\big(\mc{O}_A(N_i,0,0)(E)\big)$. Since $E$ is exceptional this is an
isomorphism $\mc{O}_{N_i\times_{K_i}N_i}(N_i,0)\Rt{\simeq}\mc{O}_{N_i\times_{K_i}N_i}(N_i,0)$
as required.
\medskip

Something similar happens for $B$. The map $\pi|_B$ projects $B$ to
$M_{ij}\times_{T^*Fl_{ij}}M_{ij}\subset K\times K$, contracting only the
exceptional divisor $E=A\cap B$. It is the basechange via
$M_{ij}\times_{T^*Fl_{ij}}M_{ij}\to K_{ij}\times_{T^*Fl_{ij}}K_{ij}$ of
$$
M_{ij}\times_{K_{ji}}M_{ij}\to K_{ij}\times_{T^*Fl_{ij}}K_{ij},
$$
which by Lemma \ref{blow2} is the blow-up of $K_{ij}\times_{T^*Fl_{ij}}K_{ij}$ in the diagonal $K_{ij}$. Since this diagonal has codimension two, $\mc{O}(E)$ is the relative canonical bundle,
which by computation is $(\omega^{-1}_{Fl_i/Fl_{ij}},N_j,N_i,0)$ (written
in terms of the four factors of $B$ (\ref{B})).

Therefore the line bundle $\mc{O}(N_i,0,0)(0,D)|_B=\mc{O}_B(N_i,N_j,N_i,0)$ is isomorphic to
$(\omega\_{Fl_i/Fl_{ij}}(N_i),0,0,0)(E)=(\omega\_{Fl/Fl_{ij}},0,0,0)(E)$.
Thus (\ref{BBB}) is $\pi_{13*}\big(\mc{O}_B(\omega\_{Fl/Fl_{ij}},0,0,0)\big)\to
\pi_{13*}\big(\mc{O}_B(\omega\_{Fl/Fl_{ij}},0,0,0)(E)\big)$, which is an isomorphism
with both equal to $\mc{O}_{M_{ij}\times_{T^*Fl_{ij}}M_{ij}}(\omega\_{Fl/Fl_{ij}},0)$
since $E\subset B$ is exceptional.
\medskip

So now (\ref{ideal}) follows from (\ref{AAA}) and (\ref{BBB}), and
(\ref{iji}) has become
$$
U_i\star U_j\star U_i\ \cong\ \mc{O}_{N_i\times_{K_i}N_i}(N_i,0)[-3]\ \oplus\
\mc{O}_{M_{ij}\times_{T^*Fl_{ij}}M_{ij}}(\omega\_{Fl/Fl_{ij}},0)[-3].
$$
But this is $U_i[-2]\oplus U_{ij}$, as claimed.
\bigskip

Finally we describe the map $\id\star\ev_j\star\id$. By (\ref{iorj}) we consider
(the shift by $[-2]$ of) the connecting homomorphism of the exact sequence
\begin{multline*}
\!\!\!\!\!0\to\mc O_{K\times N_i\times_{K_i}N_i}(0,N_i,0)\to
\mc O_{(K\times N_i\times_{K_i}N_i)\cup\_{K\times N_{ij}\times_{K_i}N_i}
(K\times N_j\times_{K_j}N_{ij}\times_{K_i}N_i)}(0,D) \\
\to\mc O_{K\times N_j\times_{K_j}N_{ij}\times_{K_i}N_i}(0,D)\to0.\!\!\!\!
\end{multline*}
We must take the derived tensor product with $\mc O_{N_i\times_{K_i}N_i}(N_i,0,0)[-1]$
and then apply $\pi_{13*}$. This derived tensor product has Tor$_0$ and Tor$_1$
terms, as shown above; the Tor$_0$ terms (i.e. the ordinary tensor product) sit inside the exact sequence
\begin{multline} \label{lastun}
\!\!\!\!0\to\mc O_{N_i\times_{K_i}N_i\times_{K_i}N_i}(N_i,N_i,0)\to
\mc O_{(N_i\times_{K_i}N_i\times_{K_i}N_i)\cup\,B}(N_i,0,0)(0,D) \\
\to\mc O_{A\cup B}(N_i,0,0)(0,D)\to0,
\end{multline}
in the above notation. Only the last term arises from a transverse intersection;
i.e. the last term is the full derived tensor product.

Now by mapping the complexes representing the full tensor products to this exact sequence of Tor$_{0\,}$s, a simple argument valid in any derived category shows the following. The connecting
homomorphism of the above exact sequence is the same as
the connecting homomorphism in the exact triangle of full tensor products,
composed with the projection from the full tensor product on the left hand
side to its Tor$_0$ term. But the map $b$ of Proposition \ref{spl} is
precisely ($\pi_{13*}$ of) this projection from $\pi_{12}^*U_i\otimes
\pi_{23}^*U_i$ to Tor$_0(\pi_{12}^*U_i,\pi_{23}^*U_i)$.

So we take $\pi_{13*}$ of the exact sequence (\ref{lastun}). We showed above
that the inclusion of $U_i[-2]$ into $U_i\star U_j\star U_i\cong U_i[-2]\oplus U_{ij}$ is induced by the inclusion of $\mc O_A(N_i,0,0)(0,D)(-E)\cong\mc O_A(N_i,0,0)$ (by (\ref{nij})) into
$\mc O_{A\cup B}(N_i,0,0)(0,D)$, whose fiber product with (\ref{lastun}) gives the exact sequence
\begin{multline*}
0\to\mc O_{N_i\times_{K_i}N_i\times_{K_i}N_i}(N_i,N_i,0)\to
\I_{E\subset N_i\times_{K_i}N_i\times_{K_i}N_i}(N_i,0,0)(0,D) \\
\to\mc O_A(N_i,0,0)\to0.
\end{multline*}
Taking $\pi_{13*}$ of this sequence gives an isomorphism
$R^0\pi_{13*}\mc O_A(N_i,0,0)\cong R^1\pi_{13*}\mc O_{N_i\times_{K_i}
N_i\times_{K_i}N_i}(N_i,N_i,0)$, which by the above working
is the morphism $b\circ(\id\star\ev_j\star\id)|_{U_i[-2]}\colon
U_i[-2]\Rt{\sim}U_i[-2]$ after shifting by $[-2]$.
\end{proof}

\begin{rmk}
One can show that the previous result is equivalent to the following fact
about Mukai flops. (We restrict attention to the $n=3$ case for sake of
exposition; in higher dimensions we just get a family version of
this case.) We let $N_{12}$ be the total space of the $\mc{O}(-1,-1)$-bundle over
$Fl(\C^3)$. This has maps $q_1,\,q_2$ to $T^*(\PP^2)^*,\,T^*\PP^2$ which are both
the blow-up of their image in the zero section. Then on $D^b(T^*\PP^2)$,
$q_{2*}\omega q_1^*q_{1*}q_2^*=\id\,\oplus\,H^*\big(\iota_{\PP^2}^*(\
\ \cdot\ )(-1)\big)\otimes\iota\_{\PP^2*}\mc{O}_{\PP^2}(-2)$,
where $\omega$ is the canonical bundle of $N_{12}$ (and therefore of the
$q_i$). This gives a perhaps more conceptual reason why $q_{1*}q_2^*\colon
D^b(T^*\PP^2)\to D^b(T^*(\PP^2)^*)$ is not an equivalence \cite{Kw, Nm}; if it were its adjoint $q_{2*}\omega q_1^*$ would be its inverse.
\end{rmk}

\begin{prop} \label{hh}
The compositions
$$
\spreaddiagramcolumns{-.4pc}
\spreaddiagramrows{-2pc}
\xymatrix{&& U_i\star U_j\star U_i \ar[drrr]^(.55){\id\star\id\star\ev_i} \\ U_{ij} \urrto\drrto &&&&& U_i\star U_j \\ 
&&  U_j\star U_i\star U_j \ar[urrr]_(.55){\ev_j\star\id\star\id}}
$$
are the same map, which we call $\tau_{ij}$.
\end{prop}

\begin{proof}
The key to this result is that there is a canonical morphism $\tau_{ij}$ in
$$
\Hom^0(U_{ij},U_i\star U_j)
=\Ext^1\!\big(\mc{O}_{M_{ij}\times_{T^*Fl_{ij}}M_{ij}}
(\omega\_{Fl/Fl_{ij}},0),\mc{O}_{N_j\times_{K_j}N_{ij}\times_{K_i}N_i}(D)\big).
$$
Letting $P$ and $Q$ denote $N_j\times_{K_j}N_{ij}\times_{K_i}N_i$ and
$M_{ij}\times_{T^*Fl_{ij}}M_{ij}$ respectively, it is the connecting
homomorphism of the standard exact sequence
\begin{equation} \label{PQ}
0\to\mc{O}_P(D)\to\mc{O}_{P\cup Q}(2N_i+2N_j,0)\to\mc{O}_Q(\omega\_{Fl/Fl_{ij}},0)\to0.
\end{equation}
Here we are using a number of easily checked facts. Firstly, on $P$, $\mc{O}(D)$ (\ref{ij}) is isomorphic to
$\mc{O}(\omega\_{Fl/Fl_{ij}},-N_i-N_j,0)$ (written with respect to the three
factors of $N_j\times_{K_j}N_{ij}\times_{K_i}N_i$). This can be proved by
expanding everything in terms of the $L_i$ using (\ref{omega}) and
$\omega\_{Fl/Fl_{ij}}\cong L_i^{-2}L_j^{-2}L_{i-1}L_{j+1}=\mc{O}(2N_i+2N_j)$,
and then moving $L_{k\ne j}$ across the first fiber product using the fact
that such line bundles are pulled back from $K_j$.

Secondly, $P\cap Q=M_{ij}\times_{K_{ij}}M_{ij}\times_{K_{ji}}M_{ij}$ is 
a divisor in both $P$ and $Q$, and as a divisor in $Q$ it defines the line
bundle $\mc{O}_Q(0,N_i+N_j,0)$. This is because $P\cap Q\subset Q$ is the pullback
of $M_{ij}\subset N_{ij}$ under $N_j\times_{K_j}N_{ij}\times_{K_i}N_i\to
N_{ij}$, and $\mc{O}_{N_{ij}}(M_{ij})$ is $\mc{O}(N_i+N_j)|_{N_{ij}}$ (\ref{nij}).
\medskip

By symmetry both maps of Proposition \ref{hh} must be the same multiple of
this canonical morphism $\tau_{ij}$, but for completeness we check that this
multiple is $1$.

Using the description (\ref{par}) of $U_i\to\mc{O}_\Delta$ as the connecting
homomorphism of the exact sequence
$$
0\to\mc{O}_\Delta\to\mc{O}_{\Delta\,\cup\,(N_i\times_{K_i}N_i)}
(N_i,0)\to\mc{O}_{N_i\times_{K_i}N_i}(N_i,0)\to0,
$$
we pull back by $\pi_{12}^*$ to $K\times K\times K$, tensor with $\pi_{23}^*
(U_i\star U_j)$ and pushdown by $\pi_{13*}$ to give $\id\star\id\star\ev_i
\colon U_i\star U_j\star U_i\to U_i\star U_j$. The proof of Proposition \ref{H}
computes this tensor product and pushdown in detail so we follow our morphism
(\ref{par}) through that proof. We see that $\id\star\id\star\ev_i$
is represented by $\pi_{13*}$ of the connecting homomorphism (shifted by $[-3]$) of the exact sequence
\begin{multline*}
\!\!\!0\to\mc{O}_{\Delta_{N_j}\times_{K_j}N_{ij}\times_{K_i}N_i}(0,N_j,N_i,0)\to
\\ \mc{O}_{(\Delta_{N_j}\times_{K_j}N_{ij}\times_{K_i}N_i)\,\cup\,(A\cup B)}
(N_i,N_j,N_i,0)\to\mc{O}_{A\cup B}(N_i,N_j,N_i,0)\to0,\!\!\!
\end{multline*}
in the notation of the proof of Proposition \ref{H}.

%Denote $\Delta_{N_j}\times_{K_j}N_{ij}\times_{K_i}N_i$ by $\delta$.
$U_{ij}\to U_i\star U_j\star U_i$ is induced by the inclusion of $\mc{O}_B
(N_i,N_j,N_i,0)(-E)$ into the last sheaf in the above exact sequence. Therefore
we pull back the exact sequence by this map to give
\begin{multline*}
\hspace{-3mm} 0\to\mc{O}_{\Delta_{N_j}\times_{K_j}N_{ij}\times_{K_i}N_i}
(0,N_j,N_i,0)\to \\ \hspace{1cm}
\mc{O}_{(\Delta_{N_j}\times_{K_j}N_{ij}\times_{K_i}N_i)\,\cup\,B}
(N_i,N_j,N_i,0)(L)\to\mc{O}_B(N_i,N_j,N_i,0)(L)\to0,
\end{multline*}
whose connecting homomorphism gives (after applying $\pi_{13*}[-3]$) the composition $U_{ij}\to U_i\star U_j\star U_i\to U_i\star U_j$ we require.
Here
$$
L=\mc{O}_{(\Delta_{N_j}\times_{K_j}N_{ij}\times_{K_i}N_i)\,\cup\,B}
(0,L_j^{-1}L_{j+1},L_j^{-1}L_i,0)
$$
is the line bundle defined by minus the divisor given by the intersection
of $A$ with $\Delta_{N_j}\times_{K_j}N_{ij}\times_{K_i}N_i$.

Pushing down as in the proof of Proposition \ref{H} gives the exact sequence
(\ref{PQ}), so the composition we require is its connecting homomorphism
-- our canonical morphism $\tau_{ij}$.

By symmetry the same is true of the other morphism.
\end{proof}

\begin{theorem} \label{YangB}
There is an isomorphism $U_i\star T_j\star
T_i\cong T_j\star T_i\star U_j$ which intertwines the evaluation maps
$\ev_i\star\id\star\id$ and $\id\star\id\star\ev_j$ from the two kernels
to $T_j\star T_i$.

Taking cones on these two maps, $T_i\star T_j\star T_i\cong T_i\star T_j\star T_i$ for $|i-j|=1$.
\end{theorem}

\begin{proof}
Expanding out $U_i\star T_j\star T_i$ gives the cone on
\begin{equation} \label{Dg}
\spreaddiagramcolumns{2.5pc}
\spreaddiagramrows{.5pc}
\diagram
U_i\star U_j\star U_i \rto^(.55){\id\star\id\star\ev_i}
\dto_{\id\star\ev_j\star\id}
& U_i\star U_j \dto^{\id\star\ev_j} \\
U_i\star U_i \rto_{\id\star\ev_i} & U_i.\!
\enddiagram
\end{equation}
If we use $\id\star\ev_i$ to split $U_i\star U_i$ as in Proposition \ref{spl},
then by Propositions \ref{H}, \ref{hh} and \ref{spl} we get
\begin{equation} \label{dg}
\spreaddiagramcolumns{1.5pc}
\spreaddiagramrows{.5pc}
\diagram
U_{ij}\oplus U_i[-2] \rto^(.55){\tau_{ij}\,\oplus\,u_1} \dto_{\id\star\ev_j\star\id}
& U_i\star U_j \dto^{\id\star\ev_j} \\
U_i\oplus U_i[-2] \rto_(.6){\id\oplus\,0} & U_i.\!
\enddiagram
\end{equation}
where the exact form of $u_1$ will not concern us.

We map
\begin{equation} \label{coneq}
\mathrm{Cone}\big(U_{ij}\Rt{\tau_{ij}}U_i\star U_j\big)
\end{equation}
into the cone on (\ref{dg}) as follows. $U_{ij}$
maps to the first term of (\ref{dg}) as the first summand. $U_i\star U_j$ maps to the top right term via the identity and to the bottom left term by
$\id\star\ev_j\,\oplus\,0$. To check this is really a morphism in $D^b(K\times
K)$ we need only check that the first vertical arrow of (\ref{dg}) acts on $U_{ij}$ as $(\id\star\ev_j)\circ\tau_{ij}\oplus\,0$. But the first factor of
this map is determined by the commutativity of (\ref{Dg}), while for the
second we note that
$$
\Hom^0(U_{ij},U_i[-2])\cong\Hom^0\Big(\mc{O}_{M_{ij}\times_{T^*Fl_{ij}}M_{ij}}
(\omega\_{Fl/Fl_{ij}},0),\mc O_{N_i\times_{K_i}N_i}(N_i,0)\Big)
$$
vanishes, since there are no morphisms $\mc O_A\to\mc O_B$ unless there is
an irreducible component of $B$ contained in $A$.

This morphism is a quasi-isomorphism, since by the second part of Proposition \ref{H} what remains is
$$ \spreaddiagramcolumns{1.5pc}
\spreaddiagramrows{.5pc}
\diagram
U_i[-2] \dto_{u_2\oplus\id} \\
U_i\oplus U_i[-2] \rto_(.6){\id\oplus\,0} & U_i,
\enddiagram
$$
for some $u_2$, by Proposition \ref{H}. But this is quasi-isomorphic to zero.

So by symmetry there is a similar quasi-isomorphism from Cone$\big(U_{ij}
\Rt{\tau_{ij}}U_i\star U_j\big)$ to $T_j\star T_i\star U_j$, making the latter isomorphic to $U_i\star T_j\star T_i$. \bigskip

Composing $\ev_i\star\id\star\id\colon U_i\star T_j\star T_i\to T_j\star T_i$ with the quasi-isomorphism from (\ref{coneq}) gives the morphism of
cones
$$
\spreaddiagramcolumns{1.5pc}
\spreaddiagramrows{.5pc}
\diagram
U_{ij} \rto^(.45){\tau_{ij}} \ar@{-->}@/_3pc/[rrr]^{\tau_{ji}} &
U_i\star U_j \ar@{-->}^(.35){\id\star\ev_j}[rrd]
\ar@{-->}@/^2.5pc/[rrr]^(.3){\ev_i\star\id} &&
U_j\star U_i \rto^(.55){\id\star\ev_i} \dto^{\ev_j\star\id} & U_j \dto^{\ev_j} \\ &&& U_i \rto_{\ev_i} & \id
\enddiagram
$$
But by similar working this is exactly the morphism of cones that comes from
composing $\id\star\id\star\ev_j\colon T_j\star T_i\star U_j\to T_j\star
T_i$ with the quasi-isomorphism from $\mathrm{Cone}\big(U_{ij}\Rt{\tau_{ij}}U_i\star U_j\big)$ to $T_j\star T_i\star U_j$.
\end{proof}

\begin{rmk} We note in passing that the faithfulness results
of \cite{KS, ST} can be used to show that our braid group action on $D^b(T^*Fl)$
is faithful. The transverse slice to a subregular nilpotent matrix gives a surface $X$ in $K=T^*Fl$ with an $A_n$ chain of $-2$-spheres $C_i$. The
$C_i$ are the transverse intersection of the subvarieties $N_i$ with the surface. It follows that the Fourier-Mukai kernels (\ref{UFM}, \ref{TFM}) for $U_i$
and $T_i$ restrict, on $X\times X\subset K\times K$, to the standard kernels
for the Dehn twists \cite{ST} about the structure sheaves $\mc{O}_{C_i}$ of the $C_i$ -- see (\ref{tri1}, \ref{tri2}) where these are reviewed.

Therefore the two types of derived category twist, and the briad group actions
they generate, are intertwined by restriction from $K$ to $X$; i.e. the following
diagram commutes.
$$
\xymatrix{
D^b(K) \rto^{T_i}\dto & D^b(K) \dto \\
D^b(X) \rto^{T_i} & D^b(X).
}$$
Thus the faithfulness of the lower action implies the same of the upper action.
\end{rmk}

\subsection{Extension to the affine braid group} \label{affine}

As noted in \cite[Example 3.9]{ST}, the braid group action of that paper
on the derived category of coherent sheaves on a surface with an $A_{n-1}$-chain
of $-2$-curves $C_i$ can be extended to an action of the \emph{affine} braid group. If we take generators of the braid group to be
$T_i\define T_{\mc{O}_{C_i}}$, the Dehn twists about the structure sheaves
$\mc{O}_{C_i}$ described in Section \ref{surfaces}, then we define the extra generator
\begin{equation} \label{T0n}
T=T_0=T_n\define T_{\mc{O}_E(1,\ldots,1)}.
\end{equation}
This is the Dehn twist about $\mc{O}_E(1,\ldots,1)$, the structure sheaf of the entire chain $E=\cup_iC_i$ twisted by the line bundle on $E$ which has degree 1 on each $C_i$. This commutes with $T_i$ for $2\le i\le n-2$ and braids with $T_1$ and $T_{n-1}$, so defines an affine braid group action.

Now it is an easy calculation that $T_{\mc{O}_{C_i}}\big(\mc{O}_{C_{i+1}\cup
\ldots\cup C_{n-1}}(0,1,\ldots,1)\big)$ is isomorphic to $\mc{O}_{C_i\cup\ldots
\cup C_{n-1}}(0,1,\ldots,1)$. Thus, by induction,
$$
\L_1T_1T_2\ldots T_{n-2}(\mc{O}_{C_{n-1}})\cong\mc{O}_E(1,1\ldots,1),
$$
where $\L_1$ is the functor of tensoring with a line bundle $L_1$ whose restriction
to $E$ is $\mc{O}_E(1,0,\ldots,0)$ (for the purposes of this motivational section we may assume that $\mc{O}_E(1,0,\ldots,0)$ extends to a line bundle $L_1$ on the surface).

Thus, by \cite[Lemma 2.11]{ST}, the twist $T=T_{\mc{O}_E(1,\ldots,1)}$
can be written as
\begin{equation} \label{TT0n}
T=\mathrm{ad}(\L_1T_1\ldots T_{n-2})T_{\mc{O}_{C_{n-1}}}=(\L_1T_1\ldots T_{n-2})
T_{n-1}(T_{n-2}'\ldots T_1'\L_1^{-1}).
\end{equation}
\smallskip

Returning to the cotangent bundle of the flag variety, the definition (\ref{T0n})
does not easily generalise to our situation, but (\ref{TT0n}) certainly does.

Define the kernel $T=T_0=T_n\in D^b(K\times K)$ to be
\begin{equation} \label{T}
\big((T_1\star T_2\star\ldots\star T_{n-2})\star T_{n-1}\star (T_{n-2}'\star \ldots\star T_2'\star T_1')\big)(L_1^{-1},L_1).
\end{equation}

From now on we will often suppress the $\star$s in similar long formulae.
We also introduce the notation that if $L$ is a line bundle on $K$ then
$\L\define\mc{O}_\Delta\otimes\pi_1^*L=\mc{O}_\Delta\otimes\pi_2^*L$ is the kernel for $L\otimes(\ \cdot\ )$. That is, $\L$ is the line bundle $L$ supported on the diagonal of $K\times K$, and we have $S\star\L=S(L,0),\ \L\star S=S(0,L)$
for any kernel $S$. Thus (\ref{T}) can be rewritten
\begin{equation} \label{TT}
T=(\L_1T_1T_2\ldots T_{n-2})T_{n-1}(T_{n-2}'\ldots T_2'T_1'\L_1^{-1}).
\end{equation}
For $i\ne j$, $L_i$ is
pulled back from $K_j$, so the line bundles $(L_i,0)$ and $(0,L_i)$
are isomorphic on restriction to $\Delta\cup(N_j\times_{K_j}N_j)$. Thus 
$T_j(0,L_i)=T_j(L_i,0)$, i.e. $\L_i\star T_j=T_j\star\L_i$ for $i\ne j$.
Using this and the braiding we have already proved, we automatically get the following.

\begin{prop}
For $2\le i\le n-2$, we have $T\star T_i\cong T_i\star T$.
\end{prop}

\begin{proof} As explained above we can commute $T_i$ past $\L_1^{-1}$ in (\ref{TT}), and then past $T_j',\,j<i-1$ using Propositions \ref{commute} and \ref{invert}. This gives
$$
T\star T_i\cong
(\L_1T_1T_2\ldots T_{n-2})T_{n-1}(T_{n-2}'\ldots T'_iT'_{i-1})T_i(T'_{i-2}\ldots
T_2'T_1'\L_1^{-1}).
$$
But $T_i$ and $T_{i-1}$ braid, so that $T'_iT'_{i-1}T_i\cong T_{i-1}T'_iT_{i-1}'$,
and then this new $T_{i-1}$ commutes past $T'_j$ and $T_j$ for $j>i$,
yielding
$$
(\L_1T_1\ldots T_{i-1}T_i)T_{i-1}(T_{i+1}\ldots T_{n-2})T_{n-1}(T_{n-2}'\ldots T'_iT_{i-1}'T'_{i-2}\ldots
T_1'\L_1^{-1}).
$$
So using $T_{i-1}T_iT_{i-1}\cong T_iT_{i-1}T_i$ we get
$$
(\L_1T_1\ldots T_{i-2}T_iT_{i-1}T_iT_{i+1}\ldots T_{n-2})T_{n-1}(T_{n-2}'\ldots T_1'\L_1^{-1}).
$$
Commuting $T_i$ past $\L_1T_1\ldots T_{i-2}$ gives $T_i\star T$, as required.
\end{proof}

Similarly the next braid relation is also a pure algebraic formality now
that we know that $T_{n-1}$ commutes with $\L_1$:

\begin{prop} \label{braid2}
$T\star T_{n-1}\star T\cong T_{n-1}\star T\star T_{n-1}$.
\end{prop}

\begin{proof}
$T\star T_{n-1}\star T$ is
$$
(\L_1T_1\ldots T_{n-2})T_{n-1}(T_{n-2}'\ldots T_1'\L_1^{-1})T_{n-1}
(\L_1T_1\ldots T_{n-2})T_{n-1}(T_{n-2}'\ldots T_1'\L_1^{-1}).
$$
We commute the middle $T_{n-1}$ past $T_{n-3}'\ldots T_1'\L_1^{-1}$ to give
$$
(\L_1T_1\ldots T_{n-2})T_{n-1}T_{n-2}'T_{n-1} T_{n-2}T_{n-1}(T_{n-2}'\ldots T_1'\L_1^{-1}).
$$
Because $T_{n-2}$ and $T_{n-1}$ braid, $T_{n-2}'T_{n-1} T_{n-2}T_{n-1}T_{n-2}'
\cong T_{n-1}$, leaving
\begin{equation} \label{TTT}
(\L_1T_1\ldots T_{n-2})T_{n-1}T_{n-1}(T_{n-3}'\ldots T_1'\L_1^{-1}).
\end{equation}
Similarly, $T_{n-1}\star T\star T_{n-1}$ is
$$
(\L_1T_1\ldots T_{n-3})T_{n-1}T_{n-2}T_{n-1}T_{n-2}'T_{n-1}(T_{n-3}'\ldots T_1'\L_1^{-1}),
$$
which by the braiding of $T_{n-2}$ and $T_{n-1}$ is
$$
(\L_1T_1\ldots T_{n-3})T_{n-2}T_{n-1}T_{n-1}(T_{n-3}'\ldots T_1'\L_1^{-1}).
$$
But this is (\ref{TTT}).
\end{proof}

The final braid relation is not quite such a formality, but requires also
the observation from (\ref{TFM}, \ref{T'FM}) that $T_i'\cong
T_i(L_iL_{i-1}^{-1},L_iL_{i+1}^{-1})$, i.e.
\begin{equation} \label{TT'}
T_i'\cong\L_i\star\L_{i+1}^{-1}\star T_i\star\L_i\star\L_{i-1}^{-1}.
\end{equation}

\begin{prop}
$T\star T_1\star T\cong T_1\star T\star T_1$.
\end{prop}

\begin{proof}
By the same reasoning as in Proposition \ref{braid2}, $\,T_1$
braids with
$$
S\define(\L_{n-1}T_{n-1}\ldots T_2) T_1(T_2'\ldots
T_{n-1}'\L_{n-1}^{-1}),
$$
i.e. $T_1\star S\star T_1\cong S\star T_1\star S$. But we know that
$T_1$ commutes with $\L_2$, so the above also equals
\begin{equation} \label{mess}
S=(\L_{n-1}T_{n-1}\ldots T_2\L_2)T_1
(\L_2^{-1}T_2'\ldots T_{n-1}'\L_{n-1}^{-1}).
\end{equation}

We can rewrite $T$ (\ref{T}) as
$$
T=(\L_1T_{n-1}'\ldots T_2')T_1(T_2\ldots T_{n-1}\L_1^{-1}),
$$
by using the braid relations first for $T_{n-1}$ and $T_{n-2}$, then $T_{n-2}$
and $T_{n-3}$ on the result, all the way down to $T_2$ and $T_1$. Commuting
$\L_1$ with $T_i$ for all $i\ge2$, this can be rewritten
\begin{equation} \label{ttt}
T=(T_{n-1}'\ldots T_2'\L_1)T_1(\L_1^{-1}T_2\ldots T_{n-1}).
\end{equation}

But by (\ref{TT'}), $T_{n-1}'\ldots T_2'\L_1=\L_{n-1}\L_n^{-1}T_{n-1}
T_{n-2}\ldots T_2\L_2$ and $L_n$ is trivial, so we can rewrite (\ref{ttt}) as (\ref{mess}). Thus $T=S$, which braids with $T_1$.
\end{proof}

\section{Actions of the braid cobordism category on derived categories
of coherent sheaves}

We now extend the braid group actions of \cite{ST} (on surfaces $X$ with $A_n$-chains of $-2$-spheres $C_i$) and Section
\ref{flagsection} (on $T^*Fl$) to actions of the braid cobordism category.
The latter case proceeds along similar lines to the former, which we tackle
now.

We are guided by the action of the braid cobordism category
on the homotopy category $\mc C(A_n)$ of Section
\ref{sec-survey}. Using the formality result of \cite{KS} it is shown in \cite{ST} that the differential graded algebra
$\RHom(\oplus_i\mc O_{C_i},\oplus_i\mc O_{C_i})$ is quasi-isomorphic to
$A_n\otimes\C$, so that its derived category of differential graded modules is equivalent is that of $A_n\otimes\C$. Under this equivalence, $\mc O_{C_i}$ is mapped to $P_i$, the twists $T_i,\,T_i'$ reviewed below (\ref{tri1}, \ref{tri2}) map to $R_i,\,R_i'$, and $\ev_i,\,\ev_i'$
become $\beta_i,\,\gamma_i$.

So the work in the next section is to lift the positive (\ref{positivemove1},
\ref{positivemove2}) and negative (\ref{negativemove1}, \ref{negativemove2})
braid moves on $\mc C(A_n)$ to $D^b(X)$. This is done in (\ref{moves})
and (\ref{-move}) respectively. The positive moves are fairly easy, the negative
moves require some technical work, in particular Proposition \ref{evev}.

\subsection{Surfaces with chains of rational curves} \label{surfaces}

Let $X$ be a smooth quasi-projective surface with an $A_n$-chain of $-2$-spheres
$C_i$: the exceptional set of the minimal resolution of an $A_n$ surface
singularity. Then $\omega_X$ is trivial in a neighbourhood of $E=\cup_iC_i$; for simplicity we will assume that $\omega_X$ is globally trivial and 
\emph{fix} a trivialisation $\theta\in H^0(\omega_X)$ (with a little more work what follows can be pushed through using a trivialisation defined only on some Zariski open subset $U\supset E$ of $X$).

A (weak) braid group action is constructed on the bounded derived category of coherent sheaves $D^b(X)$ on $X$ in \cite{ST}, with generators $T_i$ fitting
into distinguished triangles
\begin{equation} \label{tri1}
U_i(\F)\define\mc{O}_{C_i}\otimes\RHom(\mc{O}_{C_i},\F)\Rt{\ev_i}\F\to T_i(\F),
\end{equation}
and inverses $T_i'$
\begin{equation} \label{tri2}
T_i'(\F)\to\F\Rt{\ev_i'}\mc{O}_{C_i}\otimes\RHom(\mc{O}_{C_i},\F)[2]=U_i(\F)[2].
\end{equation}
Using very similar (but easier) methods to those in Section \ref{FMTs} one shows that these are represented by the Fourier-Mukai kernels $U_i\define
\mc O_{C_i}^\vee\boxtimes\mc O_{C_i}\cong\mc{O}_{C_i\times C_i}(C_i,0)[-1]$,
$$
T_i\define\mc{O}_{\Delta\cup(C_i\times C_i)}(C_i,0) \quad\mathrm{and}\quad
T_i'\define\mc{O}_{\Delta\cup(C_i\times C_i)}(L_i)
$$
on $X\times X$ (cf. (\ref{TFM}) and (\ref{T'FM})). Here $L_i$ is
$\mc{O}_{\Delta}$ glued across $\Delta_{C_i}$ to
$\mc{O}_{C_i\times C_i}(C_i,0)(-\Delta_{C_i})$, using the fact that the latter
restricts on $\Delta_{C_i}$ to $\Omega^1_{C_i}(C_i)$ which is isomorphic to
$\mc{O}_{\Delta_{C_i}}$ using adjunction and the trivialisation $\theta$ of $\omega_X|_{C_i}$.

The exact triangles (\ref{tri1}) and (\ref{tri2}) then arise from the exact
sequences of Fourier-Mukai kernels
\begin{equation} \label{ex1}
0\to\mc{O}_{\Delta}\to\mc{O}_{\Delta\cup(C_i\times C_i)}(C_i,0)
\rt{r}\mc{O}_{C_i\times C_i}(C_i,0)\to0
\end{equation}
and
\begin{equation} \label{ex2}
0\to\mc{O}_{C_i\times C_i}(C_i,0)\rt{i}\mc{O}_{\Delta\cup(C_i\times C_i)}
(L_i)\to\mc{O}_{\Delta}\to0.
\end{equation}
Taking the Yoneda product of these two extensions gives the connecting morphism
$\ev_i\circ\ev_i'\colon\mc{O}_{\Delta}\to\mc{O}_{\Delta}[2]$ in the
exact triangle
\begin{equation} \label{ii}
\mc{O}_{\Delta}[1]\to\mathrm{Cone}\Big(\mc{O}_{\Delta\cup(C_i\times C_i)}
(C_i,0)\Rt{i\circ r}\mc{O}_{\Delta\cup(C_i\times C_i)}(L_i)\Big)
\to\mc{O}_{\Delta}
\end{equation}
which arises from the exact sequence given by splicing (\ref{ex1})
and (\ref{ex2}) together:
$$
0\to\mc{O}_{\Delta}\to\mc{O}_{\Delta\cup(C_i\times C_i)}(C_i,0)\Rt{i\circ r}\mc{O}_{\Delta\cup(C_i\times C_i)}(L_i)\to\mc{O}_{\Delta}\to0.
$$
The sheaf $\mc{O}_{C_i\times C_i}(C_i,0)(-\Delta_{C_i})$ naturally injects
into both terms of the Cone in (\ref{ii}); take cokernels to give the
following.

\begin{lemma} \label{OO}
The natural transformation $\ev_i\circ\ev_i'\colon\id\to\id[2]$ is represented
by the morphism $\ev_i\circ\ev_i'\colon\mc{O}_{\Delta}\to\mc{O}_{\Delta}[2]$
in the exact triangle
$$
\mc{O}_{\Delta}[1]\to\mathrm{Cone}\Big(\mc{O}_{\Delta}(\Delta_{C_i})
\Rt{i\circ r}\mc{O}_{\Delta\cup\,2\Delta_{C_i}}\Big)
\to\mc{O}_{\Delta}\to\mc{O}_{\Delta}[2],
$$
where $2\Delta_{C_i}$ denotes the thickening of $\Delta_{C_i}$ in $C_i\times
C_i$ only; i.e. the pushforward to $X\times X$ of the subscheme of $C_i\times C_i$ defined by the ideal sheaf $\I_{\Delta_{C_i}}^{\,2}$. \hfill $\square$
\end{lemma}

It is a result of Gerstenhaber and Schack \cite{K} that there is a canonical
isomorphism
$$
\Ext^{i+j}(\mc{O}_{\Delta},\mc{O}_{\Delta})\ \cong\ \bigoplus_{i+j=k}H^i(\mc Ext^j(\mc{O}_{\Delta},\mc{O}_{\Delta}))\ \cong\ \bigoplus_{i+j=k}H^i(\Lambda^j
T_{\Delta}).
$$
In our situation, fixing $i+j=2$ and contracting with the trivialisation $\theta$ of $\omega_X$ gives
\begin{equation} \label{GS}
\Ext^2(\mc{O}_{\Delta},\mc{O}_{\Delta})\ \cong\ H^0(\Omega^2_X)\oplus
H^1(\Omega^1_X)\oplus H^2(\mc{O}_X).
\end{equation}
Later it will be important that here we have fixed an identification of $T_X$
with the normal bundle
to $\Delta$; our convention is that in the splitting $T_{X\times
X}|\_{\Delta}\cong T_X\oplus T_X$ induced by the two factors, $T_X$ maps
to the normal bundle as the image of
\begin{equation} \label{v-v}
TX\ni v\mapsto(v,-v)\in T_X\oplus T_X.
\end{equation}

As a divisor, $C_i\subset X$ defines a de Rham class $[C_i]\in H^1(\Omega^1_X)$
in the following way (which gives the usual cohomology class
$[C_i]\in H^{1,1}(X)$ if $X$ is compact). Think of $H^1(\Omega^1_X)$ as
$\Ext^1(\mc{O}_X,\Omega^1_X)$, which has a map from $\Ext^1(\mc{O}_{C_i},
\Omega^1_X)=H^0(\mc Ext^1(\mc{O}_{C_i},\Omega^1_X))=H^0(\Omega^1_X(C_i))=
\Hom(T_X,\nu\_{C_i})$. The canonical element of the latter (projection of a tangent vector to the normal bundle) defines the extension
$$
0\to\Omega^1_X\to\Omega^1_X(\log C_i)\to\mc{O}_{C_i}\to0,
$$
where the last map takes the residue of logarithmic forms on $C_i$. The map
$\mc O_X\to\mc O_{C_i}$ pulls this back to the extension
\begin{equation} \label{dR}
0\to\Omega^1_X\to\Omega^1_X(\log C_i)\times_{\mc O_{C_i}}\mc O_X\to\mc{O}_X\to0,
\end{equation}
whose extension class is $[C_i]\in H^1(\Omega^1_X)$.

\begin{prop} \label{121}
The class of $\ev_i\circ\ev_i'\in\Ext^2(\mc{O}_{\Delta},\mc{O}_{\Delta})$ of Lemma \ref{OO} can be described, via the isomorphism (\ref{GS}), as the de Rham class $[C_i]\in H^1(\Omega^1_X)\subset\Ext^2(\mc{O}_{\Delta},
\mc{O}_{\Delta})$.
\end{prop}

\begin{proof}
Dualising (\ref{dR}) shows that $[C_i]\in H^1(\Omega^1_X)=\Ext^1(T_X,\mc
O_X)$ defines the extension
\begin{equation} \label{F1}
0\to\mc O_X\to\F\to T_X\to0,
\end{equation}
where $\F$ is the dual of $\big\{(\sigma,f)\in\Omega^1_X(\log C_i)\oplus\mc O_X\colon Res_{C_i}\sigma=f|_{C_i}\in\mc O_{C_i}\big\}=\Omega^1_X(\log C_i)
\times_{\mc O_{C_i}}\mc O_X$. Thus it is easy to see that $\F$ is
\begin{equation} \label{ft}
\big\{(v,s)\in T_X\oplus\mc O(C_i)\colon \pi(v)=s|_{C_i}\big\}=
T_X\times_{\nu\_{C_i}}\mc O(C_i),
\end{equation}
where $\pi\colon T_X\to\nu\_{C_i}$ is the canonical projection to $\mc
O_{C_i}(C_i)$. It therefore also sits inside the exact sequence
\begin{equation} \label{F2}
0\to T_X(-\log C_i)\to\F\to\mc O(C_i)\to0,
\end{equation}
where $T_X(-\log C_i)\define\ker(T_X\rt{\pi}\nu\_{C_i})$ is the dual of $\Omega^1_X
(\log C_i)$.

Meanwhile $\mc Ext^1_{X\times X}(\mc O_{\Delta},\mc O_{\Delta})
\cong\Omega^1_{\Delta}$ is represented by the universal extension
\begin{equation} \label{doublediag}
0\to\Omega^1_{\Delta}\to\mc O_{2\Delta}\to\mc O_{\Delta}\to0,
\end{equation}
which contraction with $\theta\in H^0(\omega_X)$ makes into
\begin{equation} \label{2delta}
0\to T_{\Delta}\to\mc O_{2\Delta}\to\mc O_{\Delta}\to0.
\end{equation}
Therefore $[C_i]\in H^1(\mc Ext^1(\mc O_{\Delta},\mc O_{\Delta}))\subset
\Ext^2(\mc O_{\Delta},\mc O_{\Delta})$ is
represented by the splicing together of the exact sequences (\ref{2delta})
and (the pushforward to $\Delta$ of) (\ref{F1}):
\begin{equation} \label{d}
0\to\mc O_{\Delta}\to\F\to\mc O_{2\Delta}\to\mc O_{\Delta}\to0.
\end{equation}
Now $T_X(-\log C_i)$ injects into $\F$ (\ref{F2}) and under $\theta$ it injects
into $\Omega^1_X$ as $\Omega^1_X(\log C_i)(-C_i)$, with quotient
$\Omega^1_{C_i}$. Therefore it injects into $\mc O_{2\Delta}$ with quotient $\mc O_{\Delta\cup\,2\Delta_{C_i}}$, with $C_i$ thickened only inside $C_i\times C_i$ as described in Lemma \ref{OO}. So taking the quotient of the central two terms of (\ref{d}) by $T_X(-\log C_i)$ gives, by (\ref{F2}), the exact sequence
$$
0\to\mc O_{\Delta}\to\mc O_{\Delta}(\Delta_{C_i})\to\mc O_{\Delta\cup\,
2\Delta_{C_i}}\to\mc O_{\Delta}\to0.
$$
But this is exactly what represents $\ev_i\circ\ev_i'$, by Lemma \ref{OO}.
\end{proof}

So an element of $H^1(\Omega^1_X)\subset\Ext^2(\mc O_{\Delta},\mc
O_{\Delta})=\Hom^0(\mc O_{\Delta},\mc O_{\Delta}[2])$ induces, by Fourier-Mukai transform, a morphism $\F\to\F[2]$ for any $\F\in D^b(X)$. In particular, for every proper curve $D\subset X$ we induce an element of $\Ext^2(\mc O_D,\mc O_D)=H^1(\Omega^1_D)$, which is
described as one might expect.

\begin{lemma} \label{commutes} The diagram
$$ \xymatrix{
H^1(\Omega^1_X) \ar@{^(->}[r]\dto & \Ext^2(\mc O_{\Delta},\mc O_{\Delta}) \dto \\
H^1(\Omega^1_D) \ar@{=}[r] & \Ext^2(\mc O_D,\mc O_D)
}$$
commutes, where the vertical arrows are restriction to $D$ and the Fourier-Mukai action, respectively.
\end{lemma}

\begin{proof}
Working through the
Fourier-Mukai transform this is basically the statement that the induced
map on normal bundles $\nu\_{\Delta}\to\nu\_{D\times X}$ (on the intersection
$\Delta_D$) gives, after contraction with $\theta$, the restriction map
$\Omega^1_X\to\Omega^1_D$ on cotangent bundles.

More explicitly, recall (\ref{d}) that a class $\sigma\in H^1(\Omega^1_X)$ gives rise to an extension $0\to\mc O_X\to\F_\sigma\to T_X\to0$ (\ref{F1}) which can be spliced together with $0\to\Omega^1_{\Delta}\to\mc O_{2\Delta}\to\mc
O_{\Delta}\to0$ (\ref{2delta}) by pushing the former forward under the
diagonal map, and identifying $T_X$ with $\Omega^1_{\Delta}$ using $\theta$.
This gives $[\sigma]\in\Ext^2(\mc O_{\Delta},\mc O_{\Delta})$.

Applying this as a Fourier-Mukai transform to $\mc O_D$ means taking the derived tensor product of it with
$\mc O_{D\times X}$ and pushing down to the second $X$ factor. Since $D\times
X$ is transverse to $\Delta$ we take the usual tensor product, yielding
the splicing together of the exact sequences
\begin{equation} \label{11}
0\to\mc O_D\to\F_\sigma|_D\to T_X|_D\to0
\end{equation}
(pushed forward by $\Delta_D\subset X\times X$), and
\begin{equation} \label{22}
0\to\Omega^1_{\Delta}|\_D\to\mc O_{2\Delta_D\subset D\times X}\to\mc
O_{\Delta_D}\to0.
\end{equation}
The central sheaf above denotes the (pushforward to $X\times X$ of the)
thickening of $\Delta_D$ \emph{inside} $D\times X$.

Now in the first of the above exact sequences (\ref{11}), taking the kernel of the surjections of the last two terms to $\nu\_D$ gives
\begin{equation} \label{33}
0\to\mc O_D\to\F_{\sigma|\_D}\to T_D\to0,
\end{equation}
where $\F_{\sigma|_D}$ is the sheaf on $D$ given by the extension $\sigma|_D
\in H^1(\Omega^1_D)=\Ext^1(T_D,\mc O_D)$. In the second (\ref{22}), this
surjection corresponds, under $\theta$, to $\Omega^1_{\Delta}|_D\to
\Omega^1_{\Delta_D}$. This does \emph{not} lift to a map from the second term of (\ref{22}) to $\Omega^1_{\Delta_D}$, but it does once we pushdown to the
second $X$ factor. That is, the pushdown of (\ref{22}) to $X$ gives
\begin{equation} \label{44}
0\to\Omega^1_X|_D\to J^1_X(\mc O_D)\to\mc O_D\to0,
\end{equation}
where $J^1_X(\mc O_D)$ is the first jet space on $X$ of $\mc O_D$ (which is, by definition, the Fourier-Mukai transform of $\mc O_D$ by $\mc
O_{2\Delta}$). $J^1_X(\mc O_D)$ surjects onto (the pushforward to $X$ of)
$J^1_D(\mc O_D)$ (since $2\Delta$ surjects onto the pushforward to $X\times
X$ of $2\Delta_D\subset D\times D$) and the latter \emph{splits} as $\mc
O_D\oplus\Omega^1_D$ since the Atiyah class on $D$ of $\mc O_D$ vanishes
-- the trivial line bundle on $D$ clearly admits a holomorphic connection. The resulting map $J^1_X(\mc O_D)\to\Omega^1_D$ induces
via (\ref{44}) the standard surjection $\Omega^1_X|\_D\to\Omega^1_D$, and
taking the kernels of these surjections to $\Omega^1_D$ turns (\ref{44})
into
\begin{equation} \label{55}
0\to\nu^*_D\to\mc O_{2D}\to\mc O_D\to0.
\end{equation}
Thus the pushdowns of (\ref{11}) and (\ref{22}) both have commuting surjections
to $\Omega^1_D\cong\nu\_D$, so their spliced representative of
$\Ext^2(\mc O_D,\mc O_D)$ has a map to the trivial complex
$\{\Omega^1_D\rt{\id\,}\Omega^1_D\}$; taking kernels shows that it is quasi-isomorphic to the splicing of (\ref{33}) and (\ref{55}) via the
isomorphism $\theta\colon T_D\cong\nu_D^*$. But this is precisely the sequence
$$
0\to\mc O_D\to\F_{\sigma|\_D}\to\mc O_{2D}\to\mc O_D\to0
$$
defining $[\sigma|_D]\in H^1(\Omega_D^1)=\Ext^2(\mc O_D,\mc O_D)$.
\end{proof}
% We only require this for smooth irreducible curve classes $[C]$ restricted
% to smooth irreducible curves $D$, though for $X$ compact this is of course
% sufficient to prove the whole result by the Lefschetz (1,1) theorem.
% 
% There are two cases: when $C\pitchfork D$ and when $C=D$. For the former.....

So everything is governed by intersection numbers. Let
\begin{equation} \label{int}
C_{ij}\define C_i.C_j=\left\{\begin{array}{ccl} -2 && i=j, \\
1 && |i-j|=1, \\ 0 && |i-j|>1,\end{array}\right.
\end{equation}
be the $(n\times n)$ intersection matrix of the $A_n$-chain of curves $C_i$, with inverse $D_{ij}$. In fact one can calculate that
\begin{equation} \label{Cinverse}
D_{ij}=\frac{-1}{n+1}\min(i,j)\big(n-\max(i,j)+1\big).
\end{equation}
Thus the classes $D_i\define\sum_{i=1}^nD_{ij}[C_j]\in H^1(\Omega^1_X)$ are the
dual basis to the $C_i$ for the span of the $C_i$ in $H^1(\Omega^1_X)$ under the intersection pairing (\ref{int}): $D_i.C_j=\delta_{ij}$.

We are finally ready to define the extra natural transformations $\id\to\id[2]$
that we need to extend the braid group action to a representation of the
braid cobordism category.

\begin{defn}
Let
\begin{equation} \label{Xi}
X_i\colon\mc O_{\Delta}\to\mc O_{\Delta}[2]
\end{equation}
be defined by the above dual basis $D_i\in H^1(\Omega^1_X)\subset\Ext^2(\mc O_{\Delta},\mc O_{\Delta})$. That is, $X_i=\sum_{i=1}^nD_{ij}[C_j]$.

For any $\F\in D^b(X)$ we also let $X_i\colon\F\to\F[2]$ denote the morphism induced by (\ref{Xi}) by Fourier-Mukai transform.
\end{defn}

For a proper irreducible curve $C\subset
X$ let $\theta_C\colon\mc O_C\to\mc O_C[2]$ denote the generator corresponding to $1$ in $\Ext^2(\mc O_C,\mc O_C)\cong H^{1,1}(C)\cong\C$.

\begin{prop} \label{evev}
The morphisms $X_i\colon\mc O_{C_j}\to\mc O_{C_j}[2]$ are $0$ for $i\ne j$
and $\theta_j$ for $j=i$.
There are isomorphisms
\begin{eqnarray*}
\ev_i\circ\ev_i' &\cong& X_{i-1}-2X_i+X_{i+1}\colon\ \mc O_{\Delta}\to
\mc O_{\Delta}, \\
\ev_i'\circ\ev_i &\cong& X_i\boxtimes\id-\id\boxtimes X_i
\colon\ U_i\to U_i[2].
\end{eqnarray*}
\end{prop}

\begin{rmks} Here we have set $X_0=0=X_{n+1}$, and $\id\boxtimes X_i$ (for
instance) maps $\F\boxtimes\mc G\to\F\boxtimes\mc G[2]$ in $D^b(X\times X)$ by the identity on $\F$ and $X_i$ on $\mc G$. Since $U_i=\mc O_{C_i\times C_i}(C_i,0)[-1]$ is of the form $\F\boxtimes\mc G$ this explains the notation above.
\end{rmks}

\begin{proof}
The first statement follows directly from Lemma \ref{commutes} and the
definition of the $D_i$. Since $[C_i]=D_{i-1}-2D_i+D_{i+1}\in H^1(\Omega^1_X)$ by construction, the second statement is just a rephrasing of Proposition \ref{121}.

So it remains to prove the last isomorphism.
Splicing together (\ref{ex1}) and (\ref{ex2}) (in the opposite order from
before) gives the following exact sequence representing the morphism
$\ev_i'\circ\ev_i\colon\mc{O}_{C_i\times C_i}(C_i,0)\to\mc{O}_{C_i\times C_i}(C_i,0)[2]$:
\begin{multline*}
0\to\mc{O}_{C_i\times C_i}(C_i,0)\to \\
\mc{O}_{\Delta\cup(C_i\times C_i)}(L_i)\to\mc{O}_{\Delta\cup(C_i\times C_i)}(C_i,0)\to\mc{O}_{C_i\times C_i}(C_i,0)\to0.
\end{multline*}
$\mc O_{\Delta}(-\Delta_{C_i})$ injects into the central two terms; taking
their cokernels shows that $\{\ev_i'\circ\ev_i\colon U_i\to U_i[2]\}$ is quasi-isomorphic to
\begin{equation} \label{CiCi}
\Big\{\mc{O}_{C_i\times C_i}(C_i,0)(\Delta_{C_i})\to
\mc{O}_{2\Delta_{C_i}\cup\,(C_i\times C_i)}(C_i,0)\Big\}.
\end{equation}
Here $2\Delta_{C_i}$ is the thickening of $\Delta_{C_i}$ inside $C_i\times
C_i$, and the central map above is the restriction map to $\mc O_{\Delta_{C_i}}
(C_i,0)(\Delta_{C_i})\cong\mc O_{\Delta_{C_i}}\cong\I_{\Delta_{C_i}\subset
2\Delta_{C_i}}(C_i,0)$. \medskip

Now $\Ext^2(\mc{O}_{C_i\times C_i}(C_i,0),\mc{O}_{C_i\times C_i}(C_i,0))=
\Ext^2(\mc{O}_{C_i\times C_i},\mc{O}_{C_i\times C_i})$ is isomorphic to
$H^1(\nu\_{C_i\times C_i})\cong H^1(\Omega^1_{C_i\times C_i})\cong
H^1(\Omega^1_{C_i})\oplus H^1(\Omega^1_{C_i})$, with the penultimate
isomorphism given
by contraction with the two form $\theta_1+\theta_2$, where $\theta_j$ is
the pullback of $\theta$ from the $j$th factor of $X\times X$. Under this
isomorphism it is clear that the classes $D_i\oplus0$ and $0\oplus D_i$
correspond to $X_i\boxtimes\id$ and $\id\boxtimes X_i$ respectively.

These classes are, as usual, described the (pushforward
from $C_i\times C_i$ to $X\times X$ of the) extension
class $0\to\mc O_{C_i\times C_i}\to\F\to T_{C_i\times C_i}\to0$ corresponding
to an element of $H^1(\Omega^1_{C_i\times C_i})$, composed with the universal
extension $0\to\nu_{C_i\times C_i}^*\to\mc O_{2(C_i\times C_i)}\to\mc
O_{C_i\times C_i}\to0$ by using $\theta_1+\theta_2$ to identify $T_{C_i\times C_i}\cong\nu_{C_i\times C_i}^*$. That is, these classes $\mc{O}_{C_i\times C_i}\to\mc{O}_{C_i\times C_i}[2]$ are quasi-isomorphic to
\begin{equation} \label{CC}
\{\F\to\mc O_{2(C_i\times C_i)}\}.
\end{equation}

If the $H^1(\Omega^1_{C_i\times C_i})$ class is the fundamental class $[D]$
of an effective divisor $D\subset C_i\times C_i$ then we have seen that we can describe $\F$ as $T_{C_i\times C_i}\times_{\nu\_D}\mc O_{C_i\times C_i}(D)$
(\ref{ft}). Then we divide both terms in (\ref{CC}) by $T_{C_i\times C_i}(-\log D)$ and tensor with $\mc O(C_i,0)$ to show that the morphism $[D]\colon
\mc{O}_{C_i\times C_i}(C_i,0)\to\mc{O}_{C_i\times C_i}(C_i,0)[2]$ is
quasi-isomorphic to the complex
\begin{equation} \label{CiD}
\{\mc O_{C_i\times C_i}(D)(C_i,0)\to\mc O_{(C_i\times C_i)\cup\,2D}(C_i,0)\},
\end{equation}
where $2D$ is the first order thickening of $D$ in the directions perpendicular
to $T_D$ under $\theta_1+\theta_2$ (which contain one more direction other than those of $T_{C_i\times C_i}$).

Picking $D=\Delta_{C_i}\subset C_i\times C_i$ we get $\Delta_{C_i}$ thickened
inside $\Delta$ in (\ref{CiD}), which is not quite (\ref{CiCi}).
\emph{But} if we change the identification $T_{C_i\times C_i}\cong\nu_{C_i\times C_i}^*$ by using $\theta_1-\theta_2$ instead of $\theta_1+\theta_2$, then
we get $D$ thickened in the directions perpendicular
to $T_D$ under $\theta_1-\theta_2$ in (\ref{CiD}). This now describes the
$\Ext^2$ class corresponding not to $[D]\in H^1(\Omega^1_{C_i\times C_i})$
but to its image under the map multiplying its second component in its K\"unneth
decomposition by $-1$.

Applying this to $D=\Delta_{C_i}$, whose class in
$H^1(\Omega^1_{C_i\times C_i})$ has degree 1 on both factors, proves that
the $\Ext^2$ class described by the cohomology class $(D_i,-D_i)$ corresponds
to the complex (\ref{CiD}) with $D=\Delta_{C_i}$ thickened inside $C_i\times
C_i$. But this is precisely (\ref{CiCi}), so that $(X_i\boxtimes\id)-
(\id\boxtimes X_i)\cong\ev_i'\circ\ev_i$.
\end{proof}

Given a kernel $K\in D^b(X\times X)$ we denote by $l_{X_i}\colon K\to K[2]$ the following morphism (giving rise to a natural transformation between the
Fourier-Mukai functors associated to $K$ and $K[2]$):
\begin{equation} \label{li}\xymatrix{
K\cong\mc O_{\Delta}\star K\rto^(.43){X_i\star\id} &
\mc O_{\Delta}[2]\star K\cong K[2]}.
\end{equation}
Similarly $r_{X_i}\colon K\to K[2]$ is
\begin{equation} \label{ri}\xymatrix{
K\cong K\star\mc O_{\Delta}\rto^(.43){\id\star X_i} &
K\star\mc O_{\Delta}[2]\cong K[2]}.
\end{equation}
It is clear that on $\mc O_{\Delta}$, $l_{X_i}$ and $r_{X_i}$ act as the same morphism.
Given a linear combination $a=\sum_ja_jX_j$ of the $X_i$ we also set
$$
l_a\define\sum_ja_jl_{X_j} \quad\mathrm{and}\quad r_a\define\sum_ja_jr_{X_j}.
$$
Then we have the following analogue of Proposition \ref{lara}.

\begin{prop} \label{liri}
The morphisms $T_i\to T_i[2]$ and $T_i'\to T_i'[2]$ given by $l_a-
r_{a+a_i(X_{i-1}-2X_i+X_{i+1})}$ are both zero. \\
For instance for $i\ne j$, $l_{X_j}$ and $r_{X_j}$ act as isomorphic morphisms on $T_i$ and as zero on $U_i$.
\end{prop}

\begin{proof}
We first claim that on $U_i=\mc O_{C_i}(C_i)[-1]\boxtimes\mc O_{C_i}$, the
functors $l_{X_i}$ and $r_{X_i}$ act as $\id\boxtimes X_i$ and $-X_i\boxtimes\id$ respectively. For instance $r_{X_i}$ acts as
$$\xymatrix{
U_i\star\mc O_{\Delta}\rto^(.47){\id\star X_i} & U_i\star\mc O_{\Delta}[2]},
$$
which is, in the usual notation of Section \ref{FMTs} (cf. (\ref{conv})),
$$\xymatrix{
\pi_{13*}\Big[\mc O_{\Delta\times X}\otimes\mc O_{X\times C_i\times C_i}
(0,C_i,0)[-1]\rrto^(.52){\pi_{12}^*X_i\,\otimes\,\id} && \mc O_{\Delta\times X}\otimes\mc O_{X\times C_i\times C_i}(0,C_i,0)[1]\Big]},
$$
which is
$$
\pi_{13*}\Big[\mc O_{\Delta_{C_i}\times C_i}
(0,C_i,0)[-1]\To\mc O_{\Delta_{C_i}\times C_i}(0,C_i,0)[1]\Big].
$$
Here the arrow is pulled back from the first and second factors, so the result
is of the form
$$\xymatrix{
\mc O_{C_i\times C_i}(C_i,0)[-1]\rto^{Z\boxtimes\id} &
\mc O_{C_i\times C_i}(C_i,0)[1]},
$$
where, by the above working, $Z$ is
\begin{equation} \label{Zed} \xymatrix{
\pi_{1*}\Big[\mc O_{\Delta}\otimes\mc O_{X\times C_i}
(0,C_i)[-1]\rto^(.53){X_i\otimes\id} & \mc O_{\Delta}\otimes\mc O_{X\times C_i}(0,C_i)[1]\Big]}.
\end{equation}
This is where the sign comes in, as this is \emph{not} $X_i\colon
\mc O_{C_i}(C_i)[-1]\to\mc O_{C_i}(C_i)[1]$, the latter being by definition the map of Fourier-Mukai transforms
$$
\Phi_{\mc O_{\Delta}}(\mc O_{C_i}(C_i)[-1])\Rt{\Phi_{X_i}\,}
\Phi_{\mc O_{\Delta}[2]}(\mc O_{C_i}(C_i)[-1]),
$$
which is
$$\xymatrix{
\pi_{2*}\Big[\mc O_{\Delta}\otimes\mc O_{C_i\times X}
(C_i,0)[-1]\rto^(.53){X_i\otimes\id} & \mc O_{\Delta}\otimes\mc O_{C_i\times X}(C_i,0)[1]\Big]}.
$$
Comparing with (\ref{Zed}), the two differ just by swapping the two factors
of $X\times X$, which sends $X_i$ to $-X_i$ because it reverses the identification
of the normal bundle of $\Delta$ with $T_{\Delta}$ (\ref{v-v}).

For $l_{X_i}$ the working is the same but easier, without the swapping of the
$X\times X$ factors.

Now $l_a$ acts on $T_i$ as
$$\xymatrix{
\mathrm{Cone}\Big(\mc O_{C_i\times C_i}(C_i,0)[-1]
\rrto^(.7){\ev_i}\dto_{l_a} && \mc O_{\Delta}\Big)\hspace{-5mm}
\dto^{l_a} \\ \mathrm{Cone}\Big(\mc O_{C_i\times C_i}(C_i,0)[1]
\rrto^(.7){\ev_i} && \mc O_{\Delta}[2]\Big),\hspace{-1cm}}
$$
which, by the above and the first part of Proposition \ref{evev} is
$$\xymatrix{
\mathrm{Cone}\Big(\mc O_{C_i\times C_i}(C_i,0)[-1]
\rrto^(.7){\ev_i}\dto_{a_i\id\boxtimes X_i} &&
\mc O_{\Delta}\Big)\hspace{-5mm} \dto^{a=\sum_ja_jX_j} \\
\mathrm{Cone}\Big(\mc O_{C_i\times C_i}(C_i,0)[1] \rrto^(.7){\ev_i} && \mc O_{\Delta}[2]\Big).\hspace{-1cm}}
$$
Now use the homotopy $h$ which maps the top right corner to the bottom left
by $\ev_i'$, so that $dh+hd$ is $\ev_i'\circ\ev_i$ on the left and $\ev_i
\circ\ev_i'$ on the right. Adding $a_i$ times this to the given maps and using Proposition \ref{evev}, we find that $l_a\colon T_i\to T_i$ is
homotopic to
$$\xymatrix{
\mathrm{Cone}\Big(\mc O_{C_i\times C_i}(C_i,0)[-1] \rrto^(.65){\ev_i}
\dto_{a_i\id\boxtimes X_i+a_i(X_i\boxtimes\id-\id\boxtimes X_i)}
^{=\,a_iX_i\boxtimes\id} && \mc O_{\Delta}\Big)\hspace{-5mm}
\dto^{a+a_i(X_{i-1}-2X_i+X_{i+1})} \\
\mathrm{Cone}\Big(\mc O_{C_i\times C_i}(C_i,0)[1] \rrto^(.65){\ev_i} && \mc O_{\Delta}[2]\Big).\hspace{-1cm}}
$$
Since the coefficient of $X_i$ in $a+a_i(X_{i-1}-2X_i+X_{i+1})$ is $-a_i$,
this is $r_{a+a_i(X_{i-1}-2X_i+X_{i+1})}\colon T_i\to T_i$.

The proof for $T_i'$ is almost identical, using $\ev_i$ in the homotopy.
\end{proof}

We define the morphisms (between Fourier-Mukai kernels, inducing natural
transformations between the corresponding functors) associated to positive
and negative braid moves as follows.

For a \emph{positive braid move} we define the morphisms
\begin{equation} \label{moves}
\mc O_\Delta\to T_i \qquad\text{and}\qquad T_i'\to\mc O_\Delta
\end{equation}
by the first morphism of (\ref{ex1}) and the second morphism of (\ref{ex2})
respectively.

For the \emph{negative braid moves} we write $T_i=\mc{O}_{\Delta\cup(C_i\times C_i)}(C_i,0)$ as the cone on $\mc O_{C_i\times C_i}(C_i,0)[-1]\Rt{\ev_i}\mc O_\Delta$ by (\ref{ex1}), and then define $\delta_i\colon T_i\to\mc O_\Delta[2]$ by
\begin{equation} \label{-move} \xymatrix{
\mathrm{Cone}\Big(\mc O_{C_i\times C_i}(C_i,0)[-1] \rto^(.75){\ev_i} & \mc O_\Delta\Big)\!\!\!\!\!\! \dto^{X_{i+1}-X_{i-1}} \\ & \mc O_\Delta[2].
\hspace{-8mm}}
\end{equation}
To check this really defines a morphism in $D^b(X\times X)$ we need to know
that $(X_{i-1}-X_{i+1})\circ\ev_i\cong0$. In fact we claim that
\begin{equation} \label{nohom}
\Hom^0(\mc O_{C_i\times C_i}(C_i,0)[-1],\mc O_\Delta[2])=0.
\end{equation}
To prove (\ref{nohom}), note that $C_i\times C_i$ and $\Delta$ intersect with excess dimension one and excess normal bundle $\mc O_{\Delta_{C_i}}(C_i)$.
Therefore a standard Koszul resolution argument shows that the only nonzero $\mc Ext^{\,i\,}$ is
\begin{equation} \label{singlehom}
\mc Ext^1(\mc O_{C_i\times C_i}(C_i,0),\mc O_\Delta)\cong\mc O_{\Delta_{C_i}},
\end{equation}
generated by $\ev_i$. So the local-to-global spectral sequence
$H^i(\mc Ext^j)\Longrightarrow\Ext^{i+j}$ shows that $\Ext^3(\mc O_{C_i\times C_i}(C_i,0),\mc O_\Delta)=0$.

Similarly, writing $T_i'$ as the cone on $\mc O_\Delta[-1]\Rt{\ev_i'}
\mc O_{C_i\times C_i}(C_i,0)$ by (\ref{ex2}), then a chase around the diagram
proving that $T_i'\star T_i\cong\mc O_\Delta$ (i.e. a proof of braid movie
move 12, which will be given below) shows that the corresponding morphism
$\delta_i\colon\mc O_\Delta[-2]\to T_i'$ is
\begin{equation} \label{--move} \xymatrix{
\qquad\ \mc O_\Delta[-3] \dto_{X_{i+1}-X_{i-1}} \\
\!\!\!\mathrm{Cone}\Big(\mc O_\Delta[-1] \rto^(.45){\ev_i'} &
\mc O_{C_i\times C_i}(C_i,0)\Big).}
\end{equation}
% Don't need the following as invertible elemnets of $H^0(\mc O_X)$ are scalars
% for $X$ quasi-projective ?
% As a way to simplify proofs we restrict attention to a certain class of
% surfaces.
% 
% \begin{defn} We call the surface $X$ \emph{good} if it is either projective
% or the affine surface that is the minimal resolution of the $A_n$ singularity
% $\{xy=z^{n+1}\}\subset\C^3$.
% \end{defn}

\begin{theorem}
% For good $X$
The $T_i,\,T_i'$ and morphisms described above induce a projective action of the braid cobordism category on $D^b(X)$.
\end{theorem}

\begin{proof}
The braid relations were already proved in \cite{ST}. For the first 10 braid
movie moves of Section \ref{sec-bc1} we use the same trick as in the proof
of Theorem \ref{first-theorem} in Section \ref{An}. They each
describe two isomorphisms between the same two Fourier-Mukai kernels (a combination
of $T_i$s and $T_j'$s), and we would like to show that these isomorphisms
are projectively equivalent. Composing one with the inverse of the other
(via the operation $\star$ of (\ref{conv})),
we get an isomorphism from one kernel to itself which we are required to
show is a scalar. But the kernel is invertible, so composing with its inverse we get an isomorphism $\mc O_\Delta\to
\mc O_\Delta$, i.e. an invertible element of $\Hom^0(\mc O_\Delta,\mc O_\Delta)
\cong H^0(\mc O_X)$. This must be a non-zero scalar by the Nullstellensatz.

For the moves 11--15 with one \emph{positive} branch point, we use the fact
that $\Hom^0(T_i',\mc O_\Delta)\cong\C\cong\Hom^0(\mc O_\Delta,T_i)$ with generators the second morphism of (\ref{ex2}) and the first morphism of
(\ref{ex1}), respectively. This is proved in much the same way as
one computes the other degree parts of (\ref{nohom}) (or can be deduced from it). Namely (\ref{singlehom}) shows that the morphisms from $U_i$ to $\mc
O_\Delta$ are one-dimensional, generated by $\ev_i$. A long exact sequence
then gives the result we want.

By composing with invertible
Fourier-Mukai kernels we find that the morphisms from a kernel at the
top of any of the movies 11--14 to the one at the bottom are also one dimensional
and generated by the morphism represented by the movie. Therefore under a
braid movie move this can only change by an invertible scalar.

Next we deal with move 11 with one \emph{negative} branch point using the
map $\delta_j$ (\ref{--move}) by showing
that, for instance, for $|i-j|>1$ the composition
\begin{equation} \label{sunday}
T_i\cong T_i\star\mc O_\Delta\Rt{\delta_j}T_i\star T_j'[2]\cong T_j'\star T_i[2]
\end{equation}
is the same as the composition
\begin{equation} \label{monday}
T_i\cong\mc O_\Delta\star T_i\Rt{\delta_j}T_j'\star T_i[2].
\end{equation}
But the first composition (\ref{sunday}) is the map $r_{j+1}-r_{j-1}$ of the cone $U_i\Rt{\ev_i}\mc O_\Delta$ into the top line of the cone
$$ \xymatrix{
U_i[2] \rto^{\ev_i}\dto^{\ev_j'} & \mc O_\Delta[2] \dto^{\ev_j'} \\
U_i\star U_j[4] \rto^(.57){\ev_i} & U_j[4].}
$$
Since this map is zero on the $U_i$ piece (by the last part of Proposition
\ref{liri}) and equal to $l_{j-1}-l_{j+1}$ on $\mc O_\Delta$, and since
$U_i\star U_j\cong0\cong U_j\star U_i$ for $|i-j|>1$, we find it is isomorphic to the map $l_{j-1}-l_{j+1}$ of the cone $U_i\Rt{\ev_i}\mc O_\Delta$ into the top line of the cone
$$ \xymatrix{
U_i[2] \rto^{\ev_i}\dto^{\ev_j'} & \mc O_\Delta[2] \dto^{\ev_j'} \\
U_j\star U_i[4] \rto^(.57){\ev_i} & U_j[4].}
$$
Reflecting this in its diagonal gives the second composition (\ref{monday}).

For move 12 with one negative branch point we recall the quasi-isomorphism
that makes $T_i\star T_i'\cong\id$. (This is given in \cite{ST} using functors,
but that proof can be easily translated into Fourier-Mukai transforms, and
takes the exact same form as our proof in Proposition \ref{invert} since
the $U_i$s satisfy the same relations.) Namely $U_i\star U_i\cong U_i[-2]\oplus
U_i$, and $T_i\star T_i'$ is the total cone on
$$ \spreaddiagramrows{-15pt} \xymatrix{
\,U_i \rto^(.35){\ev_i'}\ddto^{\ev_i} & U_i\oplus U_i[2] \ddto^{\ev_i} \\\\
\ \mc O_\Delta \rto^(.45){\ev_i'} & U_i[2]}
\xymatrix{\\\\ \qquad\Longleftarrow\qquad}
\xymatrix{& U_i[2] \ar@{=}[dd] \\\\ \mc O_\Delta \rto^(.45){\ev_i'} & U_i[2]}
\xymatrix{\\\\ \qquad\Longrightarrow\qquad}
\spreaddiagramrows{21pt} \xymatrix{\\ \mc O_\Delta}
$$
with the double arrows the obvious quasi-isomorphisms. Then move 12 says
that the morphism from the first square to the cone
$\big(\mc O_\Delta[2]\Rt{\ev_i'}U_i[4]\big)$ given by zero on the top row and $r_{X_{i-1}-X_{i+1}}$ applied to both terms on the bottom row should be isomorphic to the morphism from the last square (i.e. just $\mc O_\Delta$) to the cone $\big(\mc O_\Delta[2]\Rt{\ev_i'}U_i[4]\big)$ given by $r_{X_{i-1}-X_{i+1}}$ from $\mc O_\Delta$ to $\mc O_\Delta[2]$. But indeed both are isomorphic to the morphism from the central square to the cone
$\big(\mc O_\Delta[2]\Rt{\ev_i'}U_i[4]\big)$ given by zero on the top row and $X_{i-1}-X_{i+1}$ applied to both terms on the bottom row.

Move 13 we defer to the proofs for $D^b(T^*Fl)$ in the next section; this
paper is quite long enough without repeating proofs twice, so we concentrate
on the harder $T^*Fl$ case -- for surfaces the same proof applies but is
much easier since there is no $U_{ij}$ term involved in proving Theorem \ref{YangB}.
As we noted before, move 14 follows from the various versions of move 13.

Finally locality moves 15 with one negative branch point follow from the
second part of Proposition \ref{liri}; we describe how in the example of
move 15 drawn in figure 7. For that we combine the morphism
$\delta_k\colon T_k\to\mc O_\Delta[2]$ of (\ref{-move}) and a standard isomorphism $\phi$ of Fourier-Mukai kernels $K\rt{\phi}K'$ (a composition of one of the braid relations)
consisting of cones of morphisms of kernels with
expressions $U_i\star\ldots\star U_j$ involving only $U_i$s with $|i-k|>1$
(and one copy of $\mc O_\Delta$, on which $\phi$ acts as the identity). We would like to show that the compositions of the following vertical maps
of horizontal cones
$$ \xymatrix{
\Big(K\star U_k \rto^(.6){\ev_k}\dto^\phi & K\dto^\phi \Big)\!\!\!\! &&
\text{and} && \Big(K\star U_k \rto^(.6){\ev_k} & K\dto^{r_{k-1}-r_{k+1}} \Big)\!\!\!\! \\
\Big(K'\star U_k \rto^(.6){\ev_k}\dto & K'\dto^{r_{k-1}-r_{k+1}} \Big)\!\!\!\! &&&&& K[2] \dto^\phi \\ \!\!\!\!\Big(0 \rto & K'[2] \Big)\!\!\!\! &&&&& K'[2]}
$$
are the same, which is the case if the following compositions are equal:
$$ \xymatrix{
K \rto^\phi & K' \rrto^(.45){r_{k-1}-r_{k+1}} && K'[2] & \text{and} &
K \rrto^(.45){r_{k-1}-r_{k+1}} && K[2]\rto^{\phi} & K'[2]}.
$$
But by Proposition \ref{liri}, $r_{k-1}-r_{k+1}$ acts as zero on every term in $K$ and $K'$ except for the $\mc O_\Delta$ term in each, on which $\phi$ acts as the identity and so commutes with $r_{k-1}-r_{k+1}$.
\end{proof}

One can show that the central ambiguity in the projective action of
the braid cobordism category is at most $\pm1\subset \C^*$,  by passing
from $D^b(X)$ to the homotopy category $\mc C(A_n\otimes\C)$ of Section
\ref{sec-survey}, as described at the beginning of this Section. Under the
equivalence \cite{KS, ST} of the derived categories of differential graded modules over $\RHom(\oplus_i\mc O_{C_i},\oplus_i\mc O_{C_i})$ and
$A_n\otimes\C$, the $\mc O_{C_i}$
are mapped to the $P_i$, $T_i$ to $R_i$, $T_i'$ to $R_i'$, and $\ev_i,\,\ev_i'$
to $\beta_i,\,\gamma_i$. The rational numbers $D_{ij}$ (\ref{Cinverse}) that
we used become integers under this equivalence, since by Proposition \ref{evev} the $X_i$ act on the $\mc O_{C_i}$ by integral generators of $H^{1,1}(C_j)
\cong\Ext^2(\mc O_{C_j},\mc O_{C_j})$.

The subcategory generated by the projective modules
$P_i$ is then equivalent to $\mc C(A_n\otimes\C)$, and the above positive
(\ref{moves}) and negative (\ref{-move}) braid moves map to those of
(\ref{positivemove1}, \ref{positivemove2}) and (\ref{negativemove1},
\ref{negativemove2}) respectively: for the positive braid moves this
is obvious, while for the negative ones it follows from Proposition \ref{evev}.
Therefore the $\C^*$ scalar ambiguities in the movie moves in $D^b(X)$ map to the scalar ambiguities $\pm1$ in $\mc C(A_n\otimes\C)$, and so are themselves
$\pm1$. Of course we expect that they can all be shown to be $+1$, but this
would require more work.

\subsection{The cotangent bundle of the flag variety}
\label{flagcobordism}

The construction of the last section goes over almost unchanged to $D^b(T^*Fl)$; we outline the steps here.
The kernels $U_i,\ T_i,\ T_i'$ from Sections \ref{flagsection} and \ref{surfaces}
are obviously analogous, as are the maps between them. Using the holomorphic
symplectic structure of $K=T^*Fl$ we find that $T_K\cong\Omega^1_K$, so that the description (\ref{GS}) of $\Ext^2(\mc O_\Delta,\mc O_\Delta)$ still holds. 
The divisors $N_i\subset K$ define de Rham classes $[N_i]\in H^1(\Omega^1_K)
\subset\Ext^2(\mc O_\Delta,\mc O_\Delta)$ by the same formulae as for the $[C_i]$ before (splicing the logarithmic one-forms sequence (\ref{dR}) to the doubled diagonal exact sequence (\ref{doublediag}) via $T_K\cong\Omega^1_K$).

Taking the same linear combinations $D_{ij}$ (\ref{Cinverse}) of these classes as before we define the $X_i\in\Ext^2(\mc O_\Delta,\mc O_\Delta)$ as $\sum_j
D_{ij}[N_j]$.

The analogue of Lemma \ref{commutes} is that the action of such classes on some $\mc O_{N_k}$ by Fourier-Mukai transform is described by the restriction of the $H^1(\Omega^1_K)$ class to $H^1(\Omega^1_{N_k/K_k})\cong
H^1(\mc O_{N_k}(N_k))\cong\Ext^2(\mc O_{N_k},\mc O_{N_k})$. This is proved
in exactly the same way by noting that the natural map $T_K\to\mc O_{N_k}(N_k)$
becomes, on contraction with the holomorphic symplectic form, the restriction
map $\Omega^1_K\to\Omega^1_{N_k/K_k}$. (The fact that $\Ext^2(\mc O_{N_k},\mc O_{N_k})\cong H^1(\mc O_{N_k}(N_k))$ follows easily from the local-to-global spectral sequence for Ext and the Leray spectral sequence for $N_k\to K_k$.)

Now in turn, we can calculate $H^1(\Omega^1_{N_k/K_k})$ by the Leray spectral sequence, using the fact that $R^1p_{k*}\Omega^1_{N_k/K_k}\cong\mc O_{K_k}$ by relative Serre duality, with all other derived functors zero. Therefore
basechange holds, and the fibers of $p_k\colon N_k\to K_k$ are compact, so to determine a class in $H^1(\Omega^1_{N_k/K_k})\cong H^0(\mc O_{K_k})$ we
need only determine its de Rham or singular cohomology class on each fiber,
for instance by intersection theory.

In particular, for $|i-j|>1$, $N_i\cap N_j=N_{ij}$ is zero on the generic
fiber of $N_j\to K_j$ (being the pullback of the divisor $K_{ij}\subset K_j$)
so $[N_i]|_{N_j}$ defines the zero class in $H^0(\mc O_{K_j})$. (So here it is crucial that
the $\Ext^2(\mc O_{N_j},\mc O_{N_j})$ class depends only on the \emph{relative},
fiberwise, $H^2$ class of $[N_i]|_{N_j}$, rather its absolute class, which
is nonzero.) For $i-j=1$, $N_i\cap N_j
=N_{ij}$ is generically a section of $N_j\to K_j$ (away from $M_{ij}\subset N_{ij}$, by Lemma \ref{blowup}) and so defines the canonical generator $1\in H^0(\mc O_{K_j})$ away from $K_{ij}\subset K_j$, and so in fact everywhere by Hartog's
theorem. Similarly, considering $i=j$, $N_i$ has degree $-2$ on each fiber of $N_i\to K_i$ so we find that $[N_i]|_{N_i}$ is represented,
via the Leray spectral sequence, by $-2\in H^0(\mc O_{K_i})$.

These numbers coincide, of course, with the intersections of the $C_i$ in
Section \ref{surfaces}. Therefore, by construction of the $D_{ij}$,
(\ref{Cinverse}) we find that the
$X_i$ act as zero on $\mc O_{N_j}$ for $j\ne i$ and as the generator of
$\Ext^2(\mc O_{N_i},\mc O_{N_i})\cong H^0(\mc O_{K_i})$ on $\mc O_{N_i}$.

This gives the correct analogue of the first part of Proposition \ref{evev}.
The last part follows almost verbatim the proof for surfaces. For the second
part the proof is also very similar, using relative de Rham or singular cohomology
classes as above instead of absolute ones. It follows immediately that the
corresponding $l_{X_i}$ and $r_{X_i}$ satisfy the same result as in Proposition \ref{liri}.

With all of these technical results done, the proofs of the braid movie moves
carry over just as before. For move 11 with one negative branch point we
no longer have $U_j\star U_i\cong0\cong U_i\star U_j$ for $|i-j|>1$, but
all we need is that they commute. But the proof of Proposition \ref{commute} shows that both $U_i\star U_j$ and $U_j\star U_i$ are isomorphic to $\mc O_{N_{ij}\times_{T^*Fl_{ij}}N_{ij}}(D_{ij})[-2]$ in the notation of (\ref{Tij}).

We have to prove the following version of move 13 with one negative branch point; that for $j=i+1$ the following two compositions
$$ \diagram
(U_i \rto^{\ev_i} & \id)\star T_j\star T_i\hspace{-18mm} \dto_{\delta_i} &&&& T_j\star T_i\star(U_j \rto^{\ev_j} & \id)\!\!\! \dto_{\delta_j} \\
& \id\star T_j\star T_i[2] \hspace{-19mm} &&&&&
T_j\star T_i\star\id[2] \hspace{11mm}
\enddiagram
$$
become isomorphic maps $T_i\star T_j\star T_i\to T_j\star T_i[2]$ upon
identifying $T_i\star T_j\star T_i\cong T_j\star T_i\star T_j$ by the second part of Theorem \ref{YangB}.

By the first part of Theorem \ref{YangB} the isomorphism $T_i\star T_j\star T_i\cong T_j\star T_i\star T_j$ comes from an isomorphism $U_i\star T_j\star T_i\cong T_j\star T_i\star U_j$ of the first two terms in the two cones
on the top line of the above diagram; as proved there this intertwines the two $\ev$ maps on the top line. Since these two terms map to zero in the
bottom line, it remains to prove that the compositions
$$ \diagram
\id\star T_j\star T_i\hspace{-15mm} \dto_{\delta_i} &&&& T_j\star T_i\star\id \hspace{15mm} \dto_{\delta_j} \\
\id[2]\star T_j\star T_i \hspace{-2cm} &&&& T_j\star T_i\star\id[2]
\hspace{11mm} \enddiagram
$$
are equal. The left hand one is $(l_{X_{i+1}}-l_{X_{i-1}})\star\id=
l_{X_{i+1}}\star\id$
by the last part of Proposition \ref{liri}. By the first part of
Proposition \ref{liri} this is $(r_{X_i}-r_{X_{i+1}})\star\id=
\id\star(l_{X_i}-l_{X_{i+1}})$.

The right hand map is $\id\star(r_{X_{j+1}}-r_{X_{j-1}})=\id\star(-r_{X_{j-1}})$
by the last part of Proposition \ref{liri}. By the first part of
Proposition \ref{liri} this is also $\id\star(l_{X_i}-l_{X_{i+1}})$.

The other negative versions of move 13 are similar; for instance showing
that the following diagram becomes commutative
$$ \diagram
\id[-2]\star T_j'\star T_i' \hspace{-23mm} \dto_{\delta_i} &&&&& \hspace{-8mm}
T_j'\star T_i'\star\id[-2] \dto_{\delta_j} \\
\!(\id \rto^(.3){\ev_i'} & U_i)\star T_j'\star T_i'  &&&& \hspace{-17mm} T_j'\star T_i'\star(\id \rto^{\ev_j'} & U_j)\!\!\!
\enddiagram
$$
on mapping each $T_j'\star T_i'$ on the left hand side by the identity to
the corresponding $T_j'\star T_i'$ on the right hand side, and mapping
$U_i\star T_j'\star T_i'$ to $T_j'\star T_i'\star U_j$ by the isomorphism
given by taking adjoints in Proposition \ref{YangB}. This commutes with maps
on both sides by the adjoint of Proposition \ref{YangB}, so we are left with
showing that the maps
$$
(l_{X_{i+1}}-l_{X_{i-1}})\star\id \quad\text{and}\quad \id\star(r_{X_{j+1}}
-r_{X_{j-1}})\colon\ T_j'\star T_i'[-2]\to T_j'\star T_i'
$$
are the same. But this follows as above by using Proposition \ref{liri}.

% $H^1(\Omega^1_{N_j/K_j})\cong H^0(\mc O_{K_j})\cong H^0(\mathrm{Sym}^*T_{Fl_j})$,
% with the grading on the total symmetric product being the one induced by
% the scaling $\C^*$-action on the fibers of $T^*Fl$. 

\noindent {\small\tt khovanov@math.columbia.edu} \\
\noindent Department of Mathematics, Columbia University, New York, NY 10027.
\smallskip \\
\noindent {\small\tt richard.thomas@imperial.ac.uk} \newline
Department of Mathematics, Imperial College, London SW7 2AZ.

\end{document}